\documentclass[12pt]{article}

\usepackage{amsmath,amsthm,amsfonts,amssymb, color}

\usepackage{graphicx}

\usepackage{subcaption}

\usepackage{tikz} 
\usetikzlibrary{shapes,arrows} 

\newtheorem{theorem}{Theorem}[section]
\newtheorem{corollary}[theorem]{Corollary}
\newtheorem{example}[theorem]{Example}
\newtheorem{proposition}[theorem]{Proposition}
\newtheorem{lemma}[theorem]{Lemma}
\newtheorem{definition}[theorem]{Definition}
\newtheorem{remark}[theorem]{Remark}

\def\cA{\mathcal{A}}
\def\cB{\mathcal{B}}
\def\cD{\mathcal{D}}
\def\cE{\mathcal{E}}
\def\cF{\mathcal{F}}

\def\cM{\mathcal{M}}

\def\cU{\mathcal{U}}

\def\bC{\mathbb{C}}
\def\bD{\mathbb{D}}
\def\bR{\mathbb{R}}
\def\bS{\mathbb{S}}
\def\bZ{\mathbb{Z}}

\def\e{\varepsilon}
\def\SD{\mathbb{S}_{\mathbb{D}}}

\def\barD{\overline{\mathbb{D}}_0}
\def\oZ{\overline{Z}}

\begin{document}

\title{Stable L\'evy motion with values in the Skorokhod space: construction and approximation}

\author{Raluca M. Balan\footnote{Corresponding author. Department of Mathematics and Statistics, University of Ottawa,
585 King Edward Avenue, Ottawa, ON, K1N 6N5, Canada. E-mail
address: rbalan@uottawa.ca} \footnote{Research supported by a
grant from the Natural Sciences and Engineering Research Council
of Canada.}\and
Becem Saidani\footnote{Department of Mathematics and Statistics, University of Ottawa,
585 King Edward Avenue, Ottawa, ON, K1N 6N5, Canada. E-mail address: bsaid053@uottawa.ca}}

\date{September 5, 2018}
\maketitle

\begin{abstract}
\noindent In this article, we introduce an infinite-dimensional analogue of the $\alpha$-stable L\'evy motion, defined as a L\'evy process $Z=\{Z(t)\}_{t \geq 0}$ with values in the space $\mathbb{D}$ of c\`adl\`ag functions on $[0,1]$, equipped with Skorokhod's $J_1$ topology. For each $t \geq 0$, $Z(t)$ is an $\alpha$-stable process with sample paths in $\bD$, denoted by $\{Z(t,s)\}_{s\in [0,1]}$. Intuitively, $Z(t,s)$ gives the value of the process $Z$ at time $t$ and location $s$ in space. This process is closely related to the concept of regular variation for random elements in $\bD$ introduced in \cite{dehaan-lin01} and \cite{hult-lindskog05}.
We give a construction of $Z$ based on a Poisson random measure, and we show that $Z$ has a modification whose sample paths are c\`adl\`ag functions on $[0,\infty)$ with values in $\bD$. Finally, we prove a functional limit theorem which identifies the distribution of this modification as the limit of the partial sum sequence $\{S_n(t)=\sum_{i=1}^{[nt]}X_i\}_{t\geq 0}$, suitably normalized and centered, associated to a sequence $(X_i)_{i\geq 1}$ of i.i.d. regularly varying elements in $\bD$.
\end{abstract}



\vspace{1mm}


\tableofcontents

\section{Introduction}

Regularly varying random variables play an important role in probability theory, being used as models for heavy-tailed observations (observations which may assume extreme values with high probability). In many applications, one is often interested in the sum of such variables. For instance, if $X_i$ denotes the number of internet transactions performed on a secure website on day $i$, it might be of interest to study the total number $\sum_{i=1}^{n}X_i$ of transactions performed on this website in $n$ days. If $(X_i)_{i\geq 1}$ are independent and identically distributed (i.i.d.) regularly varying random variables, then, with suitable normalization and centering, the partial sum process $\{S_n(t)=\sum_{i=1}^{[nt]}X_i\}_{t \geq 0}$ converges as $n \to \infty$ to the $\alpha$-stable L\'evy motion, a process which plays the same central role for heavy-tailed observations as the Brownian motion for observations with finite variance.

With the rapid advancement of technology, data is no longer observed at fixed moments of time, but continuously over a fixed interval in time or space (which we may identify with the interval $[0,1]$). If this measurement is expected to exhibit a sudden drop or increase over this fixed interval, then an appropriate model for it could be a random element in an infinite dimensional space, such as the Skorokhod space $\bD=\bD([0,1])$ of c\`adl\`ag functions on $[0,1]$ (i.e. right-continuous functions with left limits). For instance, if the number of internet transactions is observed continuously during the 24-hour duration of the day (identified with the interval $[0,1]$) and $X_i(s)$ is the number recorded at time $s$ of day $i$, then we may assume that $X_i=\{X_i(s)\}_{s \in [0,1]}$ is a process with c\`adl\`ag sample paths. Another example is when $X_i(s)$ represents the energy produced by a wind turbine on day $i$ at location $s$ of a large wind farm situated on the ocean shore, modeled by the interval $[0,1]$.
In these examples, we are interested in studying the behaviour of the partial sum process $\{\sum_{i=1}^{n}X_i(s); s\in [0,1]\}$ which gives the full information about the total number of transactions (or the total amount of energy) for $n$ days, at each time $s$ during the 24-hour period (or at each location $s$ on the shore).

The goal of this article is to study the macroscopic limit (as time gets large) of the partial sum sequences as those appearing in the previous examples, associated to i.i.d. regularly varying elements in $\bD$. It turns out that this limit is an interesting object in itself, which deserves special attention and will be call an {\em $\bD$-valued $\alpha$-stable L\'evy motion} by analogy with its $\bR^d$-valued counterpart.

Our methods were deeply inspired by Resnick's beautiful presentation of the construction of the classical $\alpha$-stable L\'evy motion with values in $\bR^d$, and of its approximation by partial sums of i.i.d. regularly varying vectors, given in \cite{resnick07}.
Its aim is to extend these results to the infinite-dimensional setting, using the
concept of regular variation for random elements in $\bD$ introduced in \cite{dehaan-lin01}, and developed further in \cite{hult-lindskog05}.
More precisely, our goals are: {\em (i)} to construct a L\'evy process $\{Z(t)\}_{t \geq 0}$ with values in $\bD$, whose marginal $Z(t)=\{Z(t,s)\}_{s \in [0,1]}$ is a c\`adl\`ag $\alpha$-stable process (with a specified distribution); {\em (ii)} to show that this process has a modification whose sample paths are c\`adl\`ag functions from $[0,\infty)$ to $\bD$ (where $\bD$ is endowed with Skorohod $J_1$-topology); and {\em (iii)} to identify this modification as the limit as $n \to \infty$ of the partial sum process $\{S_n(t)=\sum_{i=1}^{[nt]}X_i\}_{t \geq 0}$ associated to i.i.d. regularly varying random elements  $(X_i)_{i \geq 1}$ in $\bD$. We believe that this L\'evy process is a natural infinite-dimensional analogue of the $\alpha$-stable L\'evy motion with values in $\bR^d$, with which it shares several properties, like independence and stationarity of increments, self-similarity, and $\alpha$-stable marginal distributions. We should emphasize that the $\bD$-valued L\'evy motion constructed in the present article is more general than the two-parameter $\alpha$-stable L\'evy sheet introduced in \cite{resnick86} (see Appendix B).

Before we introduce the definition of a L\'evy process with values in $\bD$, we need to recall some basic facts about the space $\bD$.
We denote by $\|\cdot \|$ the supremum norm on $\bD$ given by $\|x\|=\sup_{s \in [0,1]}|x(s)|$,
and by $\SD=\{x \in \bD; \|x\|=1\}$ the unit sphere in $\bD$. With this norm, $\bD$ is a Banach space, but it is not separable. For this reason, the theory of random elements in separable Banach spaces (as presented for instance in \cite{kwapien-woyczynski92}) or the functional limit theorems mentioned in Section 5 of \cite{skorokhod57} cannot be applied to $\bD$.

We endow $\bD$ with Skorokhod's $J_1$-topology, introduced in \cite{skorokhod56}.
There are two equivalent distances which induce this topology. We denote by $d_{J_1}^0$ the distance given by
(12.16) of \cite{billingsley99}, under which $\bD$ is a Polish space (i.e. a complete separable metric space). Note that a function $x \in \bD$ has a countable set of discontinuities which we denote by ${\rm Disc}(x)$. We let $\cD$ be the Borel $\sigma$-field on $\bD$. Since $\cD$  coincides with the $\sigma$-field generated by the projections $\pi_s:\bD \to \bR,s\in [0,1]$ given by $\pi_s(x)=x(s)$, a function $X: \Omega \to \bD$ defined on a probability space $(\Omega,\cF,P)$ is a random element in $\bD$ if $X(s)$ is $\cF$-measurable for any $s \in [0,1]$. For any $s_1,\ldots,s_m \in [0,1]$, the projection $\pi_{s_1,\ldots,s_m}:\bD \to \bR^m$ is defined by $\pi_{s_1,\ldots,s_m}(x)=(x(s_1),\ldots,x(s_m))$. We refer to \cite{billingsley68, billingsley99} for more details.

The analogue of the polar-coordinate transformation is the map $T: \bD_0 \to (0,\infty) \times \SD$ given by $T(x)=\left(\|x\|,\frac{x}{\|x\|} \right)$,
where $\bD_0=\bD \verb2\2 \{0\}$. Let $\nu_{\alpha}$ be the measure on $(0,\infty]$ given by:
\begin{equation}
\label{def-nu-alpha}
\nu_{\alpha}(dr)=\alpha r^{-\alpha-1}1_{(0,\infty)}(r)dr.
\end{equation}


\begin{definition}
\label{def-Levy}
{\rm Let $\nu$ be a measure on $(\bD,\cD)$ such that $\nu(\{0\})=0$ and
\begin{equation}
\label{Levy-measure}
\overline{\nu}:=\nu \circ T^{-1}=c\nu_{\alpha} \times \Gamma_1
\end{equation}
for some $c>0$, $\alpha \in (0,2), \alpha \not=1$ and a probability measure $\Gamma_1$ on $\SD$.
A collection $\{Z(t)\}_{t \geq 0}$ of random elements in $\bD$, defined on a probability space $(\Omega,\cF,P)$ is a {\em $\bD$-valued $\alpha$-stable L\'evy motion} (corresponding to $\nu$) if\\ 
{\em (i)} $Z(0)=0$ a.s.; \\
{\em (ii)} $Z(t_2)-Z(t_1), \ldots, Z(t_K)-Z(t_{K-1})$ are independent, for any $0 \leq t_1 <\ldots<t_K$, $K\geq 3$;\\
{\em (iii)} $Z(t_2)-Z(t_1) \stackrel{d}{=}Z(t_2-t_1)$ for any $0\leq t_1<t_2$, where $\stackrel{d}{=}$ means equality in distribution;\\
{\em (iv)} for any $t>0$, $Z(t)=\{Z(t,s)\}_{s \in [0,1]}$ is an $\alpha$-stable process (with sample paths in $\bD$) such that for any $s_1,\ldots,s_m \in [0,1]$ and for any $u=(u_1,\ldots,u_m) \in \bR^m$,
\begin{align}
\label{ch-fn-Z1}
E\big(e^{iu_1 Z(t,s_1)+\ldots+i u_m Z(t,s_m)}\big)&=\exp \left\{t \int_{\bR^m}(e^{i u \cdot y}-1)\mu_{s_1,\ldots,s_m}(dy) \right\} \quad \mbox{if} \quad \alpha<1, \\
\label{ch-fn-Z2}
E\big(e^{iu_1 Z(t,s_1)+\ldots+i u_m Z(t,s_m)}\big)&=\exp \left\{t \int_{\bR^m}(e^{i u \cdot y}-1-iu \cdot y)\mu_{s_1,\ldots,s_m}(dy)\right\} \quad \mbox{if} \quad \alpha>1
\end{align}
where $y=(y_1,\ldots,y_m)$, $u \cdot y=\sum_{i=1}^{m}u_i y_i$, and $\mu_{s_1,\ldots,s_m}=\nu \circ \pi_{s_1,\ldots,s_m}^{-1}$. }
\end{definition}

From this definition, it follows that $Z(t,s)$ has an $\alpha$-stable $S_{\alpha}(t^{1/\alpha}\sigma_s,\beta_s,0)$-distribution, for some constants $\sigma_s>0$ and $\beta_s \in [-1,1]$ depending on $s$ (see in Proposition \ref{Z-stable} below). Note that property \eqref{Levy-measure} implies that $\int_{\bD_0}(\|x\|^2 \wedge 1)\nu(dx)<\infty$, by a change of variables.

\begin{remark}
\label{remark-cone}
{\rm The authors of \cite{davydov-molchanov-zuyev08} considered $\alpha$-stable L\'evy processes $\{Z(t)\}_{t \geq 0}$ with values in a normed cone $\mathbb K$ with a sub-invariant norm. By definition, these processes have independent and stationary St$\alpha$S increments, where St$\alpha$S stands for ``strictly $\alpha$-stable''. If $\alpha<1$, a $\bD$-valued $\alpha$-stable L\'evy motion (in the sense of Definition \ref{def-Levy}) is an $\alpha$-stable L\'evy process on the cone $\mathbb K=\bD$, and therefore has the series representation given by Theorem 3.10 of \cite{davydov-molchanov-zuyev08}. (Note that the
space $\bD$ equipped with $d_{J_1}^0$ is a normed cone, as specified by Definition 2.6 of \cite{davydov-molchanov-zuyev08}, and the sup-norm $\|\cdot\|$ is sub-invariant, as defined by relation (2.9) of \cite{davydov-molchanov-zuyev08}, i.e. $d_{J_1}^0(x+h,x) \leq \|h\|$ for any $x,h \in \bD$.)
}
\end{remark}

If we denote by $m_{t_1,\ldots,t_n}$ the law of $(Z(t_1),\ldots,Z(t_n))$ on $(\bD^n,\cD^n)$, then by properties {\em (i)}-{\em (iii)}, the family $\{m_{t_1,\ldots,t_n}\}$ of these laws is consistent in the sense of Kolmogorov (see Theorem 3.7 of \cite{PZ07} for a statement of Kolmogorov's consistency theorem for random elements in a Polish space).
But it is not obvious how to ensure that property {\em (iv)} also holds, i.e. it is not clear how to construct a c\`adl\`ag process $\{Z(t,s)\}_{s \in [0,1]}$  with finite-dimensional distributions specified by \eqref{ch-fn-Z1} and \eqref{ch-fn-Z2}. Our first main result will tackle precisely this problem. Moreover, we will show that the process $\{Z(t)\}_{t \geq 0}$ has a {\em modification} $\{\widetilde{Z}(t)\}_{t \geq 0}$ with sample paths in $\bD([0,\infty);\bD)$, where $\bD([0,\infty);\bD)$ is the set of functions $x:[0,\infty) \to \bD$ which are right-continuous and have left-limits with respect to $J_1$.

We introduce the following assumptions on the probability measure $\Gamma_1$.

\vspace{1mm}
{\em Assumption A.} For any $s \in [0,1]$, $\Gamma_1(\{z \in \SD; z(s)=0\})=0$.

{\em Assumption B.} For any $s \in [0,1]$, $\Gamma_1(\{z \in \SD; s \in {\rm Disc}(z)\})=0$.

\vspace{1mm}

We will prove the following result.

\begin{theorem}
\label{main1}
Suppose that Assumption A holds.\\
a) For any measure $\nu$ on $(\bD,\cD)$ such that $\nu(\{0\})=0$ and \eqref{Levy-measure} holds, there exists a $\bD$-valued $\alpha$-stable L\'evy motion $\{Z(t)\}_{t \geq 0}$ (corresponding to measure $\nu$).
\\
b) If $\alpha>1$, suppose that Assumption B holds. Then, there exists a collection $\{\widetilde{Z}(t)\}_{t \geq 0}$ of random elements in $\bD$ such that $P(Z(t)=\widetilde{Z}(t))=1$ for any $t \geq 0$, and the map $t \mapsto \widetilde{Z}(t)$ is in $\bD([0,\infty);\bD)$ with probability $1$.
\end{theorem}

We now turn to our second result, the approximation theorem.

Before speaking about regular variation on $\bD$, we need to recall some classical notions. A non-negative random variable $X$ is {\em regularly varying} of index $-\alpha$ (for some $\alpha>0$) if its tail function $\overline{F}(x)=P(X>x)$ is so (hence the name). A useful characterization of this property is expressed in terms of the vague convergence $n P(X/a_n \in \cdot) \stackrel{v}{\to} \nu_{\alpha}$ of Radon measures on the space $(0,\infty]$, for some sequence $(a_n)_{n \geq 1} \subset \bR_{+}$ with $a_n \uparrow \infty$.
This property can be extended to higher dimensions. A random vector $X$ in $\bR^d$ is regularly varying if
$n P(X/a_n \in \cdot) \stackrel{v}{\to}\mu$ on $\overline{\bR}^d_0=[-\infty,\infty]^d \verb2\2 \{0\}$, for a non-null Radon measure $\mu$ on $\overline{\bR}^d_0$ with $\mu(\overline{\bR}^d_0 \verb2\2 \bR^d)=0$ and a sequence $(a_n)_{n \geq 1} \subset \bR_{+}$ with $a_n \uparrow \infty$; or equivalently,
\begin{equation}
\label{def-RV-Rd}
nP\left(\left(\frac{|X|}{a_n},\frac{X}{|X|}\right) \in \cdot\right) \stackrel{v}{\to} c\nu_{\alpha} \times \Gamma \quad \mbox{on} \quad (0,\infty] \times \bS_d,
\end{equation}
for some $\alpha>0,c>0$ and a probability measure $\Gamma$ on the unit sphere $\bS_d=\{x \in \bR^d;|x|=1\}$ with respect to the Euclidean norm $|\cdot|$ on $\bR^d$. We refer to \cite{resnick87, resnick07} for more details.

Briefly speaking, the regular variation of a random element in $\bR^d$ reduces to the vague convergence of a sequence of Radon measures on the space $\overline{\bR}^d_0$, defined by removing $0$ and adding the $\infty$-hyperplanes. In the case of random elements in $\bD$, there is no natural analogue of an $\infty$-hyperplane. To avoid this problem, the authors of \cite{dehaan-lin01, hult-lindskog05} considered
$$\barD=(0,\infty] \times \SD.$$
Another problem is the fact that vague convergence is defined only for Radon measures on locally compact spaces with countable basis, and $\barD$ is not locally compact. This problem is solved by using the concept of $\widehat{w}$-convergence (defined in Section \ref{section-point} below). 
Note that $\barD$ is a Polish space equipped with the distance $d_{\barD}$ given by:
\begin{equation}
\label{def-d-J1}
d_{\barD}\big((r,z),(r',z')\big)=\left|\frac{1}{r}-\frac{1}{r'}\right| \wedge d_{J_1}^0(z,z'),
\end{equation}
for any $(r,z),(r',z') \in \barD$, with the convention $1/\infty=0$.
With this distance, a set of the form $(\e,\infty] \times \SD$ is bounded in $\barD$. This fact plays an important role in this article.

Since $\|\cdot\|$ is $J_1$-continuous on $\bD$, $T$ is a homeomorphism. Similarly to \cite{roueff-soulier15} (but unlike \cite{hult-lindskog05, davis-mikosch08}), we prefer not to identify $\bD_0$ with $(0,\infty) \times \SD$.
Therefore, we will not say that $\bD_0$ is a subset of $\barD$.
We are now ready to give the definition of regular variation on $\bD$.

\begin{definition}
\label{def-RV}
{\rm A random element $X=\{X(s)\}_{s \in [0,1]}$ in $\bD$ is {\em regularly varying} (and we write $X \in {\rm RV}(\{a_n\}, \overline{\nu},\barD)$) if there exist a sequence $(a_n)_{n \geq 1} \subset \bR_{+}$ with $a_n \uparrow \infty$ and a non-null boundedly finite measure $\overline{\nu}$ on $\barD$ with $\overline{\nu}(\barD \verb2\2 T(\bD_0))=0$ such that
$$nP\left( \left(\frac{\|X\|}{a_n},\frac{X}{\|X\|}\right) \in \cdot \,
\right) \stackrel{\widehat{w}}{\to} \overline{\nu} \quad \mbox{on} \quad \barD.$$
In this case, we say that $\overline{\nu}$ is the {\em limiting measure} of $X$.}
\end{definition}

Since $\overline{\nu}$ is non-null, there exists $a_0>0$ such that $\overline{\nu}((a_0,\infty) \times \SD)>0$. Without loss of generality, we assume that $a_0=1$. We let $c=\overline{\nu}((1,\infty) \times \SD)$.

By Remark 3 of \cite{hult-lindskog05}, the measure $\overline{\nu}$ in Definition \ref{def-RV} has the following property: there exists $\alpha>0$ such that $\overline{\nu}(aA)=a^{-\alpha}\overline{\nu}(A)$ for any $a>0$ and $A \in \cB(\barD)$,
where $aA=\{(ar,z); (r,z) \in A\}$. We say that $\alpha$ is the {\em index} of $X$.
In Lemma \ref{nu-product} (Appendix A), we prove that the measure $\overline{\nu}$ in Definition \ref{def-RV} must be the product measure:
\begin{equation}
\label{def-bar-nu}
\overline{\nu}=c\nu_{\alpha} \times \Gamma_1,
\end{equation}
where $\Gamma_1$ is a probability measure on $\SD$ (called the {\em spectral measure} of $X$), given by
\begin{equation}
\label{def-Gamma1}
\Gamma_1(S)=\frac{\overline{\nu}((1,\infty) \times S)}{c} \quad \mbox{for all} \ S \in \cB(\SD).
\end{equation}
Here we let $\cB(\barD)$ and $\cB(\SD)$ be the classes of Borel sets of $\barD$, respectively $\SD$.

If $X \in {\rm RV}(\{a_n\}, \overline{\nu},\barD)$, then $\|X\|$ is regularly varying of  index $-\alpha$: for any $\e>0$,
\begin{equation}
\label{RV-norm-X}
n P(\|X\|>a_n \e) \to c \e^{-\alpha}, \quad \mbox{as} \quad n \to \infty,
\end{equation}
with the same constant $c>0$ as above.
From this we infer that if $\alpha>1$, $E\|X\|<\infty$, and hence $E|X(s)|<\infty$ for all $s \in [0,1]$. In this case, we define $E[X]=\{E[X(s)]\}_{s \in [0,1]}$.

In \cite{roueff-soulier15}, it is proved that if $S_n=\sum_{i=1}^{n}X_i$, where
$(X_i)_{i\geq 1}$ are i.i.d. random elements in $\bD$ with $X_1 \in RV(\{a_n\},\overline{\nu},\barD)$ and $\alpha>1$,
then
$$\frac{1}{a_n}(S_n-E(S_n)) \stackrel{d}{\to} N \quad \mbox{in} \quad \bD,$$
where $E(S_n)=\{E(S_n(s))\}_{s \in [0,1]}$ and $N=\{N(s)\}_{s \in [0,1]}$ is an $\alpha$-stable process with sample paths in $\bD$ (whose distribution is completely identified).

We are now ready to state our second main result, which is an extension of Theorem 1.1 of \cite{roueff-soulier15} to functional convergence. We let $\bD([0,\infty);\bD)$ be the set of of c\`adl\`ag functions on $[0,\infty)$ with values in $\bD$, equipped with the Skorohod distance $d_{\infty,\bD}$ (described in Section \ref{section-D} below).

\begin{theorem}
\label{main2}
Let $X,(X_i)_{i \geq 1}$ be i.i.d. random elements in $\bD$ such that $X \in RV(\{a_n\},\overline{\nu},\bD)$. 
Let $\alpha$ be the index of $X$ and $\Gamma_1$ be the spectral measure of $X$. Suppose that $\alpha \in (0,2),\alpha \not=1$ and $\Gamma_1$ satisfies Assumptions A and B.
For any $n \geq 1$, $t \geq 0$, let $S_n(t)=\{S_n(t,s)\}_{s \in [0,1]}$, where $S_n(t,s)=a_n^{-1}\sum_{i=1}^{[nt]}X_i(s)$ for $s \in [0,1]$.
Let $\{\widetilde{Z}(t)\}_{t \geq 0}$ be the process constructed in Theorem \ref{main1}.b), which may not be defined on the same probability space as the sequence $(X_i)_{i \geq 1}$.

a) If $\alpha<1$, then
$$S_n(\cdot) \stackrel{d}{\to}\widetilde{Z}(\cdot) \quad in \quad \bD([0,\infty);\bD).$$

b) If $\alpha>1$, let $\overline{S}_n(t)=S_n(t)-E[S_n(t)]$, where $E[S_n(t)]=\{E[S_n(t,s)]\}_{s \in [0,1]}$. If
\begin{equation}
\label{AN}
\lim_{\e \downarrow 0}\limsup_{n \to \infty} \max_{k \leq [nT]} P\left(\left\|\sum_{i=1}^{k} \big(X_i1_{\{\|X_i\| \leq a_n \e \}}-E[X_i 1_{\{\|X_i\|\leq a_n \e \}}]\big) \right\|>a_n \delta\right)=0
\end{equation}
for any $\delta>0$ and $T>0$, then
$$\overline{S}_n(\cdot) \stackrel{d}{\to} \widetilde{Z}(\cdot) \quad in \quad \bD([0,\infty);\bD).$$
\end{theorem}


Assumption B is the same as Condition (A-i) of \cite{roueff-soulier15}, whereas  \eqref{AN} is a stronger form of Condition (A-ii) of \cite{roueff-soulier15}, which is needed for the functional convergence.

We use the following notation. If $(X_n)_{n \geq1}$ and $X$ are random elements in a metric space $(E,d)$, we write $X_n\stackrel{d}{\to}X$ if $(X_n)_n$ converges in distribution to $X$, and $X_n \stackrel{p}{\to} X$ if 
$P(d(X_n,X)>\e)\to0$ for all $\e>0$.

This article is organized as follows. In Section \ref{section-D} we introduce the spaces $\bD([0,1];\bD)$ and $\bD([0,\infty);\bD)$, and we study the weak convergence and tightness of probability measures on these spaces. In Sections \ref{section-construction} and \ref{section-approximation} we give the proofs of Theorems \ref{main1} and \ref{main2}, respectively. Some auxiliary results are included in Appendix \ref{appendixA}.

\section{C\`adl\`ag functions with values in $\bD$}
\label{section-D}

In this section, we introduce the spaces $\bD([0,1];\bD)$ and $\bD([0,\infty);\bD)$ of c\`adl\`ag functions defined $[0,1]$, respectively $[0,\infty)$, with values in $\bD$. These spaces are equipped with the Skorohod distance introduced in \cite{whitt80}.
We examine briefly the weak convergence of probability measures on these spaces, a topic which is
developed at length in the companion paper \cite{balan-saidani18}.

\subsection{The space $\bD([0,1];\bD)$}
\label{section-D01}

In this subsection, we introduce the space $\bD([0,1];\bD)$ and discuss some of its properties.

We begin by recalling some well-known facts about the classical Skorohod space $\bD$. We refer the reader to \cite{billingsley68, billingsley99} for more details.

The Skorohod distance $d_{J_1}$ on $\bD$ is defined as follows: for any $x,y \in \bD$,
$$d_{J_1}(x,y)=\inf_{\lambda \in \Lambda} \{\|\lambda-e\| \vee \|x-y \circ \lambda\| \},$$
where $\Lambda$ the set of strictly increasing continuous functions from $[0,1]$ onto $[0,1]$ and $e$ is the identity function on $[0,1]$. The space $\bD$ equipped with distance $d_{J_1}$ is separable, but it is not complete. There exists another distance $d_{J_1}^0$ on $\bD$, which is equivalent to $d_{J_1}$, under which $\bD$ is complete and separable. This distance is given by:
(see (12.16) of \cite{billingsley99})
\begin{equation}
\label{def-dJ1}
d_{J_1}^0(x,y)=\inf_{\lambda \in \Lambda} \{\|\lambda\|^0 \vee \|x-y \circ \lambda\| \},
\end{equation}
for any $x,y \in \bD$, where
$\|\lambda\|^0=\sup_{s<s'} \left|\log \frac{\lambda(s')-\lambda(s)}{s'-s} \right|$.
Note that
\begin{equation}
\label{d-equal-norm}
d_{J_1}(x,0)=d_{J_1}^0(x,0)=\|x\| \quad \mbox{for any} \quad x \in \bD.
\end{equation}

By relation (12.17) of \cite{billingsley99},
\begin{equation}
\label{lambda0-equiv}
\sup_{s \in [0,1]}|\lambda(s)-s| \leq e^{\|\lambda\|^0}-1 \quad \mbox{for all} \ \lambda \in \Lambda.
\end{equation}
Taking $\lambda=e$ in \eqref{def-dJ1}, we obtain:
\begin{equation}
\label{dJ1-ineq}
d_{J_1}^0(x,y)\leq \|x-y\| \quad \mbox{for all} \quad x,y \in \bD.
\end{equation}
For functions $(x_n)_{n\geq 1}$ and $x$ in $\bD$, we write $x_n \stackrel{J_1}{\to}x$ if $d_{J_1}^0(x_n,x)\to 0$.



For any $\delta\in (0,1)$, we consider the following modulus of continuity of a function $x \in \bD$:
\begin{equation}
\label{def-w-second}
w''(x,\delta)=\sup_{s_1\leq s\leq s_2, s_2-s_1 \leq \delta} \big(|x(s)-x(s_1)|\wedge |x(s_2)-x(s)| \big).
\end{equation}

We denote by $\bD([0,1];\bD)$ the set of functions $x:[0,1] \to \bD$ which are right-continuous and have left limits with respect to $J_1$. We denote by $x(t-)$ the left limit of $x$ at $t \in (0,1]$.
If $x \in \bD([0,1];\bD)$, we let $x(t,s)=x(t)(s)$ for any $t\in [0,1]$ and $s \in [0,1]$.

We let $d_{\bD}$ be {\em Skorohod distance} on $\bD([0,1];\bD)$, given by relation (2.1) of \cite{whitt80}:
\begin{equation}
\label{def-d-D}
d_{\bD}(x,y)=\inf_{\lambda \in \Lambda} \{\|\lambda-e\| \vee \rho_{\bD}(x,y \circ \lambda) \},
\end{equation}
where $\rho_{\bD}$ is the {\em uniform distance} on $\bD([0,1];\bD)$ defined by:
\begin{equation}
\label{def-rho-D}
\rho_{\bD}(x,y)=\sup_{t \in [0,1]}d_{J_1}^{0}(x(t),y(t)).
\end{equation}
Hence, $d_{\bD}(x_n,x) \to 0$ if and only if there exists a sequence $(\lambda_n)_{n\geq 1} \subset \Lambda$ such that
$$\sup_{t \in [0,1]}|\lambda_n(t)-t| \to 0 \quad \mbox{and} \quad \sup_{t \in [0,1]}d_{J_1}^{0}(x_n(\lambda_n(t)),x(t))\to 0.$$

We denote by $\|\cdot\|_{\bD}$ the {\em super-uniform norm} on $\bD([0,1];\bD)$ given by:
$$\|x\|_{\bD}=\sup_{t \in [0,1]}\|x(t)\|.$$
(By the discussion in small print on page 122 of \cite{billingsley99}, the set $\{x(t);t\in [0,1]\}$ is relatively compact in $(\bD,J_1)$, and hence, $\|x\|_{\bD}<\infty$ by Theorem 12.3 of \cite{billingsley99}.)

By relation \eqref{d-equal-norm}, it follows that for any $x \in \bD([0,1];\bD)$,
\begin{equation}
\label{dD-equal-norm}
d_{\bD}(x,0)=\rho_{\bD}(x,0)=\|x\|_{\bD}.
\end{equation}

Note that for any $x,y \in \bD([0,1];\bD)$, we have:
\begin{equation}
\label{dist-ineq1}
d_{\bD}(x,y)\leq \rho_{\bD}(x,y) \leq \|x-y\|_{\bD}.
\end{equation}

The space $\bD([0,1];\bD)$ equipped with $d_{\bD}$ is separable, but it is not complete. Similarly to the distance $d_{J_1}^0$ on $\bD$, we consider another distance $d_{\bD}^0$ on $\bD([0,1];\bD)$, given by:
\begin{equation}
\label{def-d-D-0}
d_{\bD}^0(x,y)=\inf_{\lambda \in \Lambda} \{\|\lambda\|^0 \vee \rho_{\bD}(x,y \circ \lambda) \}.
\end{equation}

The following result is similar to Theorems 12.1 and 12.2 of \cite{billingsley99}. See also Theorem 2.6 of \cite{whitt80}.

\begin{theorem}
\label{bD-CSMS}
The metrics $d_{\bD}$ and $d_{\bD}^0$ are equivalent. The space $\bD([0,1];\bD)$ is separable under $d_{\bD}$ and $d_{\bD}^0$, and is complete under $d_{\bD}^0$.
\end{theorem}

Similarly to \eqref{def-w-second}, for any $x \in \bD([0,1];\bD)$ and $\delta\in (0,1)$, we consider the following modulus of continuity:
$$w_{\bD}''(x,\delta)=\sup_{t_1\leq t\leq t_2,\, \, t_2-t_1 \leq \delta} \big(d_{J_1}^0(x(t),x(t_1))\wedge d_{J_1}^0(x(t_2),x(t)) \big).$$

The following result will be used in the proof of Theorem \ref{tightness-th} below.

\begin{lemma}
\label{dD-sum}
For any $x,y \in \bD([0,1];\bD)$, we have:
$$w_{\bD}''(x+y,\delta) \leq w_{\bD}''(x,\delta)+2\|y\|_{\bD}.$$
\end{lemma}

\noindent {\bf Proof:} Let $t_1 \leq t \leq t_2$ be such that $t_2-t_1 \leq \delta$. By triangle inequality and \eqref{dJ1-ineq},
\begin{align*}
d_{J_1}^0\big(x(t)+y(t),x(t_1)+y(t_1)\big) &\leq d_{J1}^0\big(x(t)+y(t),x(t)\big)+d_{J_1}^0\big(x(t),x(t_1)\big)+
d_{J_1}^0\big(x(t_1),x(t_1)+y(t_1)\big)\\
& \leq \|y(t)\|+d_{J_1}^0\big(x(t),x(t_1)\big)+\|y(t_1)\| \\
&\leq d_{J_1}^0\big(x(t),x(t_1)\big)+2\|y\|_{\bD}.
\end{align*}
Similarly, $d_{J_1}^0\big(x(t)+y(t),x(t_2)+y(t_2)\big) \leq d_{J_1}^0\big(x(t),x(t_2)\big)+2\|y\|_{\bD}$. If $a_1,a_2,b_1,b_2,c \in \bR$ are such that
$a_i\leq b_i+c$ for $i=1,2$, then it is easy to see that $a_1 \wedge a_2 \leq b_1 \wedge b_2+c$.
It follows that $d_{j_1}^0\big(x(t)+y(t),x(t_1)+y(t_1)\big) \wedge d_{J_1}^0\big(x(t)+y(t),x(t_2)+y(t_2)\big)$ is less than
$$d_{J_1}^0\big(x(t),x(t_1)\big) \wedge d_{J_1}^0\big(x(t),x(t_2)\big) +2\|y\|_{\bD} \leq w_{\bD}''(x,\delta)+2 \|y\|_{\bD}.$$
The conclusion follows taking the supremum over all $t_1 \leq t \leq t_2$ such that $t_2-t_1 \leq \delta$. $\Box$

\vspace{3mm}

The following result shows that the super-uniform norm is continuous on $\bD([0,1];\bD)$. Its proof if given in \cite{balan-saidani18}.

\begin{lemma}
\label{dD-cont}
If $(x_n)_{n\geq 1}$ and $x$ are functions in $\bD([0,1];\bD)$ such that $d_{\bD}(x_n,x) \to 0$ as $n \to \infty$, then $\|x_n\|_{\bD} \to \|x\|_{\bD}$ as $n \to \infty$.
\end{lemma}

We conclude this subsection with a brief discussion about finite-dimensional sets in $\bD([0,1];\bD)$, and tightness of probability measures on this space.

Let $\cD_{\bD}$ be the Borel $\sigma$-field of $\bD([0,1];\bD)$, with respect to $d_{\bD}$. It can be shown that $\cD_{\bD}$ coincides with the $\sigma$-field generated by the projections $\{\pi_t^{\bD};t \in [0,1]\}$, where $\pi_{t}^{\bD}:\bD([0,1];\bD) \to \bD$ is given by $\pi_t^{\bD}(x)=x(t)$. We equip $\bD$ with the $J_1$-topology and $\bD([0,1];\bD)$ with distance $d_{\bD}$. Then the projections $\pi_0^{\bD}$ and $\pi_1^{\bD}$ are continuous everywhere, whereas for $t \in (0,1)$, $\pi_t^{\bD}$ is continuous at $x$ if and only if $x$ is continuous at $t$. If $P$ is a probability measure on $\bD([0,1];\bD)$, we let
$T_p$ be the set of $t \in [0,1]$ such that $\pi_t^{\bD}$ is continuous almost everywhere with respect to $P$. The set $T_P$ has a countable complement, and hence is dense in $[0,1]$.
For fixed $t_1,\ldots,t_k \in [0,1]$, we consider the projection
$\pi_{t_1,\ldots,t_k}^{\bD}:\bD([0,1];\bD) \to \bD^k$ given by $\pi_{t_1,\ldots,t_k}^{\bD}(x)=(x(t_1),\ldots,x(t_k))$.

If $(P_n)_{n \geq 1}$ and $P$ are probability measures on $\bD([0,1];\bD)$ such that $P_n \stackrel{w}{\to}P$, then the following marginal convergence holds for all $t_1,\ldots,t_k \in T_P$:
\begin{equation}
\label{marginal-conv}
P_n \circ (\pi_{t_1,\ldots,t_k}^{\bD})^{-1} \stackrel{w}{\to}
P \circ (\pi_{t_1,\ldots,t_k}^{\bD})^{-1} \quad \mbox{in} \quad (\bD^k,J_1^k),
\end{equation}
where $J_1^k$ is the product of $J_1$-topologies.

The following result will be used in the proof of Theorem \ref{tightness-th} below, being the analogue of Theorem 15.3 of \cite{billingsley68} for the space $\bD([0,1];\bD)$. Its proof is given in \cite{balan-saidani18}.

\begin{theorem}[]
\label{tight-th}
A sequence $(P_n)_{n\geq 1}$ of probability measures on $\bD([0,1];\bD)$ is tight if and only if it satisfies the following three conditions:\\
(i)
$\lim_{a \to \infty} \lim_{n\to \infty} P_n(\{x;\|x\|_{\bD}>a\})=0$; \\
(ii) for any $\eta>0$ and $\rho>0$, there exist $\delta \in (0,1)$ and $n_0\geq 1$ such that for all $n\geq n_0$,
$$
\left\{
\begin{array}{ll}
(a) & P_n(\{x;\,w''(x(t),\delta)> \eta \ \mbox{for some} \ t \in [0,1]\})<\rho \\
(b) & P_n(\{x;\,|x(t,\delta)-x(t,0)|> \eta \ \mbox{for some} \ t \in [0,1]\})<\rho \\
(c) & P_n(\{x;\,|x(t,1-)-x(t,1-\delta)|> \eta \ \mbox{for some} \ t \in [0,1]\})<\rho;
\end{array}
\right.$$
(iii)  for any $\eta>0$ and $\rho>0$, there exist $\delta \in (0,1)$ and $n_0\geq 1$ such that for all $n\geq n_0$,
$$
\left\{
\begin{array}{ll}
(a) & P_n(\{x;\,w_{\bD}''(x,\delta)> \eta \})<\rho \\
(b) & P_n(\{x;\,d_{J_1}^{0}\big(x(\delta),x(0)\big)> \eta \})<\rho \\
(c) & P_n(\{x;\,d_{J_1}^{0}\big(x(1-),x(1-\delta)\big)> \eta \})<\rho.
\end{array}
\right.$$
\end{theorem}

\subsection{The space $\bD([0,\infty);\bD)$}

In this subsection, we introduce the space $\bD([0,\infty);\bD)$ and we list some of its properties.

For any fixed $T>0$, we let $\bD([0,T];\bD)$ be the set of functions $x:[0,T] \to \bD$ which are right-continuous and have left-limits with respect to $J_1$. Let $\Lambda_{T}$ be the set of strictly increasing continuous functions from $[0,T]$ onto itself. Similarly to the case $T=1$, we define the Skorohod distance on $\bD([0,T];\bD)$ by:
\begin{equation}
\label{def-d-TD}
d_{T,\bD}(x,y)=\inf_{\lambda \in \Lambda_T} \{\|\lambda-e\|_T \wedge \rho_{T,\bD}(x,y \circ \lambda) \},
\end{equation}
where $\|\cdot\|_{T}$ is the supremum norm on $\Lambda_T$, $e$ is the identity function on $[0,T]$, and $\rho_{T,\bD}$ is the uniform distance on $\bD([0,T];\bD)$ given by:
\begin{equation}
\label{def-rho-TD}
\rho_{T,\bD}(x,y)=\sup_{t \in [0,T]}d_{J_1}^{0}(x(t),y(t)).
\end{equation}

We denote by $\|\cdot\|_{T,\bD}$ the super-uniform norm on $\bD([0,T];\bD)$ given by:
$$\|x\|_{T,\bD}=\sup_{t \in [0,T]}\|x(t)\|.$$
For any $x,y \in \bD([0,T];\bD)$, we have
\begin{equation}
\label{dist-ineq}
d_{T,\bD}(x,y)\leq \rho_{T,\bD}(x,y) \leq \|x-y\|_{T,\bD}.
\end{equation}

The Skorohod distance on the space $\bD([0,\infty);\bD)$ is given by: (see (2.2) of \cite{whitt80})
\begin{equation}
\label{def-d-infty}
d_{\infty,\bD}(x,y)=\int_{0}^{\infty}e^{-t} \Big( d_{t,\bD}\big(r_t(x),r_t(y)\big)\wedge 1\Big)dt,
\end{equation}
where $r_t(x)$ is the restriction to $[0,t]$ of the function $x \in \bD([0,\infty);\bD)$.

By Theorem 2.6 of \cite{whitt80}, $\bD([0,\infty);\bD)$ equipped with distance $d_{\infty,\bD}$ is a Polish space. Its Borel $\sigma$-field $\cD_{\infty,\bD}$ coincides (by Lemma 2.7 of \cite{whitt80}) with the $\sigma$-field generated by the projections $\{\pi_t^{\bD};t  \geq 0\}$, where $\pi_{t}^{\bD}:\bD([0,\infty);\bD) \to \bD$ is given by $\pi_t^{\bD}(x)=x(t)$.

Similarly to page 174 of \cite{billingsley99}, if $(P_n)_{n \geq 1}$ and $P$ are probability measures on $\bD([0,\infty);\bD)$ such that $P_n \stackrel{w}{\to}P$ then the marginal convergence \eqref{marginal-conv} holds for all $t_1,\ldots,t_k \in T_P$, where the set $T_p$ (defined as in Section \ref{section-D01} above) has a countable complement. In fact, $P_n \stackrel{w}{\to}P$ if and only if $P_n \circ r_t^{-1} \stackrel{w}{\to} P \circ r_t^{-1}$ for any $t \in T_p$ (see also Theorem 2.8 of \cite{whitt80}).

\section{Construction: proof of Theorem \ref{main1}}
\label{section-construction}

In this section, we give the construction of an $\alpha$-stable L\'evy motion $Z=\{Z(t)\}_{t \geq 0}$ with values in $\bD$, and we show that this process has a modification with sample paths in the space of c\`adl\`ag functions from $[0,\infty)$ to $\bD$. We follow the method described in Section 5.5 of \cite{resnick07}. For each $t\geq 0$, $Z(t)$ is a random element in $\bD$ which we denote by $\{Z(t,s)\}_{s \in [0,1]}$, that is $Z(t,s)=Z(t)(s)$. Intuitively, the process $Z$ evolves in time and space: $Z(t,s)$ gives the value of this process at time $t\geq 0$ and location $s \in [0,1]$ in space.

\subsection{The compound Poisson building blocks}

In this subsection, we introduce the building blocks of the construction, and we examine their properties.

Let $N=\sum_{i\geq 1}\delta_{(T_i,R_i,W_i)}$ be a Poisson random measure on $[0,\infty) \times \barD$ of intensity ${\rm Leb} \times \overline{\nu}$, defined on a complete probability space $(\Omega,\cF,P)$, where ${\rm Leb}$ is the Lebesgue measure and $\overline{\nu}$ is given by \eqref{Levy-measure} on $(0,
\infty) \times \SD$ and $\overline{\nu}(\{\infty\}\times \SD)=0$. (Refer to Definition \ref{def-PRM} below for the definition of a Poisson random measure.)

By an extension of Proposition 5.3 of \cite{resnick07} to point processes on Polish spaces, we can represent the points $(T_i,R_i,W_i)$ as follows: $\{(T_i,R_i)\}_{i\geq 1}$ are the points of a Poisson random measure on $[0,\infty) \times (0,\infty]$ of intensity ${\rm Leb} \times \nu_{\alpha}$, and
$(W_i)_{i\geq 1}$ is an independent sequence of i.i.d. random elements in $\SD$ with law $\Gamma_1$.

Let $(\e_j)_{j \geq 0}$ be a sequence of real numbers such that $\e_j \downarrow 0$ and $\e_0=1$. Let $I_j=(\e_j,\e_{j-1}]$ for $j\geq 1$ and $I_0=(1,\infty)$. We fix $t\geq 0$ and $s \in [0,1]$. For any $j \geq 0$, we let
\begin{equation}
\label{def-Zj}
Z_j(t,s)=\int_{[0,t] \times I_j \times \SD}rz(s)N(du,dr,dz)=\sum_{T_i \leq t}R_i W_i(s)1_{\{R_i\in I_j\}}.
\end{equation}
Note that for any $j \geq 0$ and $s \in [0,1]$, $Z_j(0,s)=0$.

\begin{lemma}
\label{Zj-well-def}
a) $Z_j(t,s)$ is well-defined and $\cF$-measurable for any $j\geq 0, t \geq 0, s\in [0,1]$.
b) For any $t\geq 0$ and $j \geq 0$, the process $Z_j(t)=\{Z_j(t,s)\}_{s \in [0,1]}$ has all sample paths in $\bD$, with left limit at point $s\in (0,1]$ given by
$$Z_j(t,s-)=\int_{[0,t] \times I_j \times \SD}rz(s-)N(du,dr,dz)=\sum_{T_i \leq t}R_i W_i(s-)1_{\{R_i\in I_j\}}.$$
 \end{lemma}

\noindent {\bf Proof:} a) $Z_j(t,s)$ is well-defined since $[0,t] \times I_j \times \SD$ is a {\em bounded} set in $[0,\infty) \times \barD$ (due to definition \eqref{def-d-J1} of the metric $d_{\barD}$ on $\barD$), and the sum in \eqref{def-Zj} contains finitely many terms. $Z_j(t,s)$ is $\cF$-measurable since $N$ is a point process and the map $\mu \mapsto \mu(\overline{\pi}_s)=\int_{(0,\infty) \times \SD} rz(s)\mu(dr,dz)$ is $\cM_p([0,\infty) \times \barD)$-measurable, where $\overline{\pi}_s(r,z)=rz(s)$ (see Section \ref{section-point} below for the definition of a point process).

b) This follows by the dominated convergence theorem, whose application is justified by the fact that $\int_{[0,t] \times I_j \times \SD} r N(du,dr,dz)<\infty$. $\Box$

\vspace{3mm}

To investigate the finite dimensional distribution of process $Z_j(t)$ corresponding to points $s_1,\ldots,s_m \in [0,1]$, we consider the function $\overline \pi_{s_1,\ldots,s_m}:(0,\infty) \times \SD \to \bR^m$ given by: $$\overline{\pi}_{s_1,\ldots,s_m}(r,z)=(rz(s_1),\ldots,rz(s_m)).$$
Note that $\overline{\pi}_{s_1,\ldots,s_m}\circ T=\pi_{s_1,\ldots,s_m}$.

\begin{lemma}
\label{ch-function-Zj}
For any $j \geq 0$, $t\geq 0$ and $s_1,\ldots,s_m \in [0,1]$, the vector $(Z_j(t,s_1),\ldots,Z_j(t,s_m))$ has a compound Poisson distribution in $\bR^m$ with characteristic function:
$$E\big(e^{i\sum_{k=1}^{m}u_k Z_j(t,s_k)}\big)=\exp\left\{t \int_{I_j \times \SD}(e^{iu_1rz(s_1)+\ldots+iu_m rz(s_m)}-1)\overline{\nu}(dr,dz) \right\},$$
for any $(u_1,\ldots,u_m) \in \bR^m$.
Letting $\varphi(s)=\int_{\SD}z(s) \Gamma_1(s)$ and $\psi(s)=\int_{\SD}|z(s)|^2\Gamma_1(ds)$ for any $s \in [0,1]$, we have
$$E\big(Z_j(t,s)\big)=t\int_{I_j \times \SD}rz(s)\overline{\nu}(dr,dz)=t \varphi(s)\int_{I_j}r \nu_{\alpha}(dr)$$
$${\rm Var}\big(Z_j(t,s)\big)=t\int_{I_j \times \SD}|rz(s)|^2\overline{\nu}(dr,dz)=t \psi(s)\int_{I_j}r^2 \nu_{\alpha}(dr).$$
\end{lemma}

\noindent {\bf Proof:} We represent the restriction of $N$ to $[0,t] \times I_j \times \SD$ as $N|_{[0,t]\times I_j \times \SD}\stackrel{d}{=}\sum_{i=1}^{K}\delta_{(\tau_i,J_i,W_i)}$, where $K$ is a Poisson random variable of mean $t \overline{\nu}(I_j \times \SD)$, $(\tau_i)_{i\geq 1}$ are i.i.d. uniformly distributed on $[0,1]$, $(J_i)_{i\geq 1}$ are i.i.d. on $I_j$ of law $\nu_{\alpha}/\nu_{\alpha}(I_j)$, $(W_i)_{i\geq 1}$ are i.i.d. on $\SD$ of law $\Gamma_1$, and $K, (\tau_i)_{i\geq 1}, (J_i)_{i\geq 1}, (W_i)_{i\geq 1}$ are independent. Hence, $(Z_j(t,s_1),\ldots,Z_j(t,s_m)) \stackrel{d}{=}\sum_{i=1}^K J_i Y_i$ with $Y_i=(W_i(s_1),\ldots,W_i(s_m))$. The result follows since $\{J_iY_i\}_{i\geq 1}$ are i.i.d. vectors in $\bR^m$ with law
$$\frac{1}{\overline{\nu}(I_j \times \SD)} \overline{\nu}|_{I_j \times \SD}\circ \overline{\pi}_{s_1,\ldots,s_m}^{-1},$$
where $\overline{\nu}|_{I_j \times \SD}$ is the restriction of $\overline{\nu}$ to $I_j \times \SD$. $\Box$
\vspace{3mm}

The previous result shows that for $j\geq 1$, $Z_j(t,s)$ has finite mean and finite variance, while $Z_0(t,s)$ has infinite variance (since $\alpha<2$), but has finite mean if $\alpha>1$. Note that $\sum_{j \geq 1}{\rm Var}\big(Z_j(t,s)\big)=t\psi(s) \int_{(0,1]}r^2 \nu_{\alpha}(dr)<\infty$. Moreover, the variables $\{Z_j(t,s)\}_{j \geq 0}$ are independent, since the intervals $(I_j)_{j \geq 0}$ are disjoint. Hence by Kolmogorov's convergence criterion (see e.g. Theorem 22.6 of \cite{billingsley95}), for any $t > 0$ and $s \in [0,1]$,
$$\sum_{j\geq 1}\Big(Z_j(t,s) -E\big(Z_j(t,s)\big)\Big) \quad {\rm converges \ a.s.}$$
We denote by $\Omega_{t,s}$ the event that this series converges, with $P(\Omega_{t,s})=1$.

If $\alpha<1$,
$\sum_{j \geq 1}E\big(Z_j(t,s)\big)=c\varphi(s) \int_{0}^{1} r\nu_{\alpha}(dr)$ is finite, whereas if $\alpha>1$, $E(Z_0(t,s))=c \varphi(s) \int_1^{\infty}r \nu_{\alpha}(dr)$ is finite.
For any $t \geq 0, s \in [0,1]$ fixed, on the event $\Omega_{t,s}$ we define
\begin{align}
\label{def-Z1}
\overline{Z}(t,s)&=\sum_{j \geq 0}Z_j(t,s)
\quad \mbox{if} \quad \alpha<1, \\
\label{def-Z2}
\overline{Z}(t,s)&=
\sum_{j \geq 0}\Big(Z_j(t,s)-E\big(Z_j(t,s)\big)\Big)
\quad \mbox{if} \quad \alpha>1.
\end{align}
On the event $\Omega_{t,s}^c$, we let $\overline{Z}(t,s)=x_0$, for arbitrary $x_0 \in \bD$, in both cases $\alpha<1$ and $\alpha>1$. Note that $\overline{Z}(0,s)=0$ for all $s \in [0,1]$.

For any $s_1,\ldots,s_m \in [0,1]$, we consider the following measure on $\bR^m$:
\begin{equation}
\label{def-mu-s}
\mu_{s_1,\ldots,s_m}=\nu \circ \pi_{s_1,\ldots,s_m}^{-1}=\overline{\nu} \circ \overline{\pi}_{s_1,\ldots,s_m}^{-1},
\end{equation}

The next result identifies some essential properties of the measures $\mu_{s_1,\ldots,s_m}$.
Assumption A is needed only to guarantee that $\mu_{s_1,\ldots,s_m}(\{0\})=0$.

\begin{lemma}
\label{mu-properties}
Suppose that Assumption A holds. \\
a) For any $s_1, \ldots,s_m \in [0,1]$, $\mu_{s_1,\ldots,s_m}$ is a L\'evy measure on $\bR^m$, i.e.
$$\mu_{s_1,\ldots,s_m}(\{0\})=0 \quad \mbox{and} \quad \int_{\bR^m}(|y|^2 \wedge 1)\mu_{s_1,\ldots,s_m}(dy)<\infty.$$
b) For any $s_1,\ldots,s_m \in [0,1]$,  for any $h>0$ and for any Borel set $A \subset \bR^m$,
$$\mu_{s_1,\ldots,s_m}(hA)=h^{-\alpha}\mu_{s_1,\ldots,s_m}(A).$$
c) For any $s\in [0,1]$, the measure $\mu_s$ is given by
$$\mu_s(dy)=\big(c_s^{+}\alpha y^{-\alpha-1}+c_s^{-} \alpha (-y)^{-\alpha-1}\big)dy,$$
where $c_s^{+}=\mu_s(1,\infty)$ and $c_s^{-}=\mu_s(-\infty,-1)$.
\end{lemma}

\noindent {\bf Proof:} a) By Assumption A,
$\mu_{s_1,\ldots,s_m}(\{0\})=\overline{\nu}(\{(r,z);rz(s_1)=\ldots=rz(s_m)=0\})=0$, using the convention that $\infty \cdot 0=0$. The second property follows because
\begin{align*}
& \int_{|y|\leq 1}|y|^2 \mu_{s_1,\ldots,s_m}(dy)=c \int_{\SD}\left(\int_{0}^{\left(\sum_{i=1}^{m}|z(s_i)|^2 \right)^{-1/2}}r^2 \nu_{\alpha}(dr)\right)\sum_{i=1}^{m}|z(s_i)|^2\Gamma_1(dz)\\
& \quad =c\,\frac{\alpha}{2-\alpha}\int_{\SD}\Big(\sum_{i=1}^{m}|z(s_i)|^{2}\Big)^{\alpha/2}\Gamma_1(dz)
\leq c\, \frac{\alpha}{2-\alpha},
\end{align*}
and
$$\int_{|y|>1}\mu_{s_1,\ldots,s_m}(dy)=c \int_{\SD}\int_{\big(\sum_{i=1}^{m}|z(s_i)|^2\big)^{-1}}^{\infty}\nu_{\alpha}(dr)
\Gamma_1(dz)=c
\int_{\SD}\Big(\sum_{i=1}^{m}|z(s_i)|^2\Big)^{\alpha/2}\Gamma_1(dz).$$

b) By Fubini's theorem and the scaling property of $\nu_{\alpha}$, it can be proved that $\overline{\nu}$ has the following scaling property: for any $h>0$ and $H \in \cB(\bD_0)$, $\overline{\nu}(hH)=h^{-\alpha}\overline{\nu}(H)$, where $hH=\{(hr,z); (r,z)\in H\}$. For any $h>0$ and $A \in \cB(\bR^m)$, we have
$$\mu_{s_1,\ldots,s_m}(hA)=\overline{\nu}(\{(r,z);(rz(s_1),\ldots,rz(s_m))\in hA\})=\overline{\nu}(hH)$$
where $H=\{(r,z);(rz(s_1),\ldots,rz(s_m))\in A\}=\overline{\pi}_{s_1,\ldots,s_m}^{-1}(A)$. The conclusion follows from the scaling property of $\overline{\nu}$ mentioned above.

c) This is an immediate consequence of the scaling property in b). $\Box$

\vspace{3mm}

We denote by $S_{\alpha}(\sigma,\beta,\mu)$ the $\alpha$-stable distribution given by Definition 1.1.6 of \cite{ST94}, and
\begin{equation}
\label{def-C}
C_{\alpha}^{-1}=\frac{\Gamma(2-\alpha)}{1-\alpha}\cos \left(\frac{\pi \alpha}{2}\right).
\end{equation}

Based on the previous lemma, we obtain the following result.

\begin{proposition}
\label{Z-stable}
For any $t >0$, the process $\overline{Z}(t)=\{\overline{Z}(t,s)\}_{s \in [0,1]}$ given by \eqref{def-Z1} and \eqref{def-Z2} is $\alpha$-stable with finite-dimensional distributions given by \eqref{ch-fn-Z1} and \eqref{ch-fn-Z2}. In particular, for any $t>0$ and $s \in [0,1]$, $\overline{Z}(t,s)$ has a $S_{\alpha}(t^{1/\alpha}\sigma_s,\beta_s,0)$ distribution with parameters
\begin{equation}
\label{def-sigma-beta}
\sigma_s=C_{\alpha}^{-1}(c_s^{+}+c_s^{-}) \quad \mbox{and} \quad \beta_s=\frac{c_s^{+}-c_s^{-}}{c_s^{+}+c_s^{-}},
\end{equation}
where $c_s^+$ and $c_s^{-}$ are given in Lemma \ref{mu-properties}.c). Moreover, $\overline{Z}(t,s_k) \stackrel{d}{\to} \overline{Z}(t,s)$ as $k \to \infty$,
for any $s \in [0,1]$ and for any sequence $(s_k)_{k \geq 1}$ with $s_k \to s$ and $s_k \geq s$ for all $k \geq 1$.
\end{proposition}

\noindent {\bf Proof:} {\em Case 1: $\alpha<1$.} By Lemma \ref{ch-function-Zj} and the independence of $\{Z_j(t,s)\}_{j \geq 0}$, it follows that the characteristic function of the variable $\overline{Z}(t,s)$ is given by:
$$E\big(e^{iu \overline{Z}(t,s)}\big)=\exp \left\{t \int_{\barD}(e^{iurz(s)}-1)\overline{\nu}(dr,dz) \right\}=\exp \left\{t \int_{\bR}(e^{iuy}-1)\mu_s(dy) \right\}, \quad u \in \bR.$$
The fact that $\overline{Z}(t,s)$ has a $S_{\alpha}(t^{1/\alpha}\sigma_s,\beta_s,0)$ follows essentially from the calculations on page 568 of \cite{feller71}, using the form of the measure $\mu_s$ given in Lemma \ref{mu-properties}.c).

Similarly, it can be seen that for any $s_1,\ldots,s_m \in [0,1]$, $(\overline{Z}(t,s_1),\ldots,\overline{Z}(t,s_m))$ has characteristic function given by \eqref{ch-fn-Z1}.
The fact that $(\overline{Z}(t,s_1),\ldots, \overline{Z}(t,s_m))$ has an $\alpha$-stable distribution follows by Theorem 14.3 of \cite{sato99}, using the scaling property of the measure $\mu_{s_1,\ldots,s_m}$ given in Lemma \ref{mu-properties}.b).

The last statement follows from the fact that $E(e^{iu \overline{Z}(t,s_k)}) \to E(e^{iu \overline{Z}(t,s)})$. To see this, note that $\lim_{k \to \infty}z(s_k) = z(s)$ for any $z \in \SD$. By the dominated convergence theorem,
$$\int_{\barD}(e^{iurz(s_k)}-1)\overline{\nu}(dr,dz) \to \int_{\barD}(e^{iurz(s)}-1)\overline{\nu}(dr,dz), \quad \mbox{as}  \quad k \to \infty.$$ The application of this theorem is justified using the inequalities $|e^{iurz(s)}-1| \leq |urz(s)|$ if $r\leq 1$ and $|e^{iurz(s)}-1| \leq 2$ if $r>1$.

{\em Case 2: $\alpha>1$.} This is similar to {\em Case 1}, except that we now have centering constants. In this case, the characteristic function of $\overline{Z}(t,s)$ is given by
$$E\big(e^{iu \overline{Z}(t,s)}\big)=\exp \left\{t \int_{\barD}(e^{iurz(s)}-1-iurz(s))\overline{\nu}(dr,dz) \right\}, \quad u \in \bR.$$
The last statement follows from the fact that $E(e^{iu \overline{Z}(t,s_k)}) \to E(e^{iu \overline{Z}(t,s)})$, since
$$\int_{\barD}(e^{iurz(s_k)}-1-iurz(s))\overline{\nu}(dr,dz) \to \int_{\barD}(e^{iurz(s)}-1-iurz(s))\overline{\nu}(dr,dz).$$
The application of the dominated convergence theorem is justified using the inequalities $|e^{iurz(s)}-1-iurz(s)| \leq \frac{1}{2}|urz(s)|^2$ if $r\leq 1$ and $|e^{iurz(s)}-1-iurz(s)| \leq 2|urz(s)| $ if $r>1$.
$\Box$

\vspace{3mm}

We denote by $\bD_u([0,\infty);\bD)$ the set of functions $x:[0,\infty) \to \bD$ which are right-continuous and have left limits with respect to the uniform norm $\|\cdot\|$ on $\bD$. Clearly, $\bD_u([0,\infty);\bD)$ is a subset of $\bD([0,\infty);\bD)$.

\begin{lemma}
\label{Zj-in-Du}
For any $j\geq 0$, the process $\{Z_j(t)\}_{t \geq 0}$ has all sample paths in $\bD_{u}([0,\infty);\bD)$, with left limit at $t>0$ given by $Z_j(t-)=\{Z_j(t-,s)\}_{s \in [0,1]}$, where
$$Z_j(t-,s)=\int_{[0,t) \times I_j \times \SD}rz(s)N(du,dr,dz).$$
\end{lemma}

\noindent {\bf Proof:} We first show that the map $t \mapsto Z_j(t)$ is right-continuous in $(\bD,\|\cdot\|)$. Let $t\geq 0$ be arbitrary and $(t_n)_{n\geq 1}$ such that $t_n \to t$ and $t_n \geq t$ for all $n \geq 1$. Then
$$\|Z_j(t_n)-Z_j(t)\|=\sup_{s \in [0,1]}\left|\int_{(t,t_n] \times I_j \times \SD}rz(s)N(du,dr,dz)\right|\leq \int_{(t,t_n] \times I_j \times \SD}r N(du,dr,dz),$$
and the last integral converges to $0$ as $n \to \infty$ by the dominated convergence theorem. Next, we show that the map $t \mapsto Z_j(t)$ has left limit $Z_j(t-)$ in $(\bD,\|\cdot\|)$. Let $t>0$ be arbitrary and $(t_n)_{n\geq 1}$ such that $t_n \to t$ and $t_n \leq t$ for all $n \geq 1$. Then
$$\|Z_j(t-)-Z_j(t_n)\|=\sup_{s \in [0,1]}\left|\int_{(t_n,t) \times I_j \times \SD}rz(s)N(du,dr,dz)\right|\leq \int_{(t_n,t) \times I_j \times \SD}rN(du,dr,dz),$$
and the last integral converges to $0$ as $n \to \infty$ by the dominated convergence theorem. $\Box$

\vspace{3mm}

For any $\e>0$, $t \geq 0$ and $s \in [0,1]$, we let
\begin{equation}
\label{def-Z-e}
Z^{(\e)}(t,s)=\int_{[0,t] \times (\e,\infty)\times \SD}rz(s)N(du,dr,dz)=\sum_{T_i\leq t}R_i W_i(s)1_{\{R_i\in (\e,\infty)\}}.
\end{equation}
Using this notation, we have:
\begin{equation}
\label{def-Z-e2}
Z^{(\e_k)}(t,s)=\sum_{j=0}^{k}Z_j(t,s), \quad \mbox{for all} \quad k \geq 0.
\end{equation}

\begin{remark}
\label{Z-e-in-Du}
{\rm
Similarly to Lemma \ref{Zj-well-def} and Lemma  \ref{Zj-in-Du} for $j=0$, it can be proved that
the process $Z^{(\e)}(t)=\{Z^{(\e)}(t,s)\}_{s \in [0,1]}$ has all sample paths in $\bD$ for any $t \geq 0$, and the process $Z^{(\e)}=\{Z^{(\e)}(t)\}_{t \geq 0}$ has all sample paths in $\bD_u([0,\infty);\bD)$.
}
\end{remark}

\subsection{Construction in the case $\alpha<1$}
\label{subsection-alpha1}

In this subsection, we give the proof Theorem \ref{main1} in the case $\alpha<1$. In particular, property \eqref{Zk-converges-Z} below will be used in the proof of the approximation result (Theorem \ref{main2}.a)).

Our first result shows that for any $t>0$ fixed, 
the process $\overline{Z}(t)$ given by \eqref{def-Z1} has a c\`adl\`ag modification which can be obtained as an almost sure limit with respect to the uniform norm. Recall that $\{X(s)\}_{s \in [0,1]}$ is a {\em modification} of
$\{Y(s)\}_{s \in [0,1]}$ if $P(X(s)=Y(s))=1$ for all $s \in [0,1]$.

\begin{lemma}
\label{Zk-Z-D}
If $\alpha<1$, then for any $t\geq 0$, there exists a random element $Z(t)=\{Z(t,s)\}_{s \in [0,1]}$ in $\bD$ such that
$P(Z(t,s)=\overline{Z}(t,s))=1$ for all $s \in [0,1]$, and
$$\lim_{k \to \infty}\|Z^{(\e_k)}(t)-Z(t) \| = 0 \quad a.s.$$
\end{lemma}

\noindent {\bf Proof:} For $t=0$, we define $Z(0,s)=0$ for all $s \in [0,1]$. We consider the case $t>0$.
By \eqref{def-Zj}, $\|Z_j(t)\| \leq \sum_{i \geq 1}R_i 1_{\{R_i\in I_j\}}1_{\{T_i \leq t\}}=\int_{[0,t] \times I_j \times \SD} r N(du,dr,dz)$. Since $\alpha<1$, it follows that
\begin{align*}
E \sum_{j \geq 1}\|Z_j(t)\|\leq E \sum_{j\geq 1}\int_{[0,t] \times I_j \times \SD} r N(du,dr,dz)=t \int_{(0,1] \times \SD} r\overline{\nu}(dr,dz)<\infty,
\end{align*}
which implies that $\sum_{j \geq 1}\|Z_j(t)\|<\infty$ a.s. We denote by $\Omega_t$ the event that this series converges, with $P(\Omega_t)=1$. On the event $\Omega_t$, the sequence $\{Z^{(\e_k)}(t)=\sum_{j=0}^{k}Z_j(t)\}_{k\geq 0}$ is Cauchy in $(\bD,\|\cdot\|)$, and we denote its limit by $Z(t)$. On the event $\Omega_t^c$, we let $Z(t)=x_0$.

By Lemma \ref{Zj-well-def}.a), $Z(t,s)$ is $\cF$-measurable for any $s\in [0,1]$. Hence, $Z(t)$ is a random element in $\bD$.
On the event $\Omega_{t,s} \cap \Omega_t$, $\overline{Z}(t,s)-Z^{(\e_k)}(t,s)=\sum_{j\geq k+1}Z_j(t,s)$, and hence
$$|Z^{(\e_k)}(t,s)-\overline{Z}(t,s)|\leq \sum_{j\geq k+1}|Z_j(t,s)|\leq \sum_{j \geq k+1}\|Z_j(t)\| \to 0.$$
On the other hand, on the event $\Omega_t$, $Z^{(\e_k)}(t,s) \to Z(t,s)$ for any $s\in [0,1]$.
By the uniqueness of the limit, $Z(t,s)=\overline{Z}(t,s)$ on the event $\Omega_{t,s} \cap \Omega_t$.
$\Box$

\vspace{3mm}

{\bf The following result proves Theorem \ref{main1}.a) in the case $\alpha<1$.}

\begin{theorem}
\label{existence-Levy-process1}
If $\alpha<1$, the process $\{Z(t)\}_{t \geq 0}$ defined in Lemma \ref{Zk-Z-D} is a $\bD$-valued $\alpha$-stable L\'evy motion (corresponding to $\nu$). This process is $(1/\alpha)$-self-similar, i.e.
\begin{equation}
\label{self-similar}
\{Z(ct)\}_{t \geq 0}\stackrel{d}{=} c^{1/\alpha}\{Z(t)\}_{t \geq 0} \quad \mbox{for any} \quad c>0,
\end{equation}
where $\stackrel{d}{=}$ denotes equality of finite-dimensional distributions.
\end{theorem}

\noindent {\bf Proof:} We first show that the process $\{Z(t)\}_{t \geq 0}$ satisfies properties {\em (i)}-{\em (iv)} given in Definition \ref{def-Levy}.
Property {\em (i)} is clear. To verify property {\em (ii)}, we apply Lemma \ref{lemmaB} (Appendix A) to the space $S=\bD$ equipped with $d_{J_1}^0$.
By Lemma \ref{Zk-Z-D}, for $i=2,\ldots,K$, $X_k^{(i)}:=Z^{(\e_k)}(t_{i})-Z^{(\e_k)}(t_{i-1}) \to X^{(i)}:=Z(t_{i})-Z(t_{i-1})$ a.s. as $k \to \infty$, in $(\bD,\|\cdot\|)$, and hence also in $(\bD,J_1)$. The variables $X_k^{(2)},\ldots, X_{k}^{(K)}$ are independent for any $k$, since
 $X_k^{(i)}$ is $\cF_{t_{i-1},t_{i}}^N$-measurable and the $\sigma$-fields $\cF_{t_{i-1},t_i}^N,i=2,\ldots,K$ are independent. Here $\cF_{s,t}^N$ is the $\sigma$-field generated by $N((a,b] \times B)$ for any $s<a<b \leq t$ and $B \in \cB(\barD)$. It follows that $X^{(2)},\ldots,X^{(K)}$ are independent.

For property {\em (iii)}, we have to show that vectors $X:=(Z(t_2,s_1)-Z(t_1,s_1), \ldots Z(t_2,s_1)-Z(t_1,s_m))$  and $Y:=(Z(t_2-t_1,s_1),\ldots,Z(t_2-t_1,s_m))$ have the same distribution, for any $s_1,\ldots,s_m \in [0,1]$. By \eqref{def-Z1} and Lemma \ref{Zk-Z-D}, on the event $\Omega_{t_1,s} \cap \Omega_{t_2,s} \cap \Omega_{t_1}\cap \Omega_{t_2}$,
$$Z(t_2,s)-Z(t_1,s)=\overline{Z}(t_2,s)-\overline{Z}(t_1,s)=\sum_{j\geq 0}\big(Z_j(t_2,s)-Z_j(t_1,s)\big).$$
As in the proof of Proposition \ref{Z-stable}, it follows that the characteristic function of $X$ is
$$E(e^{iu \cdot X})=\exp\left\{(t_2-t_1)\int_{\bR^m}(e^{iu \cdot y}-1)\mu_{s_1,\ldots,s_m}(dy) \right\}, \quad u \in \bR^m,$$
which is the same as the characteristic function of $Y$. Hence $X \stackrel{d}{=}Y$. Finally, property {\em (iv)} was shown in Proposition \ref{Z-stable} for $\overline{Z}(t)$, and remains valid for its modification $Z(t)$.

To prove relation \eqref{self-similar}, we have to show that $\{Z(ct)\}_{t \geq 0} \stackrel{d}{=} \{c^{1/\alpha}Z(t)\}_{t \geq 0}$ for any $c>0$. Since both processes have stationary and independent increments, it is enough to show that $Z(ct) \stackrel{d}{=}c^{1/\alpha}Z(t)$ for any $t>0$, i.e. vectors $U=(Z(ct,s_1),\ldots,Z(ct,s_m))$ and $V=c^{1/\alpha}(Z(t,s_1),\ldots, Z(t,s_m))$ have the same distribution, for any $s_1, \ldots,s_m \in [0,1]$ and $t>0$. Let $h_c(y)=c^{1/\alpha}y$ for $y \in \bR^m$. By the scaling property of the measure $\mu_{s_1,\ldots,s_m}$ given in Lemma \ref{mu-properties}.b),
$$\mu_{s_1,\ldots,s_m}(h_c^{-1}(A))=\mu_{s_1,\ldots,s_m}(c^{-1/\alpha}A)=c\mu_{s_1,\ldots,s_m}(A),$$
for any Borel set $A\subset \bR^m$. Therefore, the characteristic function of $V$ is
$$E(e^{iu\cdot V})=\exp\left\{t\int_{\bR^m} (e^{iu \cdot y}-1) (\mu_{s_1, \ldots,s_m}\circ h_{c}^{-1})(dy) \right\}=\exp\left\{ct\int_{\bR^m} (e^{iu \cdot y}-1) \mu_{s_1, \ldots,s_m}(dy) \right\}$$
for any $u \in \bR^m$, which is the same as the characteristic function of $U$. Hence $U \stackrel{d}{=}V$. $\Box$

\vspace{2mm}

{\bf The following result proves Theorem \ref{main1}.b) in the case $\alpha<1$.}

\begin{theorem}
\label{main-conv-th1}
If $\alpha<1$ and $\{Z(t)\}_{t \geq 0}$ is the process defined in Lemma \ref{Zk-Z-D}, then there exists a collection $\{\widetilde{Z}(t)\}_{t \geq 0}$ of random elements in $\bD$, such that $P(Z(t)=\widetilde{Z}(t))=1$ for all $t\geq 0$, and for any $T>0$,
\begin{equation}
\label{Zk-converges-Z}
\sup_{t\leq T}\|Z^{(\e_{k})}(t)-\widetilde{Z}(t)\| \to 0 \quad a.s. \quad as \ k \to \infty.
\end{equation}
Moreover, the map $t \mapsto \widetilde{Z}(t)$ is in $\bD_u([0,\infty);\bD)$ a.s.
\end{theorem}

\noindent {\bf Proof:} For any $T>0$, we denote by $\bD_u([0,T];\bD)$ the set of functions $x:[0,T] \to \bD$ which are right-continuous and have left-limits with respect to the norm $\|\cdot\|$ on $\bD$. Note that $\bD_u([0,T];\bD)$ is a Banach space with respect to the super-uniform norm $\|\cdot\|_{T,\bD}$.

Using the same idea as in the proof of Theorem 5.4 of \cite{resnick07}, we will show that there exists an event $\widetilde{\Omega}$ of probability 1, on which we can say that for any $T>0$,
\begin{equation}
\label{Z-e-Cauchy}
\{Z^{(\e_k)}(\cdot)\}_{k\geq 1} \quad \mbox{is a Cauchy sequence in $\bD_u([0,T];\bD)$},
\end{equation}
where $\bD_u([0,T];\bD)$ is equipped with the norm $\|\cdot\|_{T,\bD}$. We denote by $\{\widetilde{Z}(t)\}_{t \in [0,T]}$ the limit of this sequence in $\bD_u([0,T];\bD)$ (on the event $\widetilde{\Omega}$). Relation \eqref{Zk-converges-Z} then holds by definition. Since $T>0$ is arbitrary, $\widetilde{Z}(\omega,t)$ is a well-defined element in $\bD$ for any $t\geq 0$ and $\omega \in \widetilde{\Omega}$. For $\omega \not \in \widetilde{\Omega}$, we let $\widetilde{Z}(\omega,t)=y_0$ for any $t\geq 0$, where $y_0\in \bD$ is arbitrary. For any $\omega \in \Omega$ and $t \geq 0$, $Z(\omega,t)\in \bD$ and we denote $\widetilde{Z}(\omega,t,s):=\widetilde{Z}(\omega,t)(s)$ for any $s \in [0,1]$. Clearly, $\widetilde{Z}(t,s)$ is $\cF$-measurable for any $s \in [0,1]$, being the a.s. limit of the sequence $\{Z^{(\e_k)}(t,s)\}_{k \geq 1}$
This proves that $\widetilde{Z}(t)$ is a random element in $\bD$, for any $t \geq 0$.

By Lemma \ref{lemmaA} (with $S=\bD$ equipped with the uniform norm),  the map $t \mapsto \widetilde{Z}(t)$ lies in $\bD_u([0,\infty);\bD)$ (on the event $\widetilde{\Omega}$). From relation \eqref{Zk-converges-Z} and Lemma \ref{Zk-Z-D}, we infer that $\|Z(t)-\widetilde{Z}(t)\|=0$ a.s. for any $t>0$.

It remains to prove \eqref{Z-e-Cauchy}. For this, it suffices to prove that for any $\delta>0$,
\begin{equation}
\label{Z-e-Cauchy2}
\lim_{K \to \infty} \lim_{L \to \infty}P(\max_{K< k \leq L}\|Z^{(\e_k)}-Z^{(\e_K)}\|_{T,\bD}>\delta)=0.
\end{equation}

Let $\delta>0$ be arbitrary. For any $K<k\leq L$, $t > 0$ and $s \in [0,1]$,
$$Z^{(\e_k)}(t,s)-Z^{(\e_K)}(t,s)=\int_{[0,T] \times (\e_k,e_K] \times \SD}rz(s)N(du,dr,dz)=\sum_{T_i \leq t}R_i W_i(s)1_{\{\e_k<R_i\leq \e_K\}},$$
and hence
$$\|Z^{(\e_k)}(t)-Z^{(\e_K)}(t)\|\leq \sum_{T_i \leq t}R_i 1_{\{\e_k<R_i\leq \e_K\}}=\int_{[0,t] \times (\e_k,\e_K] \times \SD}r N(du,dr,dz).$$
Taking the supremum over $t \in [0,T]$ followed by the maximum over $k$ with $K<k\leq L$, we obtain:
$$\max_{K<k\leq L}\|Z^{(\e_k)}-Z^{(\e_K)}\|_{T,\bD} \leq \int_{[0,T] \times (\e_L,\e_K] \times \SD}r N(du,dr,dz).$$
By Markov's inequality,
\begin{align*}
& P(\max_{K<k\leq L}\|Z^{(\e_k)}-Z^{(\e_K)}\|_{T,\bD}>\delta) \leq \frac{1}{\delta}\,E\left( \int_{[0,T] \times (\e_L,\e_K] \times \SD}r N(du,dr,dz)\right)\\
&\quad \quad \quad  = \frac{T}{\delta}\int_{(\e_L,\e_K] \times \SD}r \overline{\nu}(dr,dz)=\frac{T}{\delta}\int_{(\e_L,\e_K]}r \nu_{\alpha}(dr) \to 0 \quad \mbox{as} \quad K,L \to \infty,
\end{align*}
using the fact that $\int_{(\e_L,1]} r \nu_{\alpha}(dr) \to \int_{0}^1 r \nu_{\alpha}(dr)<\infty$, as $L \to \infty$. This proves \eqref{Z-e-Cauchy2}. $\Box$

\subsection{Construction in the case $\alpha>1$}
\label{subsection-alpha2}

In this subsection, we give the proof of Theorem \ref{main1} in the case $\alpha>1$. In particular, property \eqref{barZk-converges-Z} below will be used in the proof of approximation result (Theorem \ref{main2}.b)).

In this case, for any $\e>0$, $E[Z^{(\e)}(t,s)]=ct\varphi(s)\int_{\e}^{\infty}r \nu_{\alpha}(dr)$ is finite, and we denote
$$\overline{Z}^{(\e)}(t,s)=Z^{(\e)}(t,s)-E[Z^{(\e)}(t,s)],$$
where $Z^{(\e)}(t,s)$ is given by \eqref{def-Z-e}. By \eqref{def-Z-e2}, it follows that
\begin{equation}
\label{def-bar-Ze}
\overline{Z}^{(\e_k)}(t,s)=\sum_{j=0}^{k}\Big(Z_j(t,s)-E\big(Z_j(t,s)\big)\Big).
\end{equation}

\begin{remark}
\label{existence-Q}
{\rm
For any probability measure $Q$ on $(\bD,\cD)$, there exists a c\`adl\`ag process $\{Y(s)\}_{s \in [0,1]}$, defined on a probability space $(\Omega',\cF',P')$, whose law under $P'$ is $Q$. This is simply because we may take $(\Omega',\cF',P')=(\bD,\cD,Q)$ and $Y(s)=\pi_s$ for all $s \in [0,1]$. This fact will be used in the proof of Lemma \ref{Zk-Z-D2} below.
}
\end{remark}

The next result is the analogue of Lemma \ref{Zk-Z-D} for the case $\alpha>1$.
The crucial elements of its proof are: {\em (i)} tightness of the sequence $\{\overline{Z}^{(\e_k)}(t)\}_{k \geq 1}$ in $\bD$, proved 
in \cite{roueff-soulier15}; and {\em (ii)} the improved version of It\^o-Nisio theorem for random elements in $\bD$, given in \cite{basse-rosinski13}. (The original version of It\^o-Nisio theorem in $\bD$ can be found in \cite{kallenberg74}.)
Recall that in the case $\alpha>1$, the process $\overline{Z}(t)=\{\overline{Z}(t,s)\}_{s \in [0,1]}$ is given by \eqref{def-Z2}.



\begin{lemma}
\label{Zk-Z-D2}
For any $t\geq 0$, there exists a random element $Z(t)=\{Z(t,s)\}_{s \in [0,1]}$ in $\bD$ such that $P(Z(t,s)=\overline{Z}(t,s))=1$ for all $s \in [0,1]$, and
\begin{equation}
\label{conv-barZ}
\lim_{k \to \infty}\|\overline{Z}^{(\e_k)}(t)-Z(t)\|=0 \quad a.s.
\end{equation}
In particular, $E\big(Z(t,s)\big)=0$ for all $s \in [0,1]$ and $t>0$.
\end{lemma}

\noindent {\bf Proof:} For $t=0$, we define $Z(0,s)=0$ for all $s \in [0,1]$. We will assume for simplicity that $t=1$, the case of arbitrary $t>0$ being similar. To simplify the notation, in this proof we denote $\overline{Z}^{(\e_k)}=\{\overline{Z}^{(\e_k)}(s)=\overline{Z}^{(\e_k)}(1,s)\}_{s \in [0,1]}$ and $\overline{Z}=\{\overline{Z}(s)=\overline{Z}(1,s)\}_{s \in [0,1]}$.

From the last part of the proof of Theorem 2.12 of \cite{roueff-soulier15}, we know that $(\overline{Z}^{(\e_k)})_{k \geq 1}$ is tight in $(\bD,J_1)$.
By Prohorov's theorem, $(\overline{Z}^{(\e_k)})_{k \geq 1}$ is relatively compact in $(\bD,J_1)$.
Hence, there exists a subsequence $N' \subset  \bZ_{+}$ and
a probability measure $Q$ on $(\bD,\cD)$ such that $P \circ (\overline{Z}^{(\e_{k})})^{-1} \stackrel{w}{\to} Q$ as $k \to \infty,k \in N'$. By Remark \ref{existence-Q}, let $Y$ be a random element in $\bD$ with law $Q$,  defined on a probability space $(\Omega',\cF',P')$. Then, $\overline{Z}^{(\e_{k})} \stackrel{d}{\to} Y$ in $(\bD,J_1)$ as $k \to \infty,k \in N'$, which
implies that
\begin{equation}
\label{fin-dim-conv}
(\overline{Z}^{(\e_{k})}(s_1), \ldots,\overline{Z}^{(\e_{k})}(s_m)) \stackrel{d}{\to} (Y(s_1),\ldots,Y(s_m)),
\end{equation}
as $k \to \infty,k \in N'$, for any $s_1,\ldots,s_m \in T$, where $T=\{s\in (0,1); P'(s \in {\rm Disc}(Y))=0\} \cup \{0,1\}$ is dense in $[0,1]$ (see p.124 of \cite{billingsley68}).
By  \eqref{def-Z2} and \eqref{def-bar-Ze},
\begin{equation}
\label{barZs-limit}
\overline{Z}(s)=\lim_{k \to \infty}\overline{Z}^{(\e_k)}(s) \quad \mbox{a.s.} \quad \mbox{for any} \ s \in [0,1].
\end{equation}
By \eqref{fin-dim-conv} and the uniqueness of the limit, it follows that for any $s_1,\ldots,s_m \in T$,
$$(\overline{Z}(s_1),\ldots, \overline{Z}(s_m))\stackrel{d}{=}(Y(s_1),\ldots,Y(s_m)).$$

Consider now another subsequence $N'' \subset \bZ_{+}$ such that $P \circ (\overline{Z}^{(\e_{k})})^{-1} \stackrel{w}{\to} Q'$ as $k \to \infty,k \in N''$, for a probability measure $Q'$ on $(\bD,\cD)$. Let $Y'$ be a random element in $\bD$ with law $Q'$, defined on a probability space $(\Omega'',\cF'',P'')$. Let $T'=\{s\in (0,1); P''(s \in {\rm Disc}(Y))=0\} \cup \{0,1\}$. The same argument as above shows that for any $s_1,\ldots,s_m \in T'$
$$(\overline{Z}(s_1),\ldots, \overline{Z}(s_m))\stackrel{d}{=}(Y'(s_1),\ldots,Y'(s_m)).$$
Hence,
$(Y(s_1),\ldots,Y(s_m))\stackrel{d}{=}(Y'(s_1),\ldots,Y'(s_m))$ for any $s_1, \ldots, s_m \in T \cap T'$. Since $T \cap T'$ is dense in $[0,1]$ and contains $1$,
by Theorem 12.5 of \cite{billingsley99}, we conclude that $Q=Q'$. This shows that any subsequence of
$\{P \circ (\overline{Z}^{(\e_{k})})^{-1}\}_k$ which converges weakly, in fact converges weakly to $Q$.
Therefore, $P \circ (\overline{Z}^{(\e_{k})})^{-1} \stackrel{w}{\to} Q$ as $k \to \infty$, and relation \eqref{fin-dim-conv} holds as $k \to \infty$ (not only along the subsequence $N'$).


Note that $\overline{Z}^{(\e_k)}(s)=\sum_{j=0}^{k}\big(Z_j(1,s)-E(Z_j(1,s))\big)$  and $\{ X_j=Z_j(1,\cdot)-E(Z_j(1,\cdot))\}_{j \geq 0}$ are random elements in $\bD$ (by Lemma \ref{Zj-well-def}), which are independent and have mean zero.
The existence of a c\`adl\`ag process $\{Z(s)\}_{s \in [0,1]}$ such that
$\lim_{k \to \infty}\|\overline{Z}^{(\e_{k})}-Z\|=0$ a.s. will follow by Theorem 2.1.(iii) of \cite{basse-rosinski13}. Relation (2.1) of \cite{basse-rosinski13} holds, due to \eqref{fin-dim-conv}. We only have to prove that
$\{|Y(s)|\}_{s \in [0,1]}$ is uniformly integrable, which is equivalent to
$\{|\overline{Z}(s)|\}_{s \in [0,1]}$ being uniformly integrable. This will follow from the fact that:
\begin{equation}
\label{finite-p-moment}
\sup_{s \in [0,1]} E|\overline{Z}(s)|^p<\infty \quad \mbox{for any} \ 1<p<\alpha.
\end{equation}

To prove \eqref{finite-p-moment}, recall from Proposition \ref{Z-stable} that $\overline{Z}(s)$ has a $S_{\alpha}(\sigma_s,\beta_s,0)$-distribution. By Property 1.2.17 of \cite{ST94}, $E|\overline{Z}(s)|^p=\sigma_s^p (c_{\alpha,\beta_s}(p))^p$, where
\begin{align*}
(c_{\alpha,\beta_s}(p))^p&=c_p \left(1+\beta_s^2 \tan^2 \frac{\alpha \pi}{2} \right)^{p/2\alpha}
\cos\left(\frac{p}{\alpha} \arctan \Big(\beta_s \tan \frac{\alpha \pi}{2} \Big) \right)\\
&\leq c_p \left(1+ \tan^2 \frac{\alpha \pi}{2} \right)^{p/2\alpha} \quad \mbox{for all} \quad s\in [0,1],
\end{align*}
and $c_p>0$ is a constant depending only on $p$. (The form of the constant $c_{\alpha,\beta}(p)$ plays an important roles in the argument above. This constant was computed in \cite{hardin84}.)
Note that for any $s \in [0,1]$,
\begin{align*}
\sigma_s &=C_{\alpha}^{-1}(c_s^{+}+c_s^{-})=C_{\alpha}^{-1}\mu_s(\{y \in \bR;|y|>1\})=C_{\alpha}^{-1} \overline{\nu}(\{(r,z) \in (0,\infty) \times \SD; r|z(s)|>1\}) \\
& \leq C_{\alpha}^{-1}\overline{\nu}((1,\infty) \times \SD)=C_{\alpha}^{-1} c \nu_{\alpha}((1,\infty))<\infty,
\end{align*}
where for the last equality we used definition \eqref{Levy-measure} of $\overline{\nu}$.
Relation \eqref{finite-p-moment} follows. $\Box$

\vspace{3mm}

{\bf The following result proves Theorem \ref{main1}.a) in the case $\alpha>1$.}

\begin{theorem}
\label{existence-Levy-process2}
If $\alpha \in (1,2)$, the process $\{Z(t)\}_{t \geq 0}$ defined in Lemma \ref{Zk-Z-D2} is a $\bD$-valued $\alpha$-stable L\'evy motion (corresponding to $\nu$). This process is $(1/\alpha)$-self-similar, i.e. it satisfies \eqref{self-similar}.
Moreover, for any $t\geq 0$ and for any monotone sequence $(t_k)_{k \geq 0}$ with $t_{k}\downarrow t$,
\begin{equation}
\label{conv-Z-tk}
\lim_{k \to \infty}\|Z(t_k)-Z(t)\| = 0 \quad a.s.
\end{equation}
\end{theorem}

\noindent {\bf Proof:} The first two sentences are proved exactly as in the case $\alpha<1$, with obvious modifications in the form of the characteristic functions, due to centering. We only have to prove the last sentence. For this, we apply again Theorem 2.1.(iii) of \cite{basse-rosinski13} with $E=\bR$.

For any $i \geq 1$, let $X_i=Z(t_{i-1})-Z(t_i)$. By property {\em (ii)} in Definition \ref{def-Levy}, $(X_i)_{i \geq 1}$ are independent random elements in $\bD$ (with zero mean). Let $S_k=\sum_{i=1}^{k}X_i=Z(t_0)-Z(t_k)$ for all $k \geq 1$, and $Y=Z(t_0)-Z(t)$. We first show that for any $s_1, \ldots,s_m \in [0,1]$,
$$(S_k(s_1),\ldots,S_k(s_m)) \stackrel{d}{\to} (Y(s_1),\ldots,Y(s_m)) \quad \mbox{as} \ k \to \infty.$$
To see this, note that $(S_k(s_1),\ldots,S_k(s_m))
\stackrel{d}{=}(Z(t_0-t_k,s_1),\ldots,
Z(t_0-t_k,s_m))$ by property {\em (iii)} in Definition \ref{def-Levy})
(stationarity of the increments). It is now clear that we have the following convergence the characteristic functions: for any $u=(u_1,\ldots,u_m)  \in \bR^m$,
\begin{align*}
& E(e^{iu_1 S_k(s_1)+\ldots+iu_m S_k(s_m)})=\exp\left\{(t_0-t_k)\int_{\bR^m} (e^{iu \cdot y}-1-iu\cdot y)\mu_{s_1,\ldots,s_m}(dy) \right\},\\
& \quad  \quad \quad \to E(e^{iu_1 Y(s_1)+\ldots+iu_m Y(s_m)})=\exp\left\{(t_0-t)\int_{\bR^m} (e^{iu \cdot y}-1-iu\cdot y)\mu_{s_1,\ldots,s_m}(dy) \right\},
\end{align*}
as $k \to \infty$. It remains to show that $\{|Y(s)|\}_{s \in [0,1]}$ is uniformly integrable, which is equivalent to saying that $\{|Z(t_0-t,s)|\}_{s \in [0,1]}$ is uniformly integrable, by the stationarity of the increments. By the self-similarity of $\{Z(t)\}_{t \geq 0}$, $Z(t_0-t,s)\stackrel{d}{=} (t_0-t)^{1/\alpha}Z(1,s)$ for all $s \in [0,1]$. Using \eqref{finite-p-moment} and the fact that $Z(1,s)=\overline{Z}(1,s)$ a.s. for any $s \in [0,1]$, it follows that for any $1<p<\alpha$,
$$\sup_{s \in [0,1]} E|Z(t_0-t,s)|^p=(t_0-t)^{p/\alpha}\sup_{s \in [0,1]}E|Z(1,s)|^p<\infty.$$
(Recall that in \eqref{finite-p-moment} we used the notation $\overline{Z}(s)=\overline{Z}(1,s)$.) Hence, $\{|Z(t_0-t,s)|\}_{s \in [0,1]}$ is uniformly integrable. By Theorem 2.1.(iii) of \cite{basse-rosinski13}, it follows that
$S_k \to Z(t_0)-Z(t)$ a.s. in $(\bD,\|\cdot\|)$, as $k \to \infty$, which is the same as
$Z(t_k) \to Z(t)$ a.s. in $(\bD,\|\cdot\|)$, as $k \to \infty$. $\Box$

\vspace{3mm}

The following preliminary result will be used in the proof of tightness of $(\oZ^{(\e_k)})_{k \geq 1}$.
\begin{lemma}
\label{bound-EZ}
For any $\e>0$ and $T>0$,
$$E\|Z^{(\e)}\|_{T,\bD} \leq T c \frac{\alpha}{\alpha-1} \e^{1-\alpha}.$$
\end{lemma}

\noindent {\bf Proof:} By definition, for any $t\in [0,T]$ and $s \in [0,1]$, we have
$$|Z^{(\e)}(t,s)| \leq \int_{[0,t] \times (\e,\infty) \times \SD} r|z(s)|N(du,dr,dz) \leq \int_{[0,T] \times (\e,\infty) \times \SD}r N(du,dr,dz)=:Y.$$
Hence $\|Z^{(\e)}\|_{T,\bD} \leq Y$ and
$E\|Z^{(\e)}\|_{T,\bD} \leq E(Y)=T \int_{(\e,\infty) \times \SD} r \overline{\nu}(dr,dz)=T c \frac{\alpha}{\alpha-1} \e^{1-\alpha}$. $\Box$

\vspace{3mm}

The next result plays a crucial role in the proof of Theorem \ref{main1}.b) in the case $\alpha>1$.
Its proof uses some results related to sums of i.i.d. regularly varying random elements in $\bD$, which are given in Section \ref{subsection-approx2} below.

\begin{theorem}
\label{tightness-th}
If Assumption B holds, then $(\oZ^{(\e_k)})_{k \geq 1}$ 
 is tight in $\bD([0,\infty);\bD)$.
\end{theorem}

\noindent {\bf Proof:} It is enough to prove that $(\oZ^{(\e_k)})_{k \geq 1}$ is tight in $\bD([0,T];\bD)$ for any $T>0$. Without loss of generality, we assume that $T=1$. Let $P_k$ be the law of $\oZ^{(\e_k)}$. We verify that $(P_k)_{k \geq 1}$ satisfies conditions {\em (i)}-{\em (iii)} of Theorem \ref{tight-th}. To prove this, we argue as in the last part of the proof of Theorem 2.12 of \cite{roueff-soulier15}.

For condition {\em (i)}, it suffices to show that the following two relations hold:
\begin{align}
\label{cond-i1}
& \lim_{A \to \infty}P(\|\oZ^{(\e_0)}\|_{\bD}>A)=0 \quad \mbox{for all} \ \e_0>0\\
\label{cond-i2}
& \lim_{\e_0 \downarrow 0}\sup_{0<\e<\e_0} P(\|\oZ^{(\e)}-\oZ^{(\e_0)}\|_{\bD}>\eta)=0 \quad \mbox{for all} \ \eta>0.
\end{align}
To see this, let $\eta>0$ and $\rho>0$ be arbitrary. By
\eqref{cond-i2} and the fact that $\e_k \downarrow 0$, there exist $\e_0^*\in (0,1)$ and $k_0$ such that
$P(\|\oZ^{(\e_k)}-\oZ^{(\e_{0}^*)}\|_{\bD}>\eta)<\rho/2$ for any $k\geq k_0$. By \eqref{cond-i1}, there exists $A_0>0$ such that $P(\|\oZ^{(\e_0^*)}\|_{\bD}>A_0)<\rho/2$. Let $a_0=\eta+A_0$. Then, for all $k \geq k_0$,
$$P(\|\oZ^{(\e_k)}\|_{\bD}>a_0) \leq P(\|\oZ^{(\e_k)}-\oZ^{(\e_0^*)}\|_{\bD}>\eta)+
P(\|\oZ^{(\e_0^*)}\|_{\bD}>A_0)<\rho.$$
This proves that condition {\em (i)} holds.

To prove \eqref{cond-i1}, let $\e_0>0$ be arbitrary. For any $A>2\|E(Z^{(\e_0)})\|_{\bD}$,
\begin{align*}
P(\|\oZ^{(\e_0)}\|_{\bD}>A)\leq
P(\|Z^{(\e_0)}\|_{\bD}>A/2)\leq \frac{2}{A}\|E(Z^{(\e_0)})\|_{\bD}\leq \frac{2}{A}T c \frac{\alpha}{\alpha-1} \e_0^{1-\alpha},
\end{align*}
using Markov inequality and Lemma \ref{bound-EZ}. Relation \eqref{cond-i1} follows letting $A \to \infty$.

To prove \eqref{cond-i2}, we use an indirect argument. Consider a sequence $(X_i)_{i \geq 1}$ of i.i.d. regularly varying elements in $\bD$ (as given by Definition \ref{def-RV}) with limiting measure $\overline{\nu}$ given by  \eqref{Levy-measure}.
Let $S_n^{(\e)}$ be given by relation \eqref{truncated-sum} below. Similarly to Theorem \ref{centered-tr-sum-conv} below (which is based on the fact that the probability measure $\Gamma_1$ satisfies Assumptions B), it can be proved that for any $0<\e<\e_0$,
\begin{equation}
\label{conv-Sn}
S_n^{(\e)}-S_n^{(\e_0)}-E(S_n^{(\e)}-S_n^{(\e_0)}) \stackrel{d}{\to} \oZ^{(\e)}-\oZ^{(\e_0)} \quad \mbox{in} \quad \bD([0,1];\bD),
\end{equation}
where $\bD([0,1];\bD)$ is equipped with distance $d_{\bD}$. For any $t>0$ and $s \in [0,1]$, we define
$$S_n^{<\e}(t,s)=\frac{1}{a_n}\sum_{i=1}^{[nt]}X_i(s) 1_{\{\|X_i\| \leq a_n \e\}}.$$
Then $S_n^{(\e)}=S_n-S_n^{<\e}$. Hence,
$S_n^{(\e)}-S_n^{(\e_0)}=S_n^{<\e_0}-S_n^{<\e}$ and relation \eqref{conv-Sn} becomes:
$$S_n^{<\e_0}-S_n^{<\e}-E(S_n^{<\e_0}-S_n^{<\e}) \stackrel{d}{\to} \oZ^{(\e)}-\oZ^{(\e_0)} \quad \mbox{in} \quad \bD([0,1];\bD).$$
Since $\|\cdot\|_{\bD}$ is $d_{\bD}$-continuous (see Lemma \ref{dD-cont}), by the continuous mapping theorem, we have: $\|S_n^{<\e_0}-S_n^{<\e}-E(S_n^{<\e_0}-S_n^{<\e})\|_{\bD} \stackrel{d}{\to} \|\oZ^{(\e)}-\oZ^{(\e_0)}\|_{\bD}$ as $n \to \infty$. Let $\eta>0$ be arbitrary. By Portmanteau theorem,
\begin{align*}
& P( \|\oZ^{(\e)}-\oZ^{(\e_0)}\|_{\bD}>\eta) \leq \liminf_{n\to \infty}P(\|S_n^{<\e_0}-S_n^{<\e}-E(S_n^{<\e_0}-S_n^{<\e})\|_{\bD}>\eta)\\
& \quad \quad \quad \leq \limsup_{n \to \infty}
P(\|S_n^{<\e_0}-E(S_n^{<\e_0})\|_{\bD}>\eta/2)+
P(\|S_n^{<\e}-E(S_n^{<\e})\|_{\bD}>\eta/2).
\end{align*}
We take the supremum over all $\e \in (0,\e_0)$, followed by the limit as $\e_0 \downarrow 0$. We obtain that $\lim_{\e_0 \downarrow 0} \sup_{0<\e<\e_0}P( \|\oZ^{(\e)}-\oZ^{(\e_0)}\|_{\bD}>\eta)$ is less than
\begin{align*}
\lim_{\e_0 \downarrow 0} \limsup_{n \to \infty}P(\|S_n^{<\e_0}-E(S_n^{<\e_0})\|_{\bD}>\eta/2)+
\lim_{\e_0 \downarrow 0}  \sup_{0<\e<\e_0}\limsup_{n \to \infty}P(\|S_n^{<\e}-E(S_n^{<\e})\|_{\bD}>\eta/2).
\end{align*}
Since $S_n^{<\e}=S_n-S_n^{(\e)}$, both these terms are zero, by relation \eqref{AN2} below (with $T=1$). This concludes the proof of \eqref{cond-i2}.

We prove that $(P_k)_{k \geq 1}$ satisfies condition {\em (ii)} of Theorem \ref{tight-th}. Let $\eta>0$ and $\rho>0$ be arbitrary. It suffices to show that there exist $\delta \in (0,1)$ and $\e_0>0$ such that for all $\e \in (0,\e_0)$,
\begin{equation}
\label{3cond-partii}
\left\{
\begin{array}{ll}
(a) & P(w''(\oZ^{(\e)}(t,\delta) > \eta \ \mbox{for some} \ t \in [0,1])<\rho \\
(b) & P(|\oZ^{(\e)}(t,\delta)-\oZ^{(\e)}(t,0)| > \eta  \ \mbox{for some} \ t \in [0,1])<\rho \\
(c) & P(|\oZ^{(\e)}(t,1-)-\oZ^{(\e)}(t,1-\delta)| > \eta  \ \mbox{for some} \ t \in [0,1])<\rho.
\end{array}
\right.
\end{equation}
By \eqref{cond-i2}, there exists $\e_0>0$ such that
\begin{equation}
\label{bound-oZ}
P(\|\oZ^{(\e)}-\oZ^{(\e_0)}\|_{\bD}>\eta/4)<\rho/2 \quad \mbox{for all} \quad \e \in (0,\e_0).
\end{equation}
Since $\bD([0,1];\bD)$ endowed with $d_{\bD}^0$ is separable and complete (see Theorem \ref{bD-CSMS}), by Theorem 1.3 of \cite{billingsley99}, the single probability measure $P \circ (\oZ^{(\e_0)})^{-1}$ is tight. Hence, by condition {\em (ii)} of Theorem \ref{tight-th}, there exists $\delta \in (0,1)$ such that
\begin{align}
\label{ii-a}
& P(w''(\oZ^{(\e_0)}(t,\delta) > \eta/2 \ \mbox{for some} \ t \in [0,1])<\rho/2 \\
\label{ii-b}
& P(|\oZ^{(\e_0)}(t,\delta)-\oZ^{(\e_0)}(t,0)| > \eta/2  \ \mbox{for some} \ t \in [0,1])<\rho/2 \\
\label{ii-c}
& P(|\oZ^{(\e_0)}(t,1-)-\oZ^{(\e_0)}(t,1-\delta)| > \eta/2  \ \mbox{for some} \ t \in [0,1])<\rho/2.
\end{align}
Using the fact that
$$w''(x+y,\delta) \leq w''(x,\delta)+2\|y\| \quad \mbox{for all} \quad x,y\in \bD,$$
we infer that
$w''(\oZ^{(\e)}(t),\delta) \leq w''(\oZ^{(\e_0)}(t),\delta)+2\|\oZ^{(\e)}-\oZ^{(\e_0)}\|_{\bD}$, and hence $P(w''(\oZ^{(\e)}(t),\delta)> \eta \ \mbox{for some} \ t\in[0,1])$ is smaller than
$$P(w''(\oZ^{(\e_0)}(t),\delta)> \eta/2 \ \mbox{for some} \ t\in[0,1])+
P(\|\oZ^{(\e)}-\oZ^{(\e_0)}\|_{\bD}>\eta/4).$$
Part (a) of \eqref{3cond-partii} follows from \eqref{bound-oZ} and \eqref{ii-a}. Similarly, part (b) of \eqref{3cond-partii} follows from \eqref{bound-oZ} and \eqref{ii-b}, using the fact that
$$|\oZ^{(\e)}(t,\delta)-\oZ^{(\e)}(t,0)|\leq |\oZ^{(\e_0)}(t,\delta)-\oZ^{(\e_0)}(t,0)|+2\|\oZ^{(\e)}-\oZ^{(\e_0)}\|_{\bD},$$
whereas part (c) of  \eqref{3cond-partii} follows from \eqref{bound-oZ} and \eqref{ii-c}, since
$$|\oZ^{(\e)}(t,1-)-\oZ^{(\e)}(t,1-\delta)|\leq |\oZ^{(\e_0)}(t,1-)-\oZ^{(\e_0)}(t,1-\delta)|
+2\|\oZ^{(\e)}-\oZ^{(\e_0)}\|_{\bD}.$$

It remains to prove that $(P_k)_{k \geq 1}$ satisfies condition {\em (iii)} of Theorem \ref{tight-th}. Let $\eta>0$ and $\rho>0$ be arbitrary. Note that $\oZ^{(\e)}(0)=0$.
We will show that there exist $\delta \in (0,1)$ and $\e_0>0$ such that for all $\e \in (0,\e_0)$,
\begin{equation}
\label{3cond-partiii}
\left\{
\begin{array}{ll}
(a) & P(w_{\bD}''(\oZ^{(\e)},\delta)> \eta)<\rho \\
(b) & P(\|\oZ^{(\e)}(\delta)\|> \eta) <\rho \\
(c) & P\big(d_{J_1}^0 \big(\oZ^{(\e)}(1-),\oZ^{(\e)}(1-\delta)\big) > 3\eta/2 \big)<\rho.
\end{array}
\right.
\end{equation}
Let $\e_0$ be such that \eqref{bound-oZ} holds. Using again the fact that $P \circ (\oZ^{(\e_0)})^{-1}$ is tight, but invoking this time condition {\em (iii)} of Theorem \ref{tight-th}, we infer that there exists $\delta \in (0,1)$ such that
\begin{align}
\label{iii-a}
& P(w_{\bD}''(\oZ^{(\e_0)},\delta)> \eta/2)<\rho/2 \\
\label{iii-b}
& P(\|\oZ^{(\e_0)}(\delta)\|> \eta/2) <\rho/2 \\
\label{iii-c}
& P\big(d_{J_1}^0\big(\oZ^{(\e_0)}(1-)-\oZ^{(\e_0)}(1-\delta)\big) > \eta/2\big)<\rho/2.
\end{align}
By Lemma \ref{dD-sum}, 
$P(w_{\bD}''(\oZ^{(\e)},\delta)>\eta) \leq P(w_{\bD}''(\oZ^{(\e_0)},\delta)>\eta/2) +\
P(2\|\oZ^{(\e)}-\oZ^{(\e_0)}\|_{\bD}>\eta/2)<\rho$.
Part (a) of \eqref{3cond-partiii} follows using \eqref{iii-a} and \eqref{bound-oZ}.
Part (b) of \eqref{3cond-partiii} follows using \eqref{iii-b} and \eqref{bound-oZ}, since
$\|\oZ^{(\e)}(\delta)\| \leq \|\oZ^{(\e_0)}(\delta)\|+\|\oZ^{(\e)}-\oZ^{(\e_0)}\|_{\bD}$.
To see that part (c) of \eqref{3cond-partiii} holds, note that by the triangular inequality,
$d_{J_1}^0 \big(\oZ^{(\e)}(1-),\oZ^{(\e)}(1-\delta)\big)$ is smaller than
$$d_{J_1}^0 \big(\oZ^{(\e)}(1-),\oZ^{(\e_0)}(1-)\big)+d_{J_1}^0 \big(\oZ^{(\e_0)}(1-),\oZ^{(\e_0)}(1-\delta)\big)+d_{J_1}^0 \big(\oZ^{(\e_0)}(1-\delta),\oZ^{(\e)}(1-\delta)\big).$$
We treat separately these three terms. For the second term, we use \eqref{iii-c}. For the last term, we use \eqref{bound-oZ}, since this term is bounded by $\|\oZ^{(\e_0)}(1-\delta)-\oZ^{(\e)}(1-\delta)\|$ which is smaller than $\|\oZ^{(\e_0)}-\oZ^{(\e)}\|_{\bD}$. For the first term, we also use \eqref{iii-c}, since this term is bounded by
$\|\oZ^{(\e)}(1-)-\oZ^{(\e_0)}(1-)\|$ which is smaller than $\|\oZ^{(\e_0)}-\oZ^{(\e)}\|_{\bD}$. To see this, note that by Remark \ref{Z-e-in-Du},
$\oZ^{(\e)}(1-)=\lim_{\delta \to 0}\oZ^{(\e)}(1-\delta)$ in $(\bD,\|\cdot\|)$ and $\oZ^{(\e_0)}(1-)=\lim_{\delta \to 0}\oZ^{(\e_0)}(1-\delta)$ in $(\bD,\|\cdot\|)$, and hence
$$\|\oZ^{(\e)}(1-)-\oZ^{(\e_0)}(1-)\|=\lim_{\delta \to 0}\|\oZ^{(\e)}(1-\delta)-\oZ^{(\e_0)}(1-\delta) \| \leq \|\oZ^{(\e)}-\oZ^{(\e_0)}\|_{\bD}.$$
$\Box$

\vspace{3mm}

{\bf The following result proves Theorem \ref{main1}.b) in the case $\alpha>1$.}

\begin{theorem}
\label{main-conv-th2}
If $\alpha \in (1,2)$ and Assumption B holds, then there exists a collection $\{\widetilde{Z}(t)\}_{t \geq 0}$ of random elements in $\bD$ such that $P(Z(t)=\widetilde{Z}(t))=1$ for all $t\geq 0$, the map $t \mapsto \widetilde{Z}(t)$ is in $\bD([0,\infty);\bD)$, and
\begin{equation}
\label{barZk-converges-Z}
\overline{Z}^{(\e_{k})}(\cdot) \stackrel{d}{\to} \widetilde{Z}(\cdot) \quad in \quad \bD([0,\infty);\bD)
\end{equation}
as $k \to \infty$, $k \in N'$, for a subsequence $N' \subset \bZ_{+}$, where $\bD([0,\infty);\bD)$ is equipped with the Skorohod distance $d_{\infty,\bD}$ given by \eqref{def-d-infty}.
\end{theorem}

\noindent {\bf Proof:} {\em Step 1.} By Theorem \ref{tightness-th}, there exists a subsequence $N' \subset \bZ_{+}$ such that
\begin{equation}
\label{w-conv-ZY}
\overline{Z}^{(\e_{k})}(\cdot) \stackrel{d}{\to} Y(\cdot) \quad \mbox{in} \quad \bD([0,\infty);\bD),
\end{equation}
as $k \to \infty,k \in N'$, where $Y$ is a random element in $\bD([0,\infty);\bD)$, defined on a probability space $(\Omega',\cF',P')$. We prove that for any $t_1,\ldots,t_n \geq 0$,
\begin{equation}
\label{fdd-equality}
(Z(t_1),\ldots,Z(t_n))\stackrel{d}{=} (Y(t_1),\ldots,Y(t_n)) \quad \mbox{in} \ \bD^n.
\end{equation}
To see this, note that \eqref{w-conv-ZY} implies that
$(\overline{Z}^{(\e_k)}(t_1),\ldots,\overline{Z}^{(\e_k)}(t_n)) \stackrel{d}{\to} (Y(t_1),\ldots,Y(t_n))$ in $(\bD^n,J_1^n)$, for any $t_1, \ldots, t_n \in T_Y=T_{P' \circ Y^{-1}}$ (see \eqref{marginal-conv}). On the other hand, by \eqref{conv-barZ},
$(\overline{Z}^{(\e_k)}(t_1),\ldots,\overline{Z}^{(\e_k)}(t_n)) \stackrel{p}{\to} (Z(t_1),\ldots,Z(t_n))$ in $(\bD^n,J_1^n)$ for any $t_1, \ldots,t_n \geq 0$. By the uniqueness of the limit, \eqref{fdd-equality} holds for any $t_1,\ldots,t_n \in T_Y$. To see that \eqref{fdd-equality} holds for arbitrary $t_1,\ldots,t_n \geq 0$, we proceed by approximation. Since $T_Y$ is dense in $[0,\infty)$, for any $i=1,\ldots,n$, there exists a monotone sequence $(t_i^k)_{k} \subset T_Y$ such that $t_i^k \downarrow t_i$ as $k \to \infty$. By \eqref{conv-Z-tk},
$(Z(t_1^k),\ldots,Z(t_n^k)) \stackrel{p}{\to} (Z(t_1),\ldots,Z(t_n))$ in $(\bD^n,J_1^n)$ as $k \to \infty$. Since $Y$ has all sample paths in $\bD([0,\infty);\bD)$,
$(Y(t_1^k),\ldots,Y(t_n^k)) \to (Y(t_1),\ldots,Y(t_n))$ in $(\bD^n,J_1^n)$ as $k \to \infty$. Relation \eqref{fdd-equality} follows again by the uniqueness of the limit.

{\em Step 2.} Relation \eqref{fdd-equality} shows that processes $\{Z(t)\}_{t \geq 0}$ and $\{Y(t)\}_{t \geq 0}$ have the same finite-dimensional distributions. The process $\{Y(t)\}_{t \geq 0}$ has sample paths in $\bD([0,\infty);\bD)$, which is a Borel space (being a Polish space). By Lemma 3.24 of \cite{kallenberg02}, there exists a process $\{\widetilde{Z}(t)\}_{t \geq 0}$ defined on the same probability space $(\Omega,\cF,P)$, whose sample paths are in $\bD([0,\infty);\bD)$, such that $P(Z(t)=\widetilde{Z}(t))=1$ for all $t \geq 0$. In particular, $\{\widetilde{Z}(t)\}_{t \geq 0}$ has the same finite-dimensional distributions as $\{Z(t)\}_{t \geq 0}$, hence also as $\{Y(t)\}_{t \geq 0}$. Since finite-dimensional distributions uniquely determine the law, it follows that the random elements $\widetilde{Z}(\cdot)=\{\widetilde{Z}(t)\}_{t \geq 0}$ and $Y(\cdot)=\{Y(t)\}_{t \geq 0}$ have the same law in $\bD([0,\infty);\bD)$.
Relation \eqref{barZk-converges-Z} follows from \eqref{w-conv-ZY}. $\Box$

\section{Approximation: proof of Theorem \ref{main2}}
\label{section-approximation}

In this section, we show that the $\alpha$-stable L\'evy process with values in $\bD$ constructed in the Section \ref{section-construction} can be obtained as the limit (in distribution) of the partial sum sequence associated with i.i.d. regularly varying elements in $\bD$, with suitable normalization and centering. This result can be viewed as an extension of the stable functional central limit theorem (see e.g. Theorem 7.1 of \cite{resnick07}) to the case of random elements in $\bD$. The proof of this result uses the method of point process convergence, instead of the classical method based on finite dimensional convergence and tightness. A similar method was used in \cite{roueff-soulier15} for fixed time $t=1$. We extend the arguments of \cite{roueff-soulier15} to include the time variable $t>0$.

\subsection{Point processes on Polish spaces}
\label{section-point}

In this subsection, we review some basic concepts related to point processes on a Polish space, following \cite{D-VJ}. Similar concepts are considered in \cite{resnick87, resnick07} for point processes on an LCCB space (i.e. a locally compact space with countable basis).

Let $(E,d)$ be a Polish space (i.e. a complete separable metric space) and $\cE$ its Borel $\sigma$-field. A measure $\mu$ on $E$ is {\em boundedly finite} if $\mu(A)<\infty$ for all bounded sets $A \in \cE$. (Recall that a set $A$ is bounded if it is contained in an open ball.) 
We denote by $\widehat{M}_{+}(E)$ the set of all boundedly finite measures on $E$, and by $\widehat{M}_p(E)$ its subset consisting of point (or counting) measures, i.e. $\bZ_{+}$-valued measures, where $\bZ_{+}=\{0,1,2,\ldots,\}$. A measure $\mu \in \widehat{M}_p(E)$ can be represented as $\mu=\sum_{i \geq 1}\delta_{x_i}$ for some $(x_i)_{i \geq 1} \subset E$, where $\delta_x$ is the Dirac measure at $x$. In this case, $(x_i)_{i\geq 1}$ are called the atoms (or points) of $\mu$. A measure $\mu=\sum_{i \geq 1}\delta_{x_i} \in \widehat{M}_p(E)$ is {\em simple} if $\mu(\{x\}) \leq 1$ for all $x \in E$, i.e. $(x_i)_{i \geq 1}$ are distinct.

The set $\widehat{M}_{+}(E)$ is equipped with the topology of {\em $\widehat{w}$-convergence}: $\mu_n \stackrel{\widehat{w}}{\to}\mu$ on $E$ if  $\mu_n(A) \to \mu(A)$ for any bounded set $A \in \cE$ with $\mu(\partial A)=0$. By Proposition A.2.6.II of \cite{D-VJ}, this is equivalent to $\mu_n(f)\to \mu(f)$ for any $f \in \widehat{C}(E)$, where $\mu(f)=\int_{E}fd\mu$ and $\widehat{C}(E)$ is the set of bounded continuous functions $f:E \to \bR$ which vanish outside a bounded set. We denote by $\widehat{\cM}_{+}(E)$ and $\widehat{\cM}_{p}(E)$ the Borel $\sigma$-fields of $\widehat{M}_{+}(E)$, respectively $\widehat{M}_{p}(E)$. By Proposition 9.1.IV of \cite{D-VJ}, $\widehat{M}_{+}(E)$ and $\widehat{M}_{p}(E)$ are  Polish spaces, and $\widehat{\cM}_{+}(E)$ and $\widehat{\cM}_{p}(E)$ are generated by the functions $\widehat{M}_{+}(E) \ni  \mu \mapsto \mu(A), A \in \cE$, respectively $\widehat{M}_p(E) \ni \mu \mapsto \mu(A),A \in \cE$.

A {\em point process} on $E$ is a function $N:\Omega \to \widehat{M}_{p}(E)$ defined on a probability space $(\Omega,\cF,P)$, which is $\cF/\widehat{\cM}_{p}(E)$-measurable, i.e. $N(A):\Omega \to \bZ_{+}$ is $\cF$-measurable for any $A \in \cE$. The law $P \circ N^{-1}$ of $N$ is uniquely determined by the {\em Laplace functional} $L_N(f)=E(e^{-N(f)})$, for all measurable functions $f:E \to [0,\infty)$ with bounded support.

We say that a sequence $(N_n)_{n \geq 1}$ of point processes on $E$ converges in distribution to the point process $N$ on $E$ and we write $N_n \stackrel{d}{\to}N$ in $\widehat{M}_p(E)$, if $(P \circ N_n^{-1})_{n \geq 1}$ converges weakly to $P \circ N^{-1}$ as probability measures on $\widehat{M}_p(E)$. By Proposition 11.1.VIII of \cite{D-VJ}, this is equivalent to $L_{N_n}(f) \to L_N(f)$ for all continuous functions $f:E\to \bR$ vanishing outside a bounded set.

\begin{definition}
\label{def-PRM}
{\rm Let $\nu \in \widehat{M}_{+}(E)$ be arbitrary. A point process $N$ on $E$ is called a {\em Poisson random measure} on $E$ of intensity $\nu$, if for any bounded set $A \in \cE$, $N(A)$ has a Poisson distribution with mean $\nu(A)$, and for any bounded disjoint sets $A_1,\ldots,A_n \in \cE$, $N(A_1), \ldots, N(A_n)$ are independent.}
\end{definition}

The Laplace functional of a Poisson random measure $N$ of intensity $\nu$ on $E$ is:
\begin{equation}
\label{Laplace-PRM}
L_N(f)=\exp\left\{-\int_{E}(1-e^{-f(x)})\nu(dx) \right\},
\end{equation}
for all bounded measurable functions $f:E \to [0,\infty)$ with bounded support.

The following result plays a crucial role in this article. It is an extension of Proposition 3.21 of \cite{resnick87} to point processes on Polish spaces, with which shares the same proof (based on Laplace functionals). Recall that a {\em random element} in $E$ is a function $X:\Omega \to E$ defined on a probability space $(\Omega,\cF,P)$, which is $\cF/\cE$-measurable.

\begin{proposition}
\label{prop3-21}
Let $E$ be a Polish space and $\nu \in \widehat{M}_{+}(E)$ be arbitrary. For any $n\geq 1$, let $(X_{i,n})_{i \geq 1}$ be i.i.d. random elements in $E$ and $N_n=\sum_{i\geq 1}\delta_{(i/n,X_{i,n})}$. Let $N$ be a Poisson random measure on $[0,\infty) \times E$ of intensity $Leb \times \nu$, where Leb is the Lebesgue measure. Then $N_n \stackrel{d}{\to} N$ in $\widehat{M}_{p}([0,\infty) \times E)$ if and only if
$$nP(X_{1,n} \in \cdot) \stackrel{\widehat{w}}{\to} \nu \quad \mbox{on} \quad E.$$
\end{proposition}

We conclude this section with few words about finite measures. We denote by $M_f(E)$ the set of finite measures on $E$, equipped with the topology if weak convergence: $\mu_n \stackrel{w}{\to}\mu$ if $\mu_n(A) \to \mu(A)$ for any set $A \in \cE$ with $\mu(\partial A)=0$. Finally, we denote by $M_{p,f}(E)$ the set of finite point measures on $E$, equipped also with the topology of weak convergence.

\subsection{Continuity of summation functional}

In this subsection, we establish the continuity of the truncated summation functional defined on the set of point measures on $[0,\infty) \times \barD$. This will constitute an important step in the proof of our main result. The proofs contained in this subsection are extensions of those of \cite{roueff-soulier15} to point measures whose atoms include also a time variable.

We endow the spaces $[0,\infty) \times \barD$ and $[0,T] \times \bD$ with the product topologies, $\bD$ being equipped with Skorohod's $J_1$-topology.

For fixed $T>0$ and $\e>0$, we define $\Psi: \widehat{M}_{p}([0,\infty) \times \barD) \to M_{p,f}([0,T] \times \bD)$ by:
$$\Psi(m)=m|_{[0,T] \times (\e,\infty) \times \SD} \circ \psi^{-1}$$
where $m|_{[0,T] \times (\e,\infty) \times \SD}$ denotes the restriction of $m$ to $[0,T] \times (\e,\infty) \times \SD$, and the function
$\psi: [0,\infty) \times (\e,\infty) \times \SD \to [0,T] \times \bD$ is given by $\psi(t,r,z)=(t,rz)$. Note that $\Psi(m)$ is a finite measure since
$[0,T] \times (\e,\infty) \times \SD$ is a bounded set.

The application of the function $\Psi$ has a double effect on a measure $m$: it removes the atoms $(t_i,r_i,z_i)$ of $m$ whose second coordinate $r_i$ is less than $\e$ or is $\infty$, and transforms the remaining atoms using the ``inverse polar-coordinate'' map $(r,z) \mapsto rz$, while leaving the first coordinate $t_i$ of these atoms unchanged (provided that $t_i \leq T$). More precisely, if $m=\sum_{i\geq 1}\delta_{(t_i,r_i,z_i)} \in \widehat{M}_p([0,\infty) \times \barD)$ then $\Psi(m)=\sum_{t_i \leq T} \delta_{(t_i,r_i z_i)}1_{\{r_i \in (\e,\infty)\}}$.

For any $m\in \widehat{M}_p([0,\infty) \times \barD)$ and for any measurable function $f:[0,T] \times \bD \to [0,\infty)$,
\begin{equation}
\label{integral-Psi-m}
\int_{[0,T] \times \bD}f(t,x)\Psi(m)(dt,dx)=\int_{[0,T] \times (\e,\infty) \times \SD} f(t,rz)m(dt,dr,dz).
\end{equation}


\begin{lemma}
\label{cont-Psi}
The function $\Psi$ is continuous on the set $\cA$ of measures $m \in M_p([0,\infty) \times \barD)$ which satisfy the following two conditions: $$m([0,\infty) \times \{\e,\infty\} \times \SD)=0 \quad and \quad m(\{0,T\} \times (\e,\infty) \times \SD)=0.$$
(The function $\Psi=\Psi_{\e,T}$ and the set $\cA=\cA_{\e,T}$ depend on $\e$ and $T$. To simplify the writing, we drop the indices $\e$, $T$.)
\end{lemma}

\noindent {\bf Proof:} Let $E=[0,\infty) \times \barD$, $E'=[0,\infty) \times (\e,\infty) \times \SD$ and $E''=[0,T] \times \bD$. Since $E'$ is a bounded set, $\widehat{M}_{p}(E')=M_{p,f}(E')$. Note that $\Psi=\Psi_2 \circ \Psi_1$, where
 $\Psi_1: \widehat{M}_{p}(E) \to M_{p,f}(E')$ is the restriction
$\Psi_1(m)=m|_{E'}$ and $\Psi_2: M_{p,f}(E') \to M_{p,f}(E'')$ is given by $\Psi_2(m)=m \circ \psi^{-1}$.

Similarly to Proposition 3.3 of \cite{feigin-kratz-resnick}, it can be shown that $\Psi_1$ is continuous on $\cA$. The fact that $\Psi_2$ is continuous follows from the continuity of function $\psi$, exactly as in the proof of Proposition 5.6.(a) of \cite{resnick07}. $\Box$


\vspace{3mm}

\begin{definition}
\label{def-M-pf-star}
{\rm We denote by $M_{p,f}^{*}([0,T] \times \bD)$ the set of measures $\mu \in M_{p,f}([0,T] \times \bD)$ which have the following properties: {\em (i)} $\mu$ is simple;
{\em (ii)} $\mu(\{(t,x),(t',x')\})\leq 1$ for any $(t,x),(t',x') \in [0,T] \times \bD$ with $x \not=x'$ and ${\rm Disc}(x) \cap {\rm Disc}(x')\not=\emptyset$;
{\em (iii)} $\mu(\{t_0\} \times \bD) \leq 1$ for all $t_0 \in [0,T]$.}
\end{definition}

Alternatively, we can say that $M_{p,f}^{*}([0,T] \times \bD)$ is the set of finite point measures $\mu=\sum_{i=1}^{p}\delta_{x_i}$ on $[0,T] \times \bD$ which satisfy the following three conditions: (1) the points $(t_1,x_1),\ldots,(t_p,x_p)$ are distinct; (2) ${\rm Disc}(x_i)\cap {\rm Disc}(x_j)=\emptyset$ for all $i\not=j$; (3) no vertical line contains two points of $\mu$.

The next result gives the continuity of the summation functional, being the extension of Lemma 2.9 of \cite{roueff-soulier15} to our setting. Recall that $\bD([0,T];\bD)$ is the space of right-continuous functions with left limits with respect to $J_1$ (see Section \ref{section-D}).

\begin{theorem}
\label{cont-Phi}
The summation functional $\Phi: M_{p,f}([0,T] \times \bD) \to \bD([0,T];\bD)$ defined by
$$\Phi(\mu)=\Big(\sum_{t_i \leq t}x_i\Big)_{t \in [0,T]} \quad \mbox{if} \quad \mu=\sum_{i=1}^{p}\delta_{(t_i,x_i)},$$
is continuous on the set $M_{p,f}^{*}([0,T] \times \bD)$, where $\bD([0,T];\bD)$ is equipped with the metric $d_{T,\bD}$ given by \eqref{def-d-TD}.
\end{theorem}

\noindent {\bf Proof:} We use a similar argument to page 221 of \cite{resnick07}, combined with the argument of Lemma 2.9 of \cite{roueff-soulier15}.
Let $\mu=\sum_{i=1}^{p}\delta_{(t_i,x_i)} \in M_{p,f}^{*}([0,T] \times \bD)$ and $(\mu_n)_{n \geq 1} \subset M_{p,f}([0,T] \times \bD)$ be such that $\mu_n \stackrel{w}{\to} \mu$. We must prove that:
\begin{equation}
\label{Phi-convergence}
\Phi(\mu_n) \to \Phi(\mu) \quad \mbox{in} \quad \bD([0,T];\bD).
\end{equation}

Note that $\mu_n([0,T] \times \bD) \to \mu([0,T] \times \bD)=p$ implies that
$\mu_n([0,T] \times \bD)=p$ for all $n \geq n_0$ for some $n_0\geq 1$, since $\mu_n([0,T] \times \bD) \in \bZ_{+}$ for all $n$.

Since $\mu$ is simple, the atoms $(t_1,x_1), \ldots,(t_p,x_p)$ are distinct. Hence, there exists $r>0$ such that $\mu(B_r(t_i,x_i))=1$ for all $i=1,\ldots,p$, where $B_r(t_i,x_i)$ is the ball of radius $r$ and center $(t_i,x_i)$. Fix $i=1,\ldots,p$. For any $r' \in (0,r)$, $\mu(\partial B_{r'}(t_i,x_i))=0$ and hence, $\mu_n(B_{r'}(t_i,x_i)) \to \mu(B_{r'}(t_i,x_i))=1$. Therefore, for any $r' \in (0,r)$, there exists $N_i(r') \geq n_0$ such that $\mu_n(B_{r'}(t_i,x_i))=1$ for all $n \geq N_i(r')$. In particular, for $r'=r/2$ there exists $N_i:=N_{i}(r/2)$ such that $\mu_n(B_{r/2}(t_i,x_i))=1$ for all $n \geq N_i$. We infer that for any $n \geq N_i$, $\mu_n$ has exactly one atom in $B_{r/2}(t_i,x_i)$, which we denote by $(t_i^{n},x_i^{n})$. We claim that:
$$(t_i^n,x_i^n) \to (t_i,x_i) \quad \mbox{in} \quad [0,T] \times \bD, \quad
\mbox{i.e.} \quad t_i^n \to t_i \ \mbox{and} \ x_i^n \stackrel{J_1}{\to} x_i.$$
To see this, let $r' \in (0,r/2$ be arbitrary. We known that for any $n \geq N_i(r')$, $\mu_n$ has exactly one atom in $B_{r'}(t_i,x_i)$, and since $B_{r'}(t_i,x_i) \subset B_{r/2}(t_i,x_i)$, this atom must be $(t_i^n,x_i^n)$. Hence, $(t_i^n,x_i^n) \in B_{r'}(t_i,x_i)$ for any $n \geq N_i(r')$.

Let $N_0=\max_{i\leq p}N_i$.
For any $n\geq N_0$, $\mu_n=\sum_{i=1}^{p}\delta_{(t_i^{n},x_i^{n})}$ and
$\Phi(\mu_n)=(\sum_{t_i^n \leq t}x_i^n)_{t \leq T}$.

The points $t_1, \ldots,t_p$ are distinct, since $\mu$ cannot have two atoms with the same time coordinate, by property {\em (iii)} in the definition of $M_{p,f}^{*}([0,T] \times \bD)$. Pick $\delta_0>0$ such that $t_{i+1}-t_i >2 \delta_0$ for all $i=1,\ldots,p-1$. Let $\delta \in (0,\delta_0)$ be arbitrary. By the choice of $\delta_0$, the intervals $(t_{i}-\delta,t_i+\delta)$, $i=1,\ldots,p$ are non-overlapping.

By property {\em (ii)} in the definition of $M_{p,f}^{*}([0,T] \times \bD)$, ${\rm Disc}(x_i) \cap {\rm Disc}(x_j)\not=\emptyset$ for all $i \not=j$.
By Theorem 4.1 of \cite{whitt80}, it follows that $\sum_{i=1}^{k}x_i^n \stackrel{J_1}{\to}\sum_{i=1}^{k}x_i$ for all $i\leq p$. Hence, there exists $n_1(\delta) \geq N_0$ such that for all $n \geq n_1(\delta)$, $|t_k^n-t_k| \leq \delta$
and $d_{J_1}^{0}(\sum_{i=1}^{k}x_i^n, \sum_{i=1}^{k}x_i )\leq \delta$ for all $k\leq p$.
Let $\lambda_n\in\Lambda_T$ be such that $\lambda_n(t_i^n)=t_i$ for all $i=1,\ldots p$ and $\lambda_n$ is a linear function between $t_i^n$ and $t_{i+1}^n$. By relation (7.20) of \cite{resnick07}, $\|\lambda_n-e\|_{T} \leq 3 \delta$ for all $n \geq n_1(\delta)$.

Recalling definitions \eqref{def-d-TD} and \eqref{def-rho-TD} of distances $d_{T,\bD}$ and $\rho_{T,\bD}$, for any $n \geq n_1(\delta)$, we have:
\begin{align*}
& \rho_{T,\bD}(\Phi(\mu),\Phi(\mu)\circ \lambda_n^{-1})=\sup_{t \in [0,T]}d_{J_1}^{0}\big(\Phi(\mu)(t),\Phi(\mu_n)(\lambda_n^{-1}(t))\big)\\
& \quad =\sup_{t \in [0,T]}d_{J_1}^0\Big(\sum_{t_i \leq t}x_i,\sum_{\lambda_n(t_i^n) \leq t}x_i^n \Big)=\max_{k \leq p} d_{J_1}^{0}\Big(\sum_{i=1}^{k}x_i,\sum_{i=1}^{k}x_i^n \Big)<\delta,
\end{align*}
and hence $d_{T,\bD}(\Phi(\mu),\Phi(\mu_n)) \leq 3 \delta$. This concludes the proof of \eqref{Phi-convergence}. $\Box$

\vspace{3mm}

The following corollary is an immediate consequence of the previous two results.

\begin{corollary}
\label{corol-cont}
The function $Q:\widehat{M}_p([0,\infty) \times \barD) \to \bD([0,T];\bD)$ given by
$$Q(m)=\Big(\sum_{t_i \leq t}r_iz_i1_{\{r_i \in (\e,\infty)\}}\Big)_{t \in [0,T]} \quad \mbox{if} \quad m=\sum_{i\geq 1}\delta_{(t_i,r_i,z_i)},$$
is continuous on the set $\cU=\cA \cap \Psi^{-1}(M_{p,f}^*([0,T] \times \bD))$, where
$\bD([0,T];\bD)$ is equipped with the distance $d_{T,\bD}$ given by \eqref{def-d-TD}.
(The function $Q=Q_{\e,T}$ and the set $\cU_{\e,T}$ depend on $\e$ and $T$. To simplify the writing, we omit the indices $\e,T$.)
\end{corollary}

\noindent {\bf Proof:} The conclusion follows by Lemma \ref{cont-Psi} and Theorem \ref{cont-Phi}, since $Q=\Phi \circ \Psi$. $\Box$

\subsection{Convergence of truncated sums}

In this subsection, we consider a sequence $(X_i)_{i \geq 1}$ of i.i.d. regularly varying random elements in $\bD$, and we prove that the sequence $(S_n^{(\e)})_{n \geq 1}$ of truncated sums defined by:
\begin{equation}
\label{truncated-sum}
S_n^{(\e)}(t)=\frac{1}{a_n}\sum_{i=1}^{[nt]}X_i 1_{\{\|X_i\|>a_n \e \}}, \quad \mbox{for any} \ t\geq 0
\end{equation}
converges in distribution in the space $\bD([0,\infty);\bD)$ to the process $Z^{(\e)}$ given by \eqref{def-Z-e}.

The following result together with Corollary \ref{corol-cont} will allow us to apply the continuous mapping theorem. For this result, we need Assumption B.

\begin{theorem}
\label{N-in-U}
Let $N$ be a Poisson random measure on $[0,\infty) \times \barD$ of intensity ${\rm Leb} \times \overline{\nu}$, where $\overline{\nu}$ is given by \eqref{def-bar-nu}. If $\Gamma_1$ satisfies Assumption B, then  $$N \in \cU_{\e,T} \quad a.s.$$ for any $\e>0$ and $T>0$, where $\cU_{\e,T}$ is the set given in Corollary \ref{corol-cont}.
\end{theorem}

\noindent {\bf Proof:} We have to show that with probability $1$, $N$ satisfies the two conditions listed in Lemma \ref{cont-Psi}, and $\xi=\Psi_{\e,T}(N) \in M_{p,f}^{*}([0,T] \times \bD)$.

We begin with the conditions of Lemma \ref{cont-Psi}. For any $n \geq 1$, $E[N([n-1,n) \times \{\e,\infty\} \times \SD)]=c\nu_{\alpha}(\{\e,\infty\})=0$ and hence
$N([n-1,n) \times \{\e,\infty\} \times \SD)=0$ a.s. By additivity,
$N([0,\infty) \times \{\e,\infty\} \times \SD)=0$ a.s. Similarly, $N(\{0,T\} \times (\e,\infty) \times \SD)=0$ a.s. 

Next, we show that with probability $1$, $\xi$ satisfies conditions {\em (i)}-{\em (iii)} given in Definition \ref{def-M-pf-star}. First, we show that $\xi$ is a Poisson random measure on $[0,T] \times \bD$ of intensity ${\rm Leb} \times \nu^{(\e)}$ where $\nu^{(\e)}=\overline{\nu}|_{(\e,\infty) \times \SD} \circ U^{-1}$ and $U:(\e,\infty) \times \SD \to \bD$ is given by $U(r,z)=rz$. Note that $\xi$ is a point process since $N$ is a point process and $\Psi_{\e,T}$ is measurable. So, it suffices to show that the Laplace functional of $\xi$  is given by \eqref{Laplace-PRM}. Let $g:[0,T] \times \bD \to [0,\infty)$ be a bounded measurable function with bounded support. By \eqref{integral-Psi-m},
\begin{align*}
L_{\xi}(g)&=E\left[\exp\left(-\int_{[0,T] \times \bD} g d\xi\right)\right]=E\left[\exp\left(-\int_{[0,T] \times (\e,\infty) \times \SD} g(t,rz)N(dt,dr,dz) \right) \right]\\
&=\exp\left\{-\int_{[0,T] \times (\e,\infty) \times \SD}(1-e^{-g(t,rz)})dt \overline{\nu}(dr,dz) \right\}\\
&=\exp\left\{-\int_{[0,T] \times \bD}(1-e^{-g(t,x)})\, dt\, \overline{\nu}^{(\e)}(dx) \right\}.
\end{align*}
Since ${\rm Leb} \times \nu^{(\e)}$ is diffuse, $\xi$ is simple a.s. So, $\xi$ satisfies condition {\em (i)} with probability $1$.

To show that $\xi$ satisfies condition {\em (ii)} with probability $1$, we represent its points as follows. Let
$P_i=c^{1/\alpha}\Gamma_i^{-1/\alpha}$ where $\Gamma_i=\sum_{j=1}^{i}E_j$ and $(E_i)_{i \geq 1}$ are i.i.d. exponential random variables of mean $1$. Let $(W_i)_{i\geq 1}$ be an independent sequence of i.i.d. random elements in $\SD$ of law $\Gamma_1$. By the extension of Proposition 5.3 of \cite{resnick07} to Polish spaces, $\sum_{i\geq 1}\delta_{(P_i,W_i)}$ is a Poisson random measure on $(0,\infty) \times \SD$ of intensity $\overline{\nu}$, and so, $\sum_{i\geq 1}\delta_{(P_i,W_i)}1_{\{P_i>\e\}}$ is a Poisson random measure on $(\e,\infty) \times \SD$ of intensity $\overline{\nu}|_{(\e,\infty) \times \SD}$. By the extension of Proposition 5.2 of \cite{resnick07} to Polish spaces,
$\sum_{i\geq 1}\delta_{P_iW_i}1_{\{P_i>\e\}}$ is a Poisson random measure on $\bD$ of intensity $\nu^{(\e)}$. Finally, by the extension of Proposition 5.3 of \cite{resnick07}, $\xi'=\sum_{i\geq 1}\delta_{(\tau_i,P_iW_i)}1_{\{P_i>\e\}}$ is a Poisson random measure on $[0,T] \times \bD$ of intensity ${\rm Leb} \times \nu^{(\e)}$, where $(\tau_i)_{i \geq 1}$ are i.i.d. uniformly distributed on $[0,T]$, independent of $(E_i)_{i \geq 1}$ and $(W_i)_{i \geq 1}$. Hence $\xi \stackrel{d}{=}\xi'$.

Consider the event $A= \cap_{i \not=j}A_{i,j}$, where $A_{i,j}=\{{\rm Disc}(W_i) \cap {\rm Disc}(W_j) =\emptyset\}$. Let $F=\{(x,y) \in \SD \times \SD; {\rm Disc}(x) \cap {\rm Disc}(y)\not=\emptyset\}$. By Fubini's theorem and Assumption B,
$$P(A_{i,j}^c)=P((W_i,W_j) \in F)=(\Gamma_1 \times \Gamma_1)(F)=\int_{\SD}\Gamma_1(F_x)\Gamma_1(dx)=0,$$
where $F_x=\{y \in \SD; (x,y) \in F\}=\cup_{s \in {\rm Disc}(x)}\{y\in \SD; s \in {\rm Disc}(y)\}$. Hence, $P(A)=1$.

Let $B$ be the event on which $\xi(\{(t,x),(t',x')\})\leq 1$ for all $(t,x),(t',x') \in [0,T] \times \bD$ with $x\not=x'$ and ${\rm Disc}(x)\cap {\rm Disc}(x')\not=\emptyset$, and $B'$ the similar event with $\xi$ replaced by $\xi'$. Since $\xi \stackrel{d}{=}\xi'$, $P(B)=P(B')$. We claim that $A \subset B'$. \big(To see this, let $\omega \in (B')^{c}$. Then, there exist $(t,x),(t',x') \in [0,T] \times \bD$ with $x\not=x'$ and ${\rm Disc}(x)\cap {\rm Disc}(x')\not=\emptyset$ such that  $\xi'(\omega;\{(t,x),(t',x')\})\geq 2$. This means that both $(t,x)$ and $(t',x')$ are atoms of $\xi'(\omega)$. But the atoms of $\xi'(\omega)$ are of the form $(\tau_i(\omega),P_i(\omega)W_i(\omega))$ with $P_i(\omega)>\e$. Hence, there exist $i\not=j$ with $P_i(\omega)>\e$ and $P_j(\omega)>\e$ such that $(t,x)=(\tau_i(\omega),P_i(\omega)W_i(\omega))$ and $(t',x')=(\tau_j(\omega),P_j(\omega)W_j(\omega))$. This proves that $\omega \in A_{ij}^c \subset A^c$.\big) Hence, $P(B)=P(B')=P(A)=1$. This proves that $\xi$ satisfies condition {\em (ii)} with probability $1$.

Finally, to show that $\xi$ satisfies condition {\em (iii)} with probability $1$, we let $C=\cap_{i\not=j}C_{i,j}$, where $C_{i,j}=\{\tau_i \not=\tau_j\}$. Note that $P(C)=1$ since for all $i\not=j$
$$P(C_{i,j}^c)=P(\tau_i=\tau_j)=\frac{1}{T^2}\int_0^T \int_{0}^{T}1_{\{x=y\}}dxdy=0.$$
Let $D$ be the event on which $\xi(\{t_0\} \times \bD) \leq 1$ for all $t_0 \in [0,T]$, and $D'$ the similar event with $\xi$ replaced by $\xi'$.  Since $\xi \stackrel{d}{=}\xi'$, $P(D)=P(D')$. We claim that $C \subset D'$.
\big(To see this, let $\omega \in (D')^c$. Then there exists $t_0 \in [0,T]$ such that $\xi'(\omega; \{t_0\} \times \bD)\geq 2$. This means that $\xi'(\omega)$ has at least two atoms with time coordinate $t_0$. Using the form of the atoms of $\xi'(\omega)$, we infer that there exist $i\not=j$ such that $\tau_i(\omega)=\tau_j(\omega)=t_0$. This proves that $\omega \in C_{i,j}^c \subset C^c$.\big) Hence, $P(D)=P(D')=P(C)=1$. This proves that $\xi$ satisfies condition {\em (iii)} with probability $1$. $\Box$

\vspace{3mm}

The next result gives the convergence of the truncated sums of i.i.d. regularly varying elements in $\bD$.

\begin{theorem}
\label{tr-sum-conv}
Let $(X_i)_{i \geq 1}$ be i.i.d. random elements in $\bD$ such that $X_1 \in RV(\{a_n\},\overline{\nu},\barD)$. Let $\alpha$ be the index of $X$ and $\Gamma_1$ be the spectral measure of $X$. Suppose that $\alpha<2,\alpha\not=1$ and $\Gamma_1$ satisfies Assumption B. If $\{S_n^{(\e)},n\geq 1\}$ and $Z^{(\e)}$ are given by \eqref{truncated-sum}, respectively \eqref{def-Z-e}, then for any $\e>0$ and $T>0$
$$S_n^{(\e)}(\cdot) \stackrel{d}{\to} Z^{(\e)}(\cdot) \quad in \quad \bD([0,T];\bD) \quad
as \ n\to \infty,$$ where $\bD([0,T];\bD)$ is equipped with distance $d_{T,\bD}$ given by \eqref{def-d-TD}. Moreover, $P(s \in {\rm Disc}(Z^{(\e)}(t)) \linebreak for some \ t >0)=0$ for all $s \in [0,1]$ and $\e>0$.
\end{theorem}

\noindent {\bf Proof:} By Proposition \ref{prop3-21} with $E=\barD$ and $X_{i,n}=(\|X_i\|/a_n,X_i/\|X_i\|)$,
$$N_n=\sum_{i \geq 1}\delta_{\big(\frac{i}{n},\,\frac{\|X_i\|}{a_n},\,\frac{X_i}{\|X_i\|}\big)} \stackrel{d}{\to}N,$$
where $N$ is a Poisson random measure on $[0,\infty) \times \barD$ of intensity ${\rm Leb} \times \overline{\nu}$.

Note that $S_n^{(\e)}=Q(N_n)$ and $Z^{(\e)}=Q(N)$, where $Q$ is the map given in Corollary \ref{corol-cont}. By the continuous mapping theorem and Theorem \ref{N-in-U}, $S_n^{(\e)} \stackrel{d}{\to} Z^{(\e)}$ in $\bD([0,T];\bD)$.

To prove the last statement, we fix $s \in [0,1]$ and we let $\Omega_T=\cup_{t \in [0,T]}\{s \in {\rm Disc}(Z^{(\e)}(t))\}$. It is enough to prove that $P(\Omega_T)=0$ for all $T>0$. From \eqref{def-Z-e}, we see that if $W_i$ is continuous at $s$ for all $i \geq 1$, then $Z^{(\e)}(t)$ is continuous at $s$ for all $t \in [0,T]$. Hence, $\Omega_T \subset \cup_{i\geq 1} \{s \in {\rm Disc}(W_i)\}$. The fact that $P(\Omega_T)=0$ follows by Assumption B, since
$P(s \in {\rm Disc}(W_i))=\Gamma_1(\{z \in \SD; s \in {\rm Disc}(z)\})=0$. $\Box$

\subsection{Approximation in the case $\alpha<1$}

In this subsection, we prove the approximation result (Theorem \ref{main2}) in the case $\alpha<1$.

The first result shows that a certain asymptotic negligibility condition holds automatically in the case $\alpha<1$.

\begin{lemma}
\label{AN-alpha-less-1}
Let $(X_i)_{i \geq 1}$ be i.i.d. random elements in $\bD$ such that $X_1 \in RV(\{a_n\},\overline{\nu},\barD)$. Suppose that $\alpha \in (0,1)$, where $\alpha$ is the index of $X$. Let $\{S_n^{(\e)},n\geq 1\}$ be given by \eqref{truncated-sum} and $S_n(t)=a_n^{-1}\sum_{i=1}^{[nt]}X_i$ for all $t \geq 0,n\geq 1$. Then for any $\delta>0$ and $T>0$
$$\lim_{\e \downarrow 0}\limsup_{n \to \infty}P(\|S_n-S_n^{(\e)}\|_{T,\bD}>\delta)=0,$$
and in particular, $\lim_{\e \downarrow 0}\limsup_{n \to \infty}P(d_{T,\bD}(S_n,S_n^{(\e)})>\delta)=0$.
\end{lemma}

\noindent {\bf Proof:} Let $\delta>0$ and $T>0$ be arbitrary. Since $S_n(t)-S_n^{(\e)}(t)=a_n^{-1}\sum_{i=1}^{[nt]}X_i 1_{\{|X_i\| \leq a_n \e\}}$,
$$\|S_n-S_n^{(\e)}\|_{T,\bD}=\frac{1}{a_n} \max_{k \leq [nT]}\left\| \sum_{i=1}^{k}X_i 1_{\{\|X_i\|\leq a_n \e\}}\right\|\leq \frac{1}{a_n} \sum_{i=1}^{[nT]}\|X_i\|1_{\{\|X_i\| \leq a_n \e\}}.$$
By Markov's inequality,
$$P(\|S_n-S_n^{(\e)}\|_{T,\bD} >\delta)\leq  
\frac{1}{\delta a_n} [nT] \, E(\|X_1\| 1_{\{\|X_1\|\leq a_n \e\}}).$$
Since $\|X_1\|$ is regularly varying of index $\alpha<1$,
$E(\|X_1\| 1_{\{\|X_1\|\leq x \}}) \sim \frac{\alpha}{1-\alpha} x P(\|X_1\|>x)$ as $x \to \infty$, by Karamata's theorem (e.g. Theorem 2.1 of \cite{resnick07}), and hence, by \eqref{RV-norm-X},
$$\frac{n}{a_n}E(\|X_1\| 1_{\{\|X_i\|\leq a_n \e\}}) \sim \frac{\alpha}{1-\alpha}  \e nP(\|X_1\|>a_n \e) \sim \frac{\alpha}{1-\alpha} \, c\, \e^{1-\alpha} \quad \mbox{as} \quad n \to \infty.$$
Here $f(x) \sim g(x)$ as $x \to \infty$ means that $f(x)/g(x) \to 1$ as $x \to \infty$. Therefore,
$$\limsup_{n \to \infty}P(\|S_n-S_n^{(\e)}\|_{T,\bD} >\delta) \leq \frac{T}{\delta} \cdot \frac{\alpha}{1-\alpha} \, c \, \e^{1-\alpha}.$$
The conclusion follows letting $\e \downarrow 0$, and using the fact that $\alpha<1$.  $\Box$

\vspace{3mm}

{\bf Proof of Theorem \ref{main2}.a)} By Theorem 2.8 of \cite{whitt80}, it is enough to prove that $$S_n(\cdot)\stackrel{d}{\to} \widetilde{Z}(\cdot) \quad \mbox{in} \quad \bD([0,T];\bD),$$ for any $T>0$, where $\bD([0,T];\bD)$ is equipped with distance $d_{T,\bD}$.
This follows by Theorem 4.2 of \cite{billingsley68}, whose hypotheses are verified due to Theorem \ref{main-conv-th1}, Theorem \ref{tr-sum-conv} and Lemma \ref{AN-alpha-less-1}. $\Box$


\subsection{Approximation in the case $\alpha>1$}
\label{subsection-approx2}

In this subsection, we prove the approximation result (Theorem \ref{main2}) in the case $\alpha>1$.

The following result is the counterpart of Lemma \ref{AN-alpha-less-1} for the case $\alpha>1$.

\begin{lemma}
\label{AN-alpha-greater-1}
Let $(X_i)_{i \geq 1}$ be i.i.d.random elements in $\bD$ such that $X_1 \in RV(\{a_n\},\overline{\nu},\barD)$. Suppose that $\alpha \in (1,2)$, where $\alpha$ is the index of $X_1$.
Let $\{S_n^{(\e)},n\geq 1\}$  be given by \eqref{truncated-sum}. \\
For any $t \geq 0$ and $n\geq 1$, let $S_n(t)=\sum_{i=1}^{[nt]}X_i/a_n$, $$\overline{S}^{(\e)}_n(t)=S_n^{(\e)}(t)-E[S_n^{(\e)}(t)] \quad \mbox{and} \quad \overline{S}_n(t)=S_n(t)-E[S_n(t)].$$
If \eqref{AN} holds for any $\delta>0$ and $T>0$, then for any $\delta>0$ and $T>0$,
\begin{equation}
\label{AN2}
\lim_{\e \downarrow 0}\limsup_{n \to \infty}P(\|\overline{S}_n-\overline{S}_n^{(\e)}\|_{T,\bD}>\delta)=0,
\end{equation}
and in particular, $\lim_{\e \downarrow 0}\limsup_{n \to \infty}P(d_{T,\bD}(\overline{S}_n,\overline{S}_n^{(\e)})>\delta)=0$.
\end{lemma}

\noindent {\bf Proof:} Since $\overline{S}_n(t)-\overline{S}_n^{(\e)}(t)=\sum_{i=1}^{[nt]}Y_{i,n}$ with $Y_{i,n}=a_n^{-1}\big(X_i 1_{\{|X_i\| \leq a_n \e\}}-E(X_i 1_{\{|X_i\| \leq a_n \e\}})\big)$,
$$\|\overline{S}_n-\overline{S}_n^{(\e)}\|_{T,\bD}=\sup_{t \in [0,T]} \|\overline{S}_n(t)-\overline{S}_n^{(\e)}(t)\| =\max_{k \leq [nT]} \left\|\sum_{i=1}^{k} Y_{i,n}\right\|.$$
By L\'evy-Octaviani inequality, which is valid for independent random elements in a normed space (see Proposition 1.1.1 of \cite{kwapien-woyczynski92}), for any $\delta>0$,
$$P(\|\overline{S}_n-\overline{S}_n^{(\e)}\|_{T,\bD}>\delta) \leq 3 \max_{k \leq [nT]}P\left(\left\| \sum_{i=1}^{k}Y_{i,n}\right\|>\delta /3\right).$$
The conclusion follows by \eqref{AN}. $\Box$

\vspace{3mm}

To deal with the centering constants, we need to use the fact that addition is continuous in the space $\bD([0,T];\bD)$ equipped with the distance $d_{T,\bD}$. To deduce this, we cannot simply apply Theorem 4.1 of \cite{whitt80} with $(S,m)=(\bD,d_{J_1}^0)$, since we do not know if the relation $d_{J_1}^0(x+y,x'+y') \leq
d_{J_1}^0(x,x')+d_{J_1}^0(y,y')$ holds for any $x,x',y,y' \in \bD$, as required on p.78 of \cite{whitt80}. Although the general question of continuity of the addition on $\bD([0,T];\bD)$ remains open, we were able to find a weaker version of this result which is sufficient for our purposes. This is contained in the lemma below.

\begin{lemma}
\label{sum-conv-d}
Let $(f_n)_{n\geq 1} \subset \bD$ and $f \in \bD$ be such that $f_n \stackrel{J_1}{\to}f$. Consider $(y_n)_{n \geq 1} \subset \bD([0,T];\bD)$ and $y \in \bD([0,T];\bD)$ defined as follows: for any $t \in [0,T]$,
\begin{equation}
\label{def-yn-y}
y_n(t)=\frac{[nt]}{n}f_n \quad \mbox{and} \quad y(t)=tf.
\end{equation}
Then $\rho_{T,\bD}(y_n,y) \to 0$. Moreover, if $f$ is continuous, then for any sequence $(x_n)_{n\geq 1}\subset \bD([0,T];\bD)$ and $x \in \bD([0,T];\bD)$ such that $d_{T,\bD}(x_n,x)\to 0$, we have:
\begin{equation}
\label{sum-x-y}
d_{T,\bD}(x_n+y_n,x+y)\to 0.
\end{equation}
\end{lemma}

\noindent {\bf Proof:} We first prove that $\rho_{T,\bD}(y_n,y) \to 0$. Since $f_n \stackrel{J_1}{\to}f$, there exists a sequence $(\rho_n)_{n\geq 1} \subset \Lambda$ such that $\|\rho_n\|^0 \to 0$ and $\|f_n -f \circ \rho_n\|\to 0$. Let $z_n(t)=\frac{[nt]}{n}f$. Let $\e>0$ be arbitrary. Then, there exists $N_{\e}$ such that for all $n\geq N_{e}$, $\|\rho_n\|^0 <\e$ and $\|f_n-f\circ \rho_n\|<\e/T$. Hence,  for any $t \in [0,T]$ and $n\geq N_{\e}$, $\|y_n(t)-z_n(t)\circ \rho_n\|\leq t \|f_n-f\circ \rho_n\|<\e$ and
$$d_{J_1}^0(y_n(t),z_n(t))\leq \|\rho_n\|^0 \vee \|y_n(t)-z_n(t)\circ \rho_n\|<\e.$$
On the other hand, there exists $N_{\e}'$ such that, for any $t \in [0,T]$ and $n \geq N_{\e}'$,
$$d_{J_1}^0(z_n(t),y(t))\leq \|z_n(t)-y(t)\|=\left|\frac{[nt]}{n}-t\right|\cdot \|f\|\leq \frac{1}{n}\|f\|<\e.$$
This shows that
$\rho_{T,\bD}(y_n,y)=\sup_{t\in [0,T]}d_{J_1}^0(y_n(t),y(t))<2\e$ for any $n \geq N_{\e}\vee N_{\e}'$.

We now prove \eqref{sum-x-y}. For any $t \in [0,T]$, we denote $x(t)=\{x(t,s)\}_{s \in [0,1]}$, and we use a similar notation for $y(t),x_n(t)$ and $y_n(t)$.
Let $\e>0$ be arbitrary. Since $f$ is uniformly continuous, there exists $\delta_{\e}\in (0,\e)$ such that for any $s,s' \in [0,1]$ with $|s-s'|<\delta_{\e}$,
\begin{equation}
\label{diff-a}
|f(s)-f(s')|<\e.
\end{equation}

Because $d_{T,\bD}(x_n,x)\to 0$, there exists a sequence $(\lambda_n)_{n\geq1} \subset \Lambda_T$ such that $\|\lambda_n-e\|_{T} \to 0$ and $\rho_{T,\bD}(x_n \circ \lambda_n,x) \to 0$. Pick $0<\eta_{\e}<\e \wedge \ln(\delta_{\e}+1)$ arbitrary.
Then, there exists $N_{\e}^{(1)}$ such that for any $n \geq N_{\e}^{(1)}$,
$\sup_{t \in [0,T]}|\lambda_n(t)-t| <\e$ and $\sup_{t \in [0,T]}
d_{J_1}^0(x_n(\lambda_n(t)),x(t))<\eta_{\e}$.
Using definition \eqref{def-dJ1} of $d_{J_1}^0$, it follows that for any $n\geq N_{\e}^{(1)}$ and for any $t \in [0,T]$, there exists $\mu_{t}^{(n)}\in \Lambda$ such that $\|\mu_t^{(n)}\|^0<\eta_{\e}$ and
\begin{equation}
\label{bound-xn-s}
\sup_{s \in [0,1]}
|x_n\big(\lambda_n(t),\mu_t^{(n)}(s)\big)-x(t,s)|<\eta_{\e}.
\end{equation}
By inequality \eqref{lambda0-equiv} and the choice of $\eta_{\e}$,
$\sup_{s \in [0,1]}|\mu_t^{(n)}(s)-s|<e^{\eta_{\e}}-1<\delta_{\e}$.

Note that $\|f_n-f\| \to 0$, since $f_n \stackrel{J_1}{\to}f$ and $f$ is continuous. Hence, there exists $N_{\e}^{(2)}$ such that
$\sup_{s \in [0,1]}|f_n(s)-f(s)|<\e$ for any $n \geq N_{\e}^{(2)}$.
By \eqref{diff-a}, for any $n\geq N_{\e}^{(1)}\vee N_{\e}^{(2)}$,
$$|f_n(\mu_t^{(n)}(s))-f(s)|\leq |f_n(\mu_t^{(n)}(s))-f(\mu_t^{(n)}(s))|+|f(\mu_t^{(n)}(s))-f(s)|<2\e.$$

Choose $N_{\e}^{(0)}$ such that $1/n<\e$ for any $n \geq N_{\e}^{(0)}$. Then, for any $n \geq N_{\e}^{(0)}$ and $t \in [0,T]$,
$$\left|\frac{[n\lambda_n(t)]}{n}-t \right|\leq \left|\frac{[n\lambda_n(t)]}{n}-\lambda_n(t) \right|+|\lambda_n(t)-t|\leq \frac{1}{n}+\e<2\e.$$


Since $\|f_n -f\| \to 0$, it follows that $C:=\sup_{n\geq 1}\|f_n\|<\infty$.
Let $N_{\e}=N_{e}^{(0)} \vee N_{e}^{(1)} \vee N_{e}^{(2)}$.
Using the definitions of $y_n$ and $y$, it follows that for any $n \geq N_{\e}$, $t \in [0,T]$ and $s \in [0,1]$,
\begin{align*}
|y_n\big(\lambda_n(t),\mu_t^{(n)}(s) \big)-y(t,s)| \leq \left|\frac{[n\lambda_n(t)]}{n}-t \right|
|f_n(\mu_t^{(n)}(s))|+t|f_n(\mu_t^{(n)}(s))-f(s)|<2\e(C+T).
\end{align*}
and hence, by \eqref{bound-xn-s},
\begin{align*}
|(x_n+y_n)\big(\lambda_n(t),\mu_t^{(n)}(s) \big)-(x+y)(t,s)|<\eta_{\e}+2\e(C+T)<\e[1+2(C+T)].
\end{align*}
To summarize, we have proved that for any $n\geq N_{\e}$, and $t \in [0,T]$, there exists $\mu_{t}^{(n)}\in \Lambda$ such that
$\|\mu_t^{(n)}\|^0<\eta_{\e}<\e$ and
$\|(x_n+y_n)(\lambda_n(t))\circ \mu_t^{(n)}-(x+y)(t)\|<\e[1+2(C+T)]$. By definition \eqref{def-dJ1} of $d_{J_1}^0$, this implies that for any $n\geq N_{\e}$ and $t \in [0,T]$,
$$d_{J_1}^0 \big((x_n+y_n)(\lambda_n(t)), (x+y)(t)\big)<\e[1+2(C+T)].$$
Therefore, for any $n \geq N_{\e}$
$$\rho_{T,\bD}((x_n+y_n)\circ \lambda_n,x+y)=\sup_{t \in [0,T]}d_{J_1}^0\big((x_n+y_n)(\lambda_n(t)), (x+y)(t)\big)<\e[1+2(C+T)].$$
Since $\|\lambda_n-\e\|_{T}<\e$, using definition \eqref{def-d-TD} of $d_{T,\bD}$, we conclude that $d_{T,\bD}(x_n+y_n,x+y)<\e[1+2(C+T)]$ for any $n \geq N_{\e}$. $\Box$

\vspace{3mm}

\begin{remark}
\label{remark-phi-cont}
{\rm In the proof of Theorem 2.12 of \cite{roueff-soulier15}, it was shown that, in a more general context than here, the function $s \mapsto E[Z^{(\e)}(1,s)]$  is continuous on $[0,1]$. In our case, $E[Z^{(\e)}(1,s)]=c\varphi(s)\int_{\e}^{\infty}r \nu_{\alpha}(dr)$, where $\varphi(s)=\int_{\bS_D}z(s) \Gamma_1(dz)$ for all $s \in [0,1]$. The continuity of $\varphi$ can be proved directly as follows. By the dominated convergence theorem, $\varphi$ is a c\`adl\`ag function. To show that $\varphi$ is left-continuous, note that for any $s \in [0,1]$,
$$\varphi(s)-\varphi(s-)=\int_{\SD} \big(z(s)-z(s-)\big)\Gamma_1(dz)=\int_{\{z \in \SD; z\{s\}>0\}}z\{s\}\Gamma_1(dz),$$
where $z\{s\}=z(s)-z(s-)$ is the jump of $z \in \SD$ at $s$. By Assumption B, the set in the last integral above has $\Gamma_1$-measure $0$, and hence this integral is equal to $0$. }
\end{remark}

The following result gives the convergence of the centered sums.

\begin{theorem}
\label{centered-tr-sum-conv}
Let $(X_i)_{i \geq 1}$ be i.i.d. random elements in $\bD$ such that $X_1 \in RV(\{a_n\}, \overline{\nu},\barD)$.
Let $\alpha$ be the index of $X$ and $\Gamma_1$ be the spectral measure of $X$. Suppose that $\alpha\in (1,2)$ and $\Gamma_1$ satisfies Assumption B. Let $\{S_n^{(\e)},n\geq 1\}$ and $Z^{(\e)}$ be given by \eqref{truncated-sum}, respectively \eqref{def-Z-e}. For any $t\geq 0$, let $\overline{S}_n^{(\e)}(t)=S_n^{(\e)}(t)-E[S_n^{(\e)}(t)]$ and $\overline{Z}^{(\e)}(t)=Z^{(\e)}(t)-E[Z^{(\e)}(t)]$.

Then, for any $\e>0$ and $T>0$
$$\overline{S}_n^{(\e)}(\cdot) \stackrel{d}{\to} \overline{Z}^{(\e)}(\cdot) \quad \mbox{in} \quad \bD([0,T];\bD),$$
where $\bD([0,T];\bD)$ is equipped with distance $d_{T,\bD}$.
\end{theorem}

\noindent {\bf Proof:} Let $X_n=S_n^{(\e)}$ and $X=Z^{(\e)}$. For any $t \geq 0$ and $s \in [0,1]$,
$$y_n(t,s):=-E[S_n^{(\e)}(t,s)]=-\frac{[nt]}{a_n}E[X_1(s) 1_{\{\|X_1\| >a_n \e\}}]=\frac{[nt]}{n}f_n(s),$$
with $f_n(s)=-\frac{n}{a_n}E[X_1(s)1_{\{\|X_1\| >a_n \e\}}]$, and
$$y(t,s):=-E[Z^{(\e)}(t,s)]=-tc\int_{(\e,\infty) \times \SD}  rz(s)\nu_{\alpha}(dr)\Gamma_1(dz)=tf(s),$$
with $f(s)=-c\frac{\alpha}{\alpha-1}\e^{1-\alpha}\varphi(s)$ and $\varphi(s)=\int_{\bS_D}z(s) \Gamma_1(dz)$. This shows that the functions $(y_n)_{n\geq 1}$ and $y$ are of the same form as in \eqref{def-yn-y}.
By Remark \ref{remark-phi-cont}, $\varphi$ is continuous on $[0,1]$.

By Theorem \ref{tr-sum-conv}, $X_n \stackrel{d}{\to} X$ in the space $\bD([0,T];\bD)$ equipped with $d_{T,\bD}$. Since this space is separable (by Theorem \ref{bD-CSMS}), by Skorohod's embedding theorem (Theorem 6.7 of \cite{billingsley99}), there exist random elements $(X_n')_{n\geq 1}$ and $X'$ defined on a probability space $(\Omega',\cF',P')$ such that $X_{n}'\stackrel{d}{=}X_n$ for all $n$, $X' \stackrel{d}{=}X'$ and $d_{T,\bD}(X_n', X') \to 0$ a.s. By Lemma \ref{sum-conv-d}, it follows that
$$d_{T,\bD}(X_n'+y_n,X'+y) \to 0 \quad \mbox{a.s.}$$
This implies that $d_{T,\bD}(X_n'+y_n,X'+y) \to 0$ in probability (and in distribution). By
Corollary to Theorem 3.1 of \cite{billingsley99} (and using again the fact that $\bD([0,T];\bD)$ equipped with $d_{T,\bD}$ is a separable space), we infer that $X_n'+y_n \stackrel{d}{\to}X'+y$ in $\bD([0,T];\bD)$ equipped with $d_{T,\bD}$. Since $(y_n)_{n\geq 1}$ and $y$ are deterministic, $X_n+y_n\stackrel{d}{=}X_n'+y_n$ for any $n$, and $X+y\stackrel{d}{=}X+y$. It follows that $X_n+y_n \stackrel{d}{\to}X+y$ in $\bD([0,T];\bD)$ equipped with $d_{T,\bD}$. $\Box$

\vspace{3mm}

{\bf Proof of Theorem \ref{main2}.b)}
This follows by Theorem 4.2 of \cite{billingsley68} whose hypotheses are verified due to Theorem \ref{main-conv-th2}, Lemma \ref{AN-alpha-greater-1} and Theorem \ref{centered-tr-sum-conv}. $\Box$

\section{Simulations}
\label{section-simulation}

In this section, we simulate the sample paths of a $\bD$-valued $\alpha$-stable L\'evy motion
using Theorem \ref{main2}, by focusing on two examples of a regularly varying process $X$ in $\bD$.

\begin{example}
\label{ex-Levy-sheet}
{\rm The simplest example of a regularly varying process $X=\{X(s)\}_{s \in [0,1]}$ in $\bD$ is the $\alpha$-stable L\'evy motion, which can be simulated using the stable central limit theorem. We recall briefly this result below. Let $\xi,(\xi_j)_{j\geq 1}$ be i.i.d. regularly varying random variables in $\bR$, i.e.
\begin{equation}
\label{RV}
P(|\xi|>x)=x^{-\alpha}L(x) \quad \mbox{and} \quad \lim_{x \to \infty}\frac{P(\xi>x)}{P(|\xi|>x)}=p,
\end{equation}
for some $\alpha \in (0,2)$, $p \in [0,1]$ and a slowly varying function $L$. Let $(a_n)_{n\geq 1}$ be a sequence of real numbers with $a_n \uparrow \infty$ such that $nP(|\xi|>a_n) \to 1$ as $n \to \infty$, i.e. $a_n^{\alpha} \sim nL(a_n)$ as $n \to \infty$.
Condition \eqref{RV} is equivalent to the vague convergence $nP(\xi/a_n \in \cdot) \stackrel{v}{\to} \nu_{\alpha,p}$ in $\overline{\bR}_0$, where
\begin{equation}
\label{def-nu-ap}
\nu_{\alpha,p}(dz)=\big(p\alpha z^{-\alpha-1}1_{(0,\infty)}(z)+q \alpha (-z)^{-\alpha-1}1_{(-\infty,0)}(z) \big)dz
\end{equation}
with $q=1-p$. In other words, for any $x>0$,
$$\lim_{n \to \infty}n P\left(\frac{\xi}{a_n}>x\right)=p x^{-\alpha} \quad \mbox{and} \quad \lim_{n \to \infty}n P\left(\frac{\xi}{a_n}<-x\right)=q x^{-\alpha}.$$
In this case, we write $\xi \in RV(\{a_n\},\nu_{\alpha,p},\overline{\bR}_0)$.
In particular, if
\begin{equation}
\label{cond-L}
\lim_{x \to \infty}L(x)= C>0,
\end{equation}
then $a_n^{\alpha}  \sim Cn$.
We assume that $\alpha\not=1$. Let $\mu=0$ if $\alpha<1$ and $\mu=E(\xi)$ if $\alpha>1$. A classical result, which can be deduced for instance from Theorem 2.7 of \cite{skorokhod57}, states that
\begin{equation}
\label{stable-CLT}
\frac{1}{a_n} \sum_{j=1}^{[n\cdot]} (\xi_j-\mu) \stackrel{d}{\to} X(\cdot) \quad \mbox{in} \quad \bD
\end{equation}
where $X=\{X(s)\}_{s \in [0,1]}$ is an $\alpha$-stable L\'evy motion, with $X(1)$ having a $S_{\alpha}(\sigma_{\alpha},\beta,0)$-distribution. Here $\sigma_{\alpha}^{\alpha}=C_{\alpha}^{-1}$ with $C_{\alpha}$ given by \eqref{def-C}, and $\beta=p-q$. By Property 1.2.15 of \cite{ST94},
$\lim_{x \to \infty}x^{\alpha}P(X(1)>x)=p$ and $\lim_{x \to \infty}x^{\alpha}P(X(1)<-x) = q$. If $L$ satisfies \eqref{cond-L}, this implies that
$X(1) \in RV(\{a_n\},C\nu_{\alpha,p},\overline{\bR}_0)$, since
$$nP\left(\frac{X(1)}{a_n}>x \right)=(n a_n^{-\alpha}) \cdot (a_n x)^{\alpha} P(X(1)>a_n x) \cdot x^{-\alpha} \to C p x^{-\alpha}$$
as $n \to \infty$, and similarly, $nP\left(\frac{X(1)}{a_n}<-x \right) \to C q x^{-\alpha}$. By Lemma 2.1 of \cite{hult-lindskog07}, it follows that $X \in RV(\{a_n\},\overline{\nu},\overline{\bD}_0)$ for a boundedly finite measure $\overline{\nu}$ on $\overline{\bD}_0$. Note that the normalizing sequence $\{a_n\}_n$ for the regular variation of $X$ in $\bD$ is {\em the same} as for $\xi$, if $L$ satisfies \eqref{cond-L}.
In the simulations, we take $a_n=(Cn)^{1/\alpha}$, where $C$ is given by \eqref{cond-L}.

In view of \eqref{stable-CLT}, for any $s \in [0,1]$,
$X(s) \approx \frac{1}{a_m} \sum_{j=1}^{[ms]}(\xi_j-\mu)$, when $m$ is large.

Next, we consider $n$ i.i.d. copies of $X$. For this, let $(\xi_{ij})_{i,j \geq 1}$ be i.i.d. copies of $\xi$. When $m$ is large, we have the following approximations for any $s \in [0,1]$:
$$X_i(s) \approx \frac{1}{a_m}\sum_{j=1}^{[ms]} (\xi_{ij}-\mu), \quad \mbox{for all} \quad i=1,\ldots,n.$$
By Theorem \ref{main2}, the following approximation gives a $\bD$-valued $\alpha$-stable L\'evy motion $Z$:
$$Z(t,s) \approx \frac{1}{a_n} \sum_{i=1}^{[nt]}X_i(s) \approx \frac{1}{a_n a_m} \sum_{i=1}^{[nt]} \sum_{j=1}^{[ms]}(\xi_{ij}-\mu),$$
for any $t,s \in [0,1]$, when $n$ and $m$ are large. (By Theorem \ref{stable-sheet-th2} below, this approximation yields in fact an {\em $\alpha$-stable L\'evy sheet}, which is an example of a $\bD$-valued $\alpha$-stable L\'evy motion, according to Theorem \ref{stable-sheet-th1} below.)

 We consider 5 examples of regularly varying random variables $\xi$ which satisfy \eqref{cond-L}:

{\em (i)} $\xi \sim {\rm Pareto}(\alpha)$, i.e. $\xi$ has density $f(x)=\alpha x^{-\alpha-1}$ if $x>1$; then $L(x)=1$;

{\em (ii)} $\xi$ has a two-sided Pareto distribution, i.e. $\xi$ has density given by $f(x)=p\alpha x^{-\alpha-1}$ if $x>1$ and $f(x)=q \alpha (-x)^{-\alpha-1}$ if $x<-1$, for $p \in (0,1)$ and $q=1-p$; then $L(x)=1$;

{\em (iii)} $\xi \sim \mbox{Fr\'echet}(\alpha)$, i.e. $\xi$ has density $f(x)=\alpha x^{-\alpha-1}e^{-x^{-\alpha}}$ if $x>0$; then $L(x) =x^{\alpha}(1-e^{-x^{-\alpha}}) \to 1$ as $x \to \infty$;

{\em (iv)} $\xi \sim \mbox{Burr}(a,b)$ with $a,b>0$, i.e. $\xi$ has density $f(x)=ab x^{b-1} (1+x^b)^{-a-1}$ for $x>0$; in this case $\alpha=ab$ and $L(x)=(1+x^{-b})^a \to 1$ as $x \to \infty$;


{\em (v)} $\xi \sim S_{\alpha}(\sigma,\beta,\mu)$; in this case $L(x) \to C:=C_{\alpha} \sigma^{\alpha}$ as $x\to \infty$.
}
\end{example}

The following pictures are the 3-dimensional plots of
$(t_k,s_l,Z(t_k,s_l))$ for $k=1,\ldots,n$ and $l=1,\ldots,m$, with $t_k=k/n$ and $s_l=l/m$, when $n=400$ and $m=250$.

\begin{figure}[h!]
    \centering
    \begin{subfigure}[b]{0.48\textwidth}
        \includegraphics[width=\textwidth]{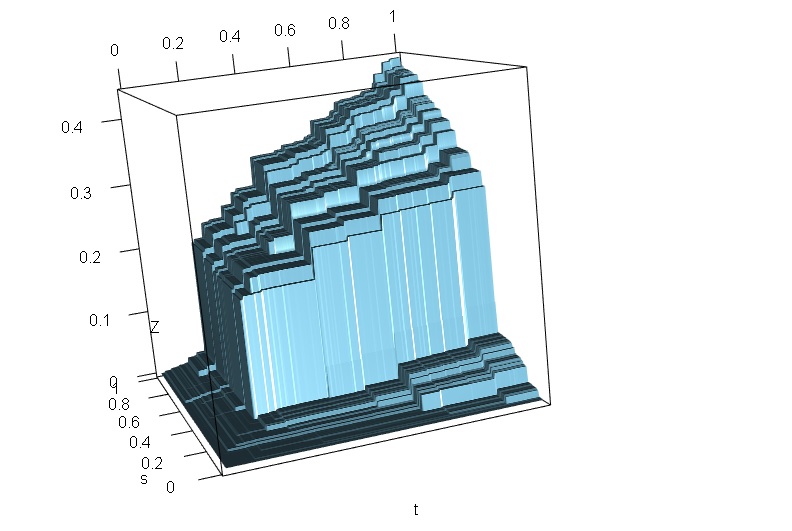}
        \caption{$\alpha=0.5$}
        \label{fig:Pareto1}
    \end{subfigure}
    ~ 
    \begin{subfigure}[b]{0.48\textwidth}
        \includegraphics[width=\textwidth]{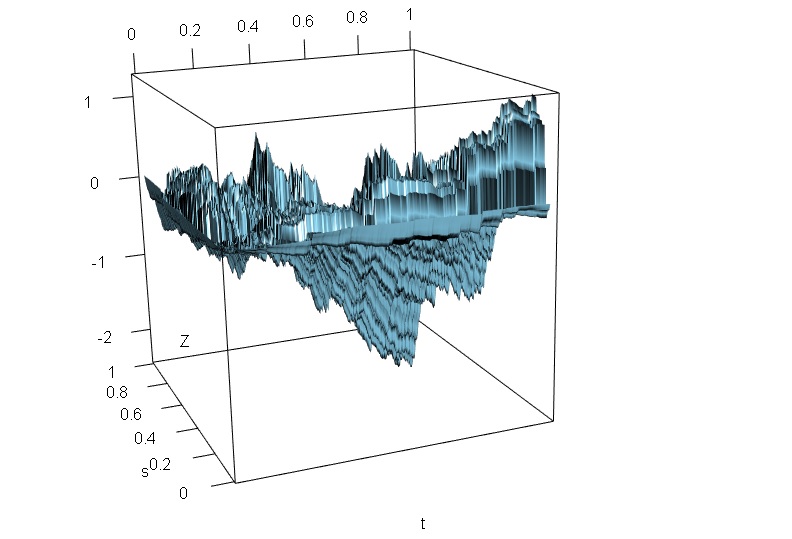}
        \caption{$\alpha=1.5$}
        \label{fig:Pareto2}
    \end{subfigure}
    ~ 
    \caption{$\alpha$-stable L\'evy sheet based on Pareto distribution}
    \label{fig:Pareto}
\end{figure}

\pagebreak

\begin{figure}[h!]
    \centering
    \begin{subfigure}[b]{0.48\textwidth}
        \includegraphics[width=\textwidth]{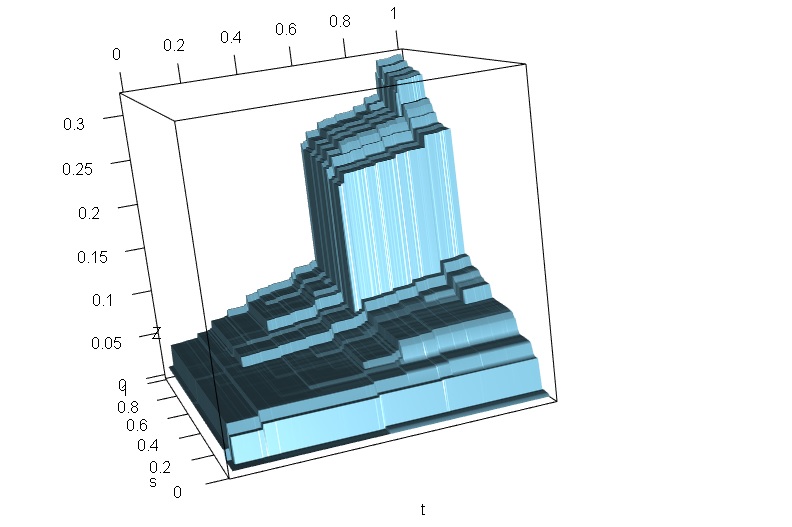}
        \caption{$\alpha=0.5$}
        \label{fig:frechet2}
    \end{subfigure}
    ~ 
    \begin{subfigure}[b]{0.48\textwidth}
        \includegraphics[width=\textwidth]{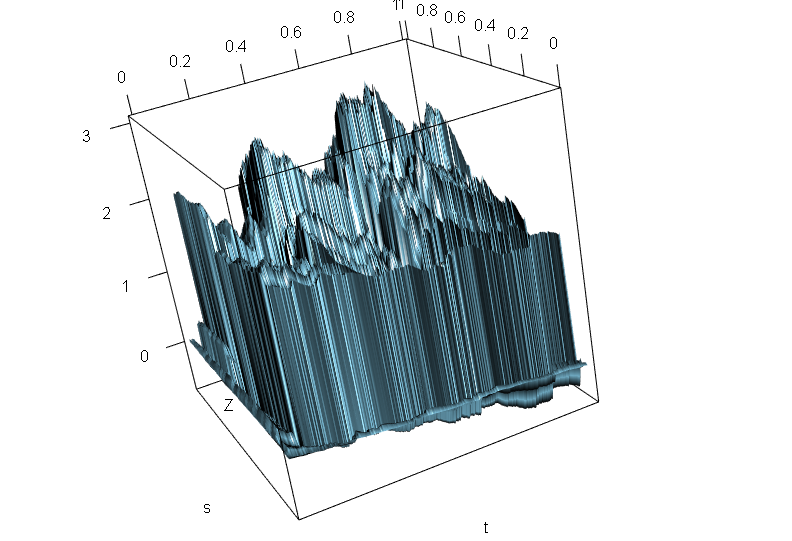}
        \caption{$\alpha=1.5$}
        \label{fig:frechet2}
    \end{subfigure}
    ~ 
    \caption{$\alpha$-stable L\'evy sheet based on Fr\'echet distribution}
    \label{fig:frechet}
\end{figure}

\begin{example}
\label{ex-series}
{\rm In this example, $X=\{X(s)\}_{s \in [0,1]}$ is a regularly varying random element in $\bD$ given by a series, as explained in Example 4.1 of \cite{davis-mikosch08}. Let $Y,(Y_j)_{j \geq 1}$ be i.i.d. random elements in the space $\bC=\bC([0,1])$ of continuous functions on $[0,1]$, such that
\begin{equation}
\label{cond-Y}
0<C_{Y,\alpha}:=E\Big(\sup_{s \in [0,1]}|Y(s)|^{\alpha}\Big)<\infty
\end{equation}
for some $\alpha \in (0,2)$. Let $(\e_j)_{j \geq 1}$ be i.i.d. random variables which take values $1$ and $-1$ with probability $1/2$, and $\Gamma_j=\sum_{i=1}^{j}E_i$ where $(E_i)_{i\geq 1}$ are i.i.d. exponential random variables of mean $1$. Assume that $(Y_j)_{j\geq 1}, (\e_j)_{j\geq 1}$ and $(E_j)_{j\geq 1}$ are independent. By Theorem 1.4.2 of \cite{ST94}, for any $s\in [0,1]$, the series
\begin{equation}
\label{series}
X(s)=\sum_{j\geq 1}\e_j \Gamma_j^{-1/\alpha}Y_j(s) \quad \mbox{converges a.s.}
\end{equation}
and has a $S_{\alpha}(\sigma_s,0,0)$-distribution, with $\sigma_s^{\alpha}=C_{\alpha}^{-1}E|Y(s)|^{\alpha}$ and $C_{\alpha}$ given by \eqref{def-C}. Moreover, the process $X=\{X(s)\}_{s \in [0,1]}$ has sample paths in $\bC$, and is regularly varying in $\bD$. More precisely, $X \in RV(\{a_n\},\overline{\nu},\overline{\bD}_0)$ with sequence $(a_n)_n$ chosen such that $a_n^{\alpha} \sim n C_{Y,\alpha}$, and limiting measure $\overline{\nu}$ specified by (4.3) of \cite{davis-mikosch08}.

In the simulation below, we truncate the series in \eqref{series} by considering only the first $K$ terms (for $K$ large), and we take $Y=W$ where $W=\{W(s)\}_{s \in [0,1]}$ is the Brownian motion. (The fact that $W$ satisfies condition \eqref{cond-Y} is proved in Appendix \ref{section-appC}.) We simulate $K$ i.i.d. copies of $W$ using Donsker theorem. Let $\xi,(\xi_{jk})_{j,k \geq 1}$ be i.i.d. random variables with mean 0 and variance 1.
When $m$ is large, $W_j(s) \approx \frac{1}{\sqrt{m}} \sum_{k=1}^{[ms]}\xi_{jk}$ for any $j=1,\ldots,K$, and $X(s) \approx \sum_{j=1}^{K} \e_j \Gamma_j^{-1/\alpha} W_j(s) \approx \frac{1}{\sqrt{m}} \sum_{j=1}^{K}\sum_{k=1}^{[ms]} \e_j \Gamma_j^{-1/\alpha}  \xi_{jk}$ for any $s \in [0,1]$.

Next, we consider $n$ i.i.d. copies of $X$. Let $(\e_{ij})_{i,j\geq 1}$ be i.i.d. copies of $\e_1$, $(E_{ij})_{i,j \geq 1}$ i.i.d. copies of $E_1$ and $(\xi_{ijk})_{i,j,k\geq 1}$ i.i.d. copies of
$\xi$. Let $\Gamma_{ij}=\sum_{k=1}^{j}E_{ik}$. We take $a_n=(nC_{W,\alpha})^{1/\alpha}$ where $C_{W,\alpha}$ is computed by approximation. By Theorem \ref{main2},
$$Z(t,s) \approx \frac{1}{a_n}\sum_{i=1}^{[nt]}X_i(s) \approx \frac{1}{a_n \sqrt{m}} \sum_{i=1}^{[nt]} \sum_{k=1}^{[ms]} \sum_{j=1}^{K}  \e_{ij} \Gamma_{ij}^{-1/\alpha}  \xi_{ijk}$$
is an approximation of a $\bD$-valued $\alpha$-stable L\'evy motion, when $n$, $m$ and $K$ are large.

The following pictures are the 3-dimensional plots of
$(t_k,s_l,Z(t_k,s_l))$ for $k=1,\ldots,n$ and $l=1,\ldots,m$, with $t_k=k/n$ and $s_l=l/m$, when $n=400$ and $m=250$.

\begin{figure}[h!]
    \centering
    \begin{subfigure}[b]{0.48\textwidth}
        \includegraphics[width=\textwidth]{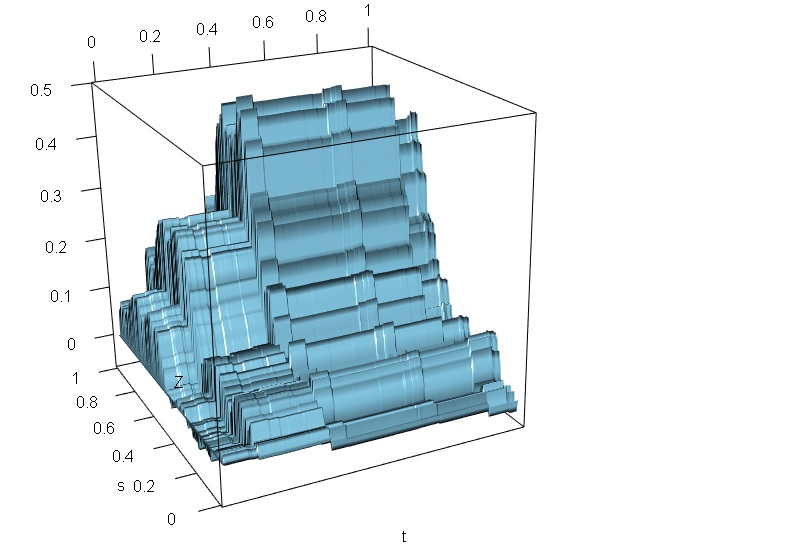}
        \caption{$\alpha=0.5$}
        \label{fig:series1}
    \end{subfigure}
    ~ 
    \begin{subfigure}[b]{0.48\textwidth}
        \includegraphics[width=\textwidth]{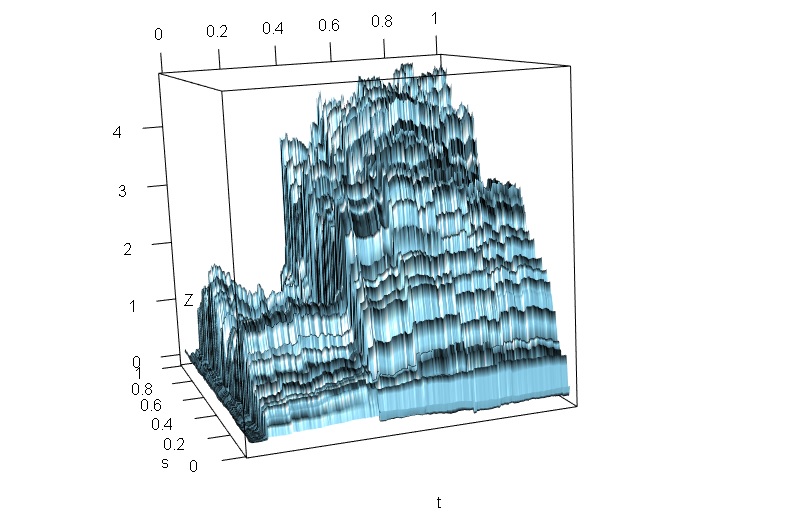}
        \caption{$\alpha=1.5$}
        \label{fig:series2}
    \end{subfigure}
    ~ 
    \caption{$\bD$-valued $\alpha$-stable L\'evy motion based on a regularly process in $\bD$ given by  series \eqref{series} in which $(Y_j)_{j\geq 1}$ are i.i.d. Brownian motions}
    \label{fig:series}
\end{figure}

}
\end{example}

\appendix

\section{Some auxiliary results}
\label{appendixA}

In this section, we include some auxiliary results which are used in this article.

The first result shows that the measure $\overline{\nu}$ which appears in the definition of regularly variation for random elements in $\bD$ must be of product form. This result is probably well-known. We include its proof since we could not find it in the literature.

\begin{lemma}
\label{nu-product}
If $c=\overline{\nu}((1,\infty)\times \SD)>0$, then the measure $\overline{\nu}$ in Definition \ref{def-RV} must be of the product from \eqref{def-bar-nu}, with probability measure $\Gamma_1$ given by \eqref{def-Gamma1}.
\end{lemma}

\noindent{\bf Proof:} Let ${\cal P}$ be the class of sets $A_{r,S}=(r,\infty) \times S$ with $r>0$ and $S \in \cB(\SD)$. Note that $A_{r,s}=T(V_{r,s})$ where $V_{r,S}=\{x \in  \bD; \|x\|>r, \frac{x}{\|x\|} \in S\}$. The sets $V_{r,S}$ have the scaling property $aV_{r,S}=V_{ar,S}$ for any $a>0$. To see that the sets $A_{r,S}$ have a similar property, we define $aA=\{(ar,z); (r,z) \in A\}$
for any $a>0$ and $A \in \cB(\barD)$. Then
 $$aA_{r,S}=\{(as,z); as>ar,z \in \SD\}=A_{ar,S}.$$
In particular, $A_{r,S}=rA_{1,S}$. By the scaling property of $\overline{\nu}$ and definition \eqref{def-Gamma1} of $\Gamma_1$,
$$\overline{\nu}(A_{r,S})=r^{-\alpha}\overline{\nu}(A_{1,S})=r^{-\alpha}c\Gamma_1(S)=(c \nu_{\alpha} \times \Gamma_1)(A_{r,S}).$$
Hence, when restricted to $(0,\infty) \times \SD$, the measures $\overline{\nu}$ and $c \nu_{\alpha} \times \Gamma_1$ coincide for sets in the class ${\cal P}$. Since ${\cal P}$ is a $\pi$-system which generates the Borel $\sigma$-field of $(0,\infty) \times \SD$ (with respect to distance $d_{\barD}$), it follows that $\overline{\nu}=c \nu_{\alpha} \times \Gamma_1$ on $(0,\infty) \times \SD$. Finally, these measures coincide on the entire space $\barD$ since they are zero on $\{\infty\} \times \SD$. $\Box$

\vspace{3mm}

The next result is an extension of Lemma 5.2 of \cite{resnick07} to the case of functions with values in an arbitrary metric space.

\begin{lemma}
\label{lemmaA}
Let $(S,d)$ be a complete metric space. We denote by $\bD([0,\infty);S)$ the set of functions $x:[0,\infty) \to S$ which are right-continuous and have left-limits (with respect to $d$). If $(x_n)_{n\geq 1}$ is a sequence in $\bD([0,\infty);S)$ and the function $x:[0,\infty) \to S$ is such that
\begin{equation}
\label{unif-conv}
\sup_{t\leq T}d(x_n(t),x(t)) \to 0 \quad \mbox{as} \quad n \to \infty,
\end{equation}
for any $T>0$, then $x \in \bD([0,\infty);S)$.
\end{lemma}

\noindent {\bf Proof:} We first prove that $x$ is right-continuous. Let $t\geq 0$ be arbitrary and $(t_k)_{k \geq 1}$ be such that $t_k \to t$ and $t_k \geq t$ for all $k$. Let $T>0$ be such that $t_k \in [0,T]$ for all $k$. Let $\e>0$ be arbitrary. By \eqref{unif-conv}, there exists $n_0$ such that
$d(x_{n_0}(t),x(t))<\e$ for all $t \leq T$. Since $x_{n_0}$ is right-continuous at $t$, there exists $K_{\e}$ such that $d(x_{n_0}(t_k),x_{n_0}(t))<\e$ for all $k\geq K_{\e}$. By the triangle inequality, for any $k\geq K_{\e}$,
$$d(x(t_k),x(t)) \leq d(x(t_k),x_{n_0}(t_k)) +d(x_{n_0}(t_k),x_{n_0}(t))+d(x_{n_0}(t),x(t))<3\e.$$

Next, we prove that $x$ has left limit at $t>0$. Let $(t_k)_{k \geq 1}$ be such that $t_k \to t$ and $t_k < t$ for all $k$. Let $T>0$ be such that $t_k \in [0,T]$ for all $k$. Let $\e>0$ be arbitrary. Choose $n_0$ as above. Since $x_{n_0}(t_{k}) \to x(t-)$, $\{x_{n_0}(t_k)\}_{k}$ is a Cauchy sequence. Then, there exists $L_{\e}$ such that $d(x_{n_0}(t_k),x_{n_0}(t_l))<\e$ for all $k,l\geq L_{\e}$. Then
$$d(x(t_k),x(t_l)) \leq d(x(t_k),x_{n_0}(t_k)) +d(x_{n_0}(t_k),x_{n_0}(t_l))+d(x_{n_0}(t_l),x(t_l))<3\e,$$
for any $k,l\geq L_{\e}$, and hence $\{x(t_k)\}_{k}$ is a Cauchy sequence. Since $S$ is complete, there exists $l=\lim_{k \to \infty}x(t_k)$. We must show that $l$ does not depend on $(t_k)_k$. Let
$l'=\lim_{k \to \infty}x(t_k')$, where
$(t_k')_k$ is another sequence such that $t_k' \to t$ and $t_k'< t$ for all $k$. Since both sequences $\{x_{n_0}(t_k)\}_{k}$ and $\{x_{n_0}(t_k')\}_{k}$ converge to $x_{n_0}(t-)$, there exists $M_{\e}$ such that
$d(x_{n_0}(t_k),x_{n_0}(t_k'))<\e$ for all $k\geq M_{\e}$. Hence,
$$d(x(t_k),x(t_k')) \leq d(x(t_k),x_{n_0}(t_k)) +d(x_{n_0}(t_k),x_{n_0}(t_k'))+d(x_{n_0}(t_k'),x(t_k'))<3\e,$$
for any $k\geq M_{\e}$. This proves that $l=l'$. $\Box$

\vspace{3mm}

The following result is probably well-known. We include its proof since we could not find it in the literature.

\begin{lemma}
\label{lemmaB}
Let $(S,d)$ be a separable metric space. Let $X_n^{(1)}, \ldots, X_n^{(k)}$ and $X^{(1)},\ldots, X^{(k)}$ be random elements in $S$ defined on a probability space $(\Omega,\cF,P)$, such that $d(X_n^{(i)},X^{(i)}) \to 0$ a.s. for any $i=1,\ldots,k$. If $X_n^{(1)},\ldots,X_n^{(k)}$ are independent for any $n \geq 1$, then $X^{(1)},\ldots,X^{(k)}$ are independent.
\end{lemma}

\noindent {\bf Proof:} We assume for simplicity that $k=2$, the general case being similar. To simplify the notation, we let $X_n=X_n^{(1)}$ and $Y_n=X_n^{(2)}$. Clearly, $d(X_n,X) \stackrel{P}{\to} 0$  and $d(Y_n,Y) \stackrel{P}{\to} 0$. Note that the space $S \times S$ equipped with the product metric is separable and $(X_n,X)$ is a random element in $S \times S$ (see p.225 of \cite{billingsley68}).
By Corollary to Theorem 3.1 of \cite{billingsley99}, $X_n \stackrel{d}{\to} X$ and $Y_n \stackrel{d}{\to} Y$. By Theorem 3.2 of \cite{billingsley68},
\begin{equation}
\label{w-conv1}
(P\circ X_n^{-1}) \times (P \circ Y_n^{-1}) \stackrel{w}{\to}(P \circ X^{-1}) \times (P \circ Y^{-1}) \quad \mbox{on}  \quad S \times S.
\end{equation}
On the other hand, $(X_n,Y_n) \to (X,Y)$ a.s. with respect to the product distance in $S \times S$. Hence, again by Corollary to Theorem 3.1 of \cite{billingsley99}, $(X_n,Y_n) \stackrel{d}{\to} (X,Y)$ in $S \times S$, i.e.
\begin{equation}
\label{w-conv2}
P\circ (X_n,Y_n)^{-1} \stackrel{w}{\to}P \circ (X,Y)^{-1} \quad \mbox{on}  \quad S \times S.
\end{equation}
Finally, $P\circ (X_n,Y_n)^{-1}=(P\circ X_n^{-1}) \times (P \circ Y_n^{-1})$ for any $n\geq 1$, since $X_n$ and $Y_n$ are independent for any $n \geq 1$. The fact that $P \circ (X,Y)^{-1}=(P \circ X^{-1}) \times (P \circ Y^{-1})$ follows from \eqref{w-conv1} and \eqref{w-conv2}, by the uniqueness of the limit. $\Box$

\section{The $\alpha$-stable L\'evy sheet}

In this section, we show that the $\alpha$-stable L\'evy sheet can be viewed as an example of a $\bD$-valued $\alpha$-stable L\'evy motion restricted to the time interval $[0,1]$.

First, we recall briefly the construction of the $\alpha$-stable L\'evy sheet, as described in Section 4.8 of \cite{resnick86}. Let $M=\sum_{i\geq 1}\delta_{(T_i,S_i,J_i)}$ be a Poisson random measure on $[0,\infty) \times [0,\infty) \times \overline{\bR}_0$ of intensity ${\rm Leb} \times {\rm Leb} \times \nu_{\alpha,p}$, where
$\nu_{\alpha,p}$ is given by \eqref{def-nu-ap},
for some $\alpha \in (0,2),\alpha \not=1$ and $p \in [0,1]$, with $q=1-p$.
Let $(\e_j)_{j\geq 0}$ be a sequence of real numbers such that $\e_j \downarrow 0$ and $\e_0=1$. Let $I_j=(\e_j,\e_{j-1}]$ for $j\geq 1$ and $I_0=(1,\infty)$. For any $t,s\in [0,1]$ and $j\geq 0$, let
$$L_j(t,s)=\int_{[0,t] \times [0,s] \times \Gamma_j}zM(du_1,du_2,dz)=\sum_{i\geq 1}J_i 1_{\{J_i \in \Gamma_j\}}1_{\{T_i \leq t, S_i \leq s\}}.$$
Note that $L_j(t,s)$ is a compound Poisson random variable with characteristic function
$$E[e^{iu L_j(t,s)}]=\exp \left\{ts \int_{\Gamma_j} (e^{iuz} -1) \nu_{\alpha,p}(dz) \right\}, \quad u \in \bR.$$
By Kolmogorov's criterion, the series
$\sum_{j\geq 1}\big(L_j(t,s)-E(L_{j}(t,s))\big)$ converges a.s., since  ${\rm Var}\big(L_j(t,s) \big)=ts \int_{\Gamma_j}z^2 \nu_{\alpha,p}(dz)$ for any $j\geq 1$ and $\int_{|z| \leq 1}z^2 \nu_{\alpha,p}(dz)<\infty$.

We define $\overline{L}(t,s)=\sum_{j\geq 0}L_j(t,s)$ if $\alpha<1$ and $\overline{L}(t,s)=\sum_{j\geq 0}\big( L_j(t,s)-E(L_j(t,s))\big)$ if $\alpha>1$. It can be proved that there exists a process $\{L(t,s)\}_{(t,s) \in [0,1]^2}$ with sample paths in $\bD([0,1]^2)$ such that $L(t,s)=\overline{L}(t,s)$ a.s. for any $t,s \in [0,1]$, and
\begin{align}
\label{sup-conv}
& \sup_{(t,s)\in [0,1]^2}|L^{(\e_k)}(t,s)-L(t,s)| \to 0 \quad \mbox{a.s.} \quad \mbox{if} \quad \alpha<1, \\
& \sup_{(t,s)\in [0,1]^2}|\overline{L}^{(\e_k)}(t,s)-L(t,s)| \to 0 \quad \mbox{a.s.} \quad \mbox{if} \quad \alpha>1,
\end{align}
where $L^{(\e_k)}(t,s)=\sum_{j=0}^k L_j(t,s)$ and $\overline{L}^{(\e_k)}(t,s)=L^{(\e_k)}(t,s)-E(L^{(\e_k)}(t,s))$ (if $\alpha>1$). Here $\bD([0,1]^2)$ is the space of functions $x:[0,1]^2 \to \bR$  which are continuous at any point $(t,s)$ when this point is approached from the upper right quadrant, and have limits when the point is approached from the other three quadrants. Moreover,
$$E[e^{iu L(t,s)}]=\exp\left\{ts \int_{\bR}(e^{iuz}-1)\nu_{\alpha,p}(dz) \right\} \quad \mbox{if} \quad \alpha<1, $$
$$E[e^{iu L(t,s)}]=\exp\left\{ts \int_{\bR}(e^{iuz}-1-iuz)\nu_{\alpha,p}(dz) \right\} \quad \mbox{if} \quad \alpha>1.$$
Consequently, $L(t,s)$ has a $S_{\alpha}\big((ts)^{1/\alpha}C_{\alpha}^{-1},\beta,0\big)$-distribution with $\beta=p-q$ and $C_{\alpha}$ given by \eqref{def-C}.
The process $\{L(t,s)\}_{(t,s) \in [0,1]^2}$ is called an {\em $\alpha$-stable L\'evy sheet}. Note that both processes $\{L(t,s)\}_{t \in [0,1]}$ and $\{L(t,s)\}_{s \in [0,1]}$ are $\alpha$-stable L\'evy motions with paths in $\bD$.

\begin{theorem}
\label{stable-sheet-th1}
Let $L(t)=\{L(t,s)\}_{s \in [0,1]}$ for any $t \in [0,1]$. The process $\{L(t)\}_{t \in [0,1]}$ is an $\bD$-valued $\alpha$-stable L\'evy motion (according to Definition \ref{def-Levy}).
\end{theorem}

\noindent {\bf Proof:} We show that $\{L(t)\}_{t \in [0,1]}$ satisfies conditions {\em (i)}-{\em (iv)} of Definition \ref{def-Levy}. We assume that $\alpha<1$, the case $\alpha>1$ being similar. Clearly $L(0)=0$, so property {\em (i)} holds.

For property {\em (ii)}, note that by \eqref{sup-conv}, $L^{(\e_k)}(t_i) \to L(t_i)$ a.s. in $(\bD,\|\cdot\|)$ as $k \to \infty$ for $i=1,\ldots,K$, and hence $L^{(\e_k)}(t_i)-L^{(\e_k)}(t_{i-1}) \to L(t_i)-L(t_{i-1})$ a.s. in $(\bD,\|\cdot\|)$ as $k \to \infty$, for any $i=2,\ldots,K$. By Lemma \ref{lemmaB}, $L(t_i)-L(t_{i-1}),i=2,\ldots,K$ are independent, since $L^{(\e_k)}(t_i)-L^{(\e_k)}(t_{i-1}),i=2,\ldots,K$ are independent for any $k$.

To verify property {\em (iii)}, we observe that for any $t_1<t_2$ and $s \in [0,1]$,
$$L(t_2,s)-L(t_1,s)=\overline{L}(t_1,s)-\overline{L}(t_2,s)=\sum_{j\geq 0} \int_{(t_1,t_2] \times [0,s] \times \Gamma_j} z M(du_1,du_2,dz) \quad \mbox{a.s.}$$
From this, it can be proved that $L(t_2)-L(t_1)=\{L(t_2,s)-L(t_1,s)\}_{s \in [0,1]}$ is an $\alpha$-stable L\'evy motion with characteristic function
$$E[e^{iu (L(t_2,s)-L(t_1,s))}]=\exp \left\{ (t_2-t_1)s \int_{\bR} (e^{iuz}-1)\nu_{\alpha,p}(dz)\right\}, \quad \mbox{for all} \ u \in \bR.$$
On the other hand, $L(t_2-t_1)=\{L(t_2-t_1,s)\}_{s \in [0,1]}$ is also an $\alpha$-stable L\'evy motion with the same characteristic function. Hence, $L(t_2)-L(t_1) \stackrel{d}{=}L(t_2-t_1)$.

To verify property {\em (iv)}, we assume first that $t=1$. The process $L(1)=\{L(1,s)\}_{s \in [0,1]}$ is an $\alpha$-stable L\'evy motion, so it is an $\alpha$-stable process. It follows that for any $s_1,\ldots,s_m \in [0,1]$,
$(L(1,s_1),\ldots,L(1,s_m))$ has an $\alpha$-stable distribution in $\bR^m$ with L\'evy measure $\mu_{s_1,\ldots,s_m}$:
$$E\big(e^{iu_1 L(1,s_1)+\ldots+iu_m L(1,s_m)} \big)=\exp \left\{ \int_{\bR^m}(e^{iu \cdot y}-1)\mu_{s_1 \ldots,s_m}(dy)\right\}, \quad u=(u_1,\ldots,u_m)\in \bR^m.$$
In particular, $(L(1,s_1),\ldots,L(1,s_m))$ is regularly varying with limiting measure $\mu_{s_1,\ldots,s_m}$.

On the other hand, by Lemma 2.1 of \cite{hult-lindskog07}, $L(1)$ is regularly varying in $\bD$ (in the sense of Definition \ref{def-RV}), i.e. $L(1) \in RV(\{a_n\},\overline{\nu},\overline{\bD}_0)$ for a boundedly finite measure $\overline{\nu}$ on $\overline{\bD}_0$ with $\overline{\nu} (\overline{\bD}_0 \verb2\2 T(\bD_0))=0$. Moreover, $\overline{\nu}=c\nu_{\alpha} \times \Gamma_1$ for some $c>0$ and a probability measure $\Gamma_1$ on $\bS_{\bD}$. Let $\nu=\overline{\nu} \circ S^{-1}$, where $S:(0,\infty) \times \bS_{\bD} \to \bD_0$ is the inverse of the map $T$, i.e. $S(r,z)=rz$. By Theorem 8 of \cite{hult-lindskog05},
$(L(1,s_1),\ldots,L(1,s_m))$ is regularly varying with limiting measure $\nu_{s_1,\ldots,s_m}=\nu \circ \pi_{s_1,\ldots,s_m}^{-1}$. By the unicity of the limit, $\mu_{s_1,\ldots,s_m}=\nu_{s_1,\ldots,s_m}$.
Finally, property {\em (iv)} for general $t$ follows using the scaling property of $\mu_{s_1,\ldots,s_m}$ and the fact that $\{L(t,s)\}_{s \in [0,1]} \stackrel{d}{=}\{t^{1/\alpha}L(1,s)\}_{s \in [0,1]}$. $\Box$

\vspace{3mm}

In relation with the simulation procedure described in Example \ref{ex-Levy-sheet}, we include the following result, which can be proved using the same argument as in Section 48 of \cite{resnick86}.

\begin{theorem}
\label{stable-sheet-th2}
Let $\xi,(\xi_{ij})_{i,j \geq 1}$ be i.i.d. regularly varying random variables, i.e. $$n P\left(\frac{\xi}{a_n} \in \cdot \right) \stackrel{v}{\to} \nu_{\alpha,p} \quad \mbox{in} \quad \overline{\bR}_0,$$
for some $a_n \uparrow \infty$, $\alpha \in (0,2),\alpha\not=1$ and $p \in [0,1]$, where $\nu_{\alpha,p}$ is given by \eqref{def-nu-ap}.  For any $t,s \in [0,1]$, let
$T_{n,m}(t,s)=a_n^{-1}a_m^{-1} \sum_{i=1}^{[nt]} \sum_{j=1}^{[ms]} (\xi_{ij}-\mu)$, where $\mu=0$ if $\alpha<1$ and $\mu=E(\xi)$ if $\alpha>1$. Let $L=\{L(t,s)\}_{(t,s) \in [0,1]^2}$. Then
$$T_{n,m} \stackrel{d}{\to}L \quad in \quad \bD([0,1]^2), \quad as \quad n,m\to \infty.$$
\end{theorem}

\section{A result about Brownian motion}
\label{section-appC}

In this section, we include a result about the Brownian motion which is used in Example \ref{ex-series}. This result is probably well-known. We include its proof since we could not find it in the literature.

\begin{lemma}
\label{BM-th}
Let $W=\{W(s)\}_{s \in [0,1]}$ be the Brownian motion. Then,
$$E\Big(\sup_{s \in [0,1]}|W(s)|^{\alpha}\Big)<\infty \quad for \ all \ \alpha >0.$$
\end{lemma}

\noindent {\bf Proof:} Let $W^+(t)=\max(W(t),0)$ and $W^{-}(t)=\max(-W(t),0)$.
For any $x>0$,
\begin{align*}
\{\sup_{t \in [0,1]} |W(t)|>x\} & \subset \{\sup_{t \in [0,1]} W^+(t)>x/2\} \cup \{\sup_{t \in [0,1]} W^{-}(t)>x/2\}\\
&=\{\sup_{t \in [0,1]} W(t)>x/2\} \cup \{\sup_{t \in [0,1]} (-W(t))>x/2\}.
\end{align*}
Note that $\{-W(t)\}_{t \in [0,1]} \stackrel{d}{=} \{W(t)\}_{t \in [0,1]}$. By reflection principle for the Brownian motion,
$$P(\sup_{t \in [0,1]} |W(t)|>x) \leq 2 P(\sup_{t \in [0,1]} W(t)>x/2)=4P(W(1)>x/2) \leq 4P(|W(1)|>x/2).$$
Hence,
\begin{align*}
E\Big(\sup_{t \in [0,1]} |W(t)|^{\alpha}\Big)& =\int_0^{\infty}P(\sup_{t \in [0,1]} |W(t)|>x^{1/\alpha})dx \leq
4 \int_0^{\infty}P(|W(1)|>(x/2)^{1/\alpha})dx\\
&=8E|W(1)|^{\alpha}<\infty.
\end{align*}
$\Box$

\vspace{3mm}

{\bf Acknowledgement.} We would like to thank Fran\c{c}ois Roueff, Gennady Samorodnitsky and Philippe Soulier for useful discussions, and for drawing their attention to reference \cite{davydov-molchanov-zuyev08} regarding $\alpha$-stable L\'evy processes on cones (see Remark \ref{remark-cone}). We are also grateful to Thomas Mikosch for the proof of Lemma \ref{BM-th}, to Xiao Liang for his help with the simulations, and to Adam Jakubowski for reading the manuscript.

\end{document}